\documentclass[reqno]{amsart}

\usepackage{amsmath}
\usepackage{amscd}
\usepackage{amssymb}
\usepackage{amstext}
\usepackage{amsmath}
\usepackage[all]{xy}
\usepackage{yhmath}
\usepackage{mathrsfs}
\usepackage{bbding}
\usepackage{graphicx}

\def\E{\ifmmode{\mathbb E}\else{$\mathbb E$}\fi} 
\def\N{\ifmmode{\mathbb N}\else{$\mathbb N$}\fi} 
\def\R{\ifmmode{\mathbb R}\else{$\mathbb R$}\fi} 
\def\Q{\ifmmode{\mathbb Q}\else{$\mathbb Q$}\fi} 
\def\C{\ifmmode{\mathbb C}\else{$\mathbb C$}\fi} 
\def\H{\ifmmode{\mathbb H}\else{$\mathbb H$}\fi} 
\def\Z{\ifmmode{\mathbb Z}\else{$\mathbb Z$}\fi} 
\def\P{\ifmmode{\mathbb P}\else{$\mathbb P$}\fi} 
\def\T{\ifmmode{\mathbb T}\else{$\mathbb T$}\fi} 
\def\SS{\ifmmode{\mathbb S}\else{$\mathbb S$}\fi} 
\def\DD{\ifmmode{\mathbb D}\else{$\mathbb D$}\fi} 
\def\K{\ifmmode{\mathbb K}\else{$\mathbb K$}\fi}

\theoremstyle{theorem}
\newtheorem{thm}{Theorem}[section]
\newtheorem{cor}[thm]{Corollary}

\newtheorem{prop}[thm]{Proposition}

\newtheorem{conj}[thm]{Conjecture}

\theoremstyle{definition}
\newtheorem{defn}[thm]{Definition}
\newtheorem{rem}[thm]{Remark}

\numberwithin{equation}{section}

\title[Homological algebra and moduli spaces]{Homological algebra and moduli spaces in
topological field theories}
\author{\bf Kenji Fukaya}
\begin{document}
\maketitle
\tableofcontents

\section{Numerical invariants: Donaldson and Gromov.}

The relations between topology and linear partial differential 
equations are classical going back to 19 century and 
a recent high point is Atiyah-Singer index theorem.
The relations between non-linear differential equations 
and geometry are very hot and active topics which become 
big trends after several important discoveries in the 1970's.
More recently relations between moduli spaces of the solutions of non-linear differential equations 
and `topology' are discovered and become an important topic.
\par
In his famous paper \cite{Don1},  Donaldson used the moduli space 
of the solutions of  Yang-Mills equation to obtain novel restrictions on the 
intersection forms of closed 4-dimensional manifolds.
In subsequent papers \cite{Don2,Don3}, 
the moduli space is used to obtain an invariant of a closed smooth 
4-dimensional manifold.
\par
Let $X$ be a closed $4$-dimensional manifold and $\mathcal E_X \to X$  an 
$SU(2)$ bundle.\footnote{One can use $SO(3)$ bundle also.}
A connection $A$ of $\mathcal E_X$ is said to be an ASD (anti-self-dual) connection, if
it satisfies the equation:
\begin{equation}\label{ASDeq}
F_A + *_X F_A = 0.
\end{equation}
Here $F_A$ is the curvature of $A$ and $*_X$ is the Hodge $*$ operator.
Note that $*_X$ depends on the choice of a Riemannian metric on $X$.
The set of solutions of (\ref{ASDeq}) 
is invariant under the action of the gauge transformation group
(the set of automorphisms of the principal bundle $\mathcal E_X \to X$).
We denote by $\mathcal M(X;\mathcal E_X)$ the set of gauge equivalence 
classes of the solutions of (\ref{ASDeq}).
The space $\mathcal M(X;\mathcal E_X)$ has the following nice properties.
\begin{enumerate}
\item[(ASD1)]
The moduli space $\mathcal M(X;\mathcal E_X)$ is generically a smooth 
manifold outside the set of points corresponding to the reducible connections.\footnote{that is, the 
connection $A$ such that the set of gauge transformations which preserves 
$A$ has positive dimension.} 

\item[(ASD2)]
The space  $\mathcal M(X;\mathcal E_X)$ has a nice compactification 
called the Uhlenbeck compactification.
\item[(ASD3)]
The `cobordism class' of $\mathcal M(X;\mathcal E_X)$ 
is independent of  various choices.
\end{enumerate}

In the case of $SU(2)$ bundles, a reducible connection 
is induced from a $U(1)$ connection. The curvature of an anti-self-dual 
$U(1)$ connection is a harmonic two form $h$ with $*_{X}h = -h$.
So the negative eigen-space on the second cohomology group $H^2(X;\R)$
(with respect to the intersection form) is related to the 
singularity of  $\mathcal M(X;\mathcal E_X)$. 
We denote by $b_2^+$ and $b_2^-$ the dimensions of positive and negative 
eigen-spaces on the second cohomology group $H^2(X;\R)$
with respect to the intersection form.
If $b_2^+ = 0$ (and $b_2^- \ne 0$) then for any Riemannian metric 
on $X$ there exists a 
harmonic 2 form $h$ representing integral cohomology classes 
and satisfying $*_{X}h = -h$.
It implies that there exists a reducible connection 
for any Riemannian metric. In general, the set of 
Riemannian metrics for which 
$\mathcal M(X;\mathcal E_X)$ contains a reducible connection has codimension  $b_2^+$. Thus 
$\mathcal M(X;\mathcal E_X)$ is `closer' to a non-singular space 
(manifold)  
if  $b_2^+$ is large. 
In fact the Donaldson invariant is well-defined 
if $b_2^+ \ge 2$ but is not defined 
if $b_2^+ = 0$.
The case $b_2^+=1$ is a borderline case where 
the invariant exists but depends on the `chamber' in which the Riemannian 
metric is contained.
Namely it depends on the metric however we can control 
the way how it changes when a family of Riemannian metrics crosses the `wall'
and moves from one chamber to the other.
The case  $b_2^+=1$ is used in \cite{Don2} to obtain 
the first example of a pair of closed $4$-dimensional manifolds 
which are homeomorphic but are not diffeomorphic.
\par
In the simplest case, that is, when the (virtual) dimension $\mathcal M(X;\mathcal E_X)$
is zero, Donaldson invariant is a number.
In general Donaldson introduced and used a certain 
cohomology class to cut down the space 
$\mathcal M(X;\mathcal E_X)$ and obtain a number.
He used a map
$$
\nu: H_*(X;\Z) \to H^{4-*}(\mathcal B^*(X;\mathcal E_X);\Z).
$$
Here $\mathcal B^*(X;\mathcal E_X)$ is the set of 
gauge equivalence classes of the irreducible connections of $\mathcal E_X$.
Donaldson invariant can be regarded as a polynomial 
on $H_*(X;\Z)$ and written as
\begin{equation}\label{forDona}
Q(a_1,\dots,a_k) = \int_{\mathcal M(X;\mathcal E_X)} \nu(a_1) \wedge \dots \nu(a_k).
\end{equation}
(ASD3) `implies' that this number is independent of the Riemannian metric 
etc. and becomes an invariant of a smooth 4-dimensional manifold.

A certain delicate  (dimension counting type) argument is necessary 
to understand how the reducible connections and  
the infinity of the Uhlenbeck compactification
affect the well-defined-ness of the integral (\ref{forDona}). 
\par\medskip
Gromov \cite{Gr} introduced the method of pseudo-holomorphic curve 
to symplectic geometry.
For a symplectic manifold $(X,\omega)$ Gromov considered 
a compatible almost complex structure $J$, that is, 
a tensor $J: TX \to TX$ such that $J^2 = -1$ and $g(V,W):= \omega(V,JW)$
becomes a Riemannian metric.
Then a pseudo-holomorphic curve is a map $u$ from a Rieman surface $(\Sigma,j_{\Sigma})$
to $(X,J)$ such that
\begin{equation}\label{phceq}
J \circ Du = Du \circ j_{\Sigma}.
\end{equation}
The equation (\ref{phceq}) is also written as 
\begin{equation}\label{phceq2}
\overline{\partial}u = 0.
\end{equation}
Here it is very important that $J$ is an {\it almost} complex structure 
which is not necessary integrable. 
If $J$ is an integrable complex structure we can take 
a complex coordinate of $X$ so that $J = \sqrt{-1}$ is constant.
Then writing $u = (u^1,\dots,u^n)$ by the complex coordinate,
(\ref{phceq}) becomes
\begin{equation}\label{linCR}
\frac{\partial u^i}{\partial y} = \sqrt{-1}\frac{\partial u^i}{\partial x},
\end{equation}
where $z = x + \sqrt{-1}y$ is a complex coordinate of the Rieman surface 
$(\Sigma,j_{\Sigma})$.
The equation (\ref{linCR}) is a linear partial differential 
equation. On the other hand, (\ref{phceq}) is non-linear.
An important observation by Gromov is that the non-linear 
partial differential equation (\ref{phceq}) can have lots of solutions 
in the case when the domain has complex dimension one.
In fact if $(Y,j)$ is a complex two dimensional manifold,
typically there is no non-constant map 
$u : Y \to X$ satisfying $J \circ Du = Du \circ j$ 
in the case when $J$ is not integrable.
\par
Another important point is in the case when the almost 
complex structure $J$ is compatible with a certain 
symplectic structure\footnote{Actually a slightly weaker 
condition, that $J$ is tamed by a symplectic structure, is enough.}
the moduli space of solutions of the equation (\ref{phceq}) 
has a nice compactification.
In fact we have the following properties.
We fix a non-negative integer $g$ and a positive number $E$
and for a symplectic manifold $X$ with compatible almost 
complex structure $J$, we denote by 
$\overset{\circ}{\mathcal M}_g(X,J;E)$ the set of 
pairs $((\Sigma,j),u)$ where 
$(\Sigma,j)$ is a Rieman surface of genus $g$
and $u: \Sigma \to X$ satisfies the equation (\ref{phceq}).
We also require
$$
\int_{\Sigma} u^* \omega \le E.
$$

\begin{enumerate}
\item[(PHC1)]
The moduli space $\overset{\circ}{\mathcal M}_g(X,J;E)$ is generically a smooth 
manifold.
\item[(PHC2)]
The space  $\overset{\circ}{\mathcal M}_g(X,J;E)$  has a nice compactification.\footnote{
Gromov (and also \cite{McSa94} etc.) used a compactification which he called the 
moduli space of cusp curves. This compactification works for the purpose of Gromov's paper \cite{Gr} and also 
in the semi-positive case but not likely works in the general case.
Later Kontsevich introduced a {\it different} compactification,  
the stable map compactification (whose origin is in algebraic geometry), 
which is now widely used in symplectic geometry also.
The stable map topology on ${\mathcal M}_g(X,J;E)$ is defined 
in \cite[Definition 10.3]{FO}} 
\item[(PHC3)]
The `cobordism class' of the compactification ${\mathcal M}_g(X,J;E)$ 
is independent of the various choices especially 
of the choice of the almost complex structure.
\end{enumerate}

These properties are similar to the properties 
(ASD1), (ASD2), (ASD3).
Gromov used them to prove that, for any almost 
complex structure $J$ on $\C P^2$ which is compatible 
with the standard symplectic structure, 
there exists a pseudo-holomorphic curve $u: S^2 \to (\C P^2,J)$
whose homology class is the generator of $H_2(\C P^2)$.
This fact has an important application which is 
called non-sqeezing theorem.
\begin{thm}{\rm (\cite{Gr})}
If there exists a map $u$ from $\overset{\circ}{D^4}(r)$ (the open $r$-ball in $\C^2$) 
to $D^2(1) \times \C$ such that $u^*\omega = \omega$
(where $\omega$ is the standard symplectic form) then 
$r<1$.
\end{thm}
It implies, for  example, that $\overset{\circ}{D^4}(1)$ is not symplectomorphic 
to $\overset{\circ}{D^2}(2) \times \overset{\circ}{D^2}(1/2)$.
This is one of the first results which show 
the `existence of global symplectic geometry'.
\par
Around the same time, physicists working on string theory 
studied `topological' version of string theory 
and the `invariant' obtained by integrating 
a certain cohomology class on the moduli space 
which can be regarded as a compactification of $\overset{\circ}{\mathcal M}_g(X,J;E)$.
(See for example \cite{witten2}.)
One point which is not so clear from physicists'
point of view is the fact that 
such `invariant' is one of symplectic structure
and not one of almost complex structure (or complex structure).
On the other hand, various important properties of 
the invariant (obtained from the moduli space of 
pseudo-holomorphic curves) are discovered by physicists.
Among them the associativity of the product structure\footnote{the product 
structure of the quantum cohomology ring} 
obtained from a pseudo-holomorphic map $S^2 \to X$ 
is  very important (see \cite{VA}).
\par
Ruan \cite{Ru1} and Ruan-Tian \cite{RT1,RT2} 
established the theory of invariants obtained 
from the moduli space of 
pseudo-holomorphic curves.\footnote{They assumed a certain 
positivity assumption, which was removed later in the year 1996 
by groups of mathematicians.}
After McDuff-Salamon's lucid exposition  \cite{McSa94} appeared
this theory becomes popular among differential and symplectic geometers.

\section{Floer homology.}
\label{sec:Floer}

Floer homology was discovered in the 1980's (by A. Floer) 
in two areas. One is gauge theory and the other is 
symplectic geometry.
Floer's work is a development of three important works 
in those areas.
\begin{enumerate}
\item
Casson invariant of 3-dimensional manifolds (See \cite{AM})  and Taubes' work \cite{Tau} 
to relate it to gauge theory.
\item
Conley-Zehnder's proof \cite{CZ} of Arnold's conjecture 
for tori. (There was a related work \cite{Rb} before that.)
\item
Witten's work \cite{Wi} which relates Morse theory to a 
(supersymmetric) quantum field theory.
\end{enumerate}

All of these are famous and important works. We mention them briefly.
\par
For a 3-dimensional manifold $M^3$, which is a homology 3-sphere, Casson defined an integer valued invartiant, 
the Casson invariant, which is morally the `number' of flat 
$SU(2)$-connections on $M^3$.
The `virtual' dimension of the moduli space of flat 
$SU(2)$-connections on $M^3$ is $0$. However because of 
failure of transversality the number of flat connections in the naive sense 
may be infinite. Also, as in the case of intersection theory in 
differential geometry or topology, the `number' should be 
counted with sign. 
Casson's way to count the number of flat $SU(2)$-connections 
on $M$ uses Heegaard splitting of $M$ into the union of two 
handle bodies $M = H^1_g \cup_{\Sigma_g} H^2_g$.
Here $\Sigma_g$ is an oriented 2-dimensional manifold of genus $g$ and 
$H_g \cong H_g^1 \cong H_g^2$ are the handle bodies which bound $\Sigma_g$.
The moduli space $R(\Sigma_g)$ of flat $SU(2)$-connections
on the trivial bundle on $\Sigma_g$ is a singular space 
of dimension $6g-6$. 
The moduli spaces $R(H_g^i)$ of flat $SU(2)$-connections
on the trivial bundle on the handle bodies $H_g^i$
become  subspaces of $R(\Sigma)$ of dimension $3g-3$.
Casson defined Casson invariant $Z(M)$ 
as the `intersection number' 
$$
Z(M): = R(H_g^1) \cdot R(H_g^2) \in \mathbb Z
$$
of two $3g-3$ dimensional subspaces in the $6g-g$ dimensional 
space $R(\Sigma_g)$.
The points to be worked out are the following:
\begin{enumerate}
\item
The intersection number is well-defined 
even though $R(\Sigma_g)$
and $R(H_g^i)$ have singularities.
\item
The number $Z(M)$ is independent of the choices of Heegaard splitting
$M \cong H_g^1 \cup_{\Sigma_g} H_g^2$ and $g$. It becomes an 
invariant of the 3-dimensional manifold $M$.
\end{enumerate}
Casson proved them in the case when $H(M;\mathbb Z) = H(S^3;\mathbb Z)$.
Note that this assumption implies that the intersection points 
$R(H_g^1) \cap R(H_g^2)$ are not singular points of $R(\Sigma_g)$
unless it is the point corresponding to the trivial connection.

Taubes' work \cite{Tau} gave an alternative construction.
In place of using Heegaard splitting Taubes studied the 
set of all connections $\mathcal B(M;SU(2))$ (modulo the 
gauge transformation group) and regard 
the condition $F_A = 0$ (the curvature is $0$) as 
a differential equation. 
In the case when the set of flat connections on $M$ is not 
transversal this equation is not transversal.
Taubes then perturbed the equation $F_A = 0$
so that after adding an appropriate perturbation term $\mu(A)$ 
the set of solutions of the equation $F_A + \mu(A) = 0$ 
becomes isolated.
He then counted its order. 
This construction works under the assumption $H(M;\mathbb Z) = H(S^3;\mathbb Z)$
since otherwise there are reducible connections other than the trivial 
one, which causes trouble.
Taubes then proved that the invariant by such a count is 
equal to one obtained by Heegaard splitting.
\par\medskip
Defining an invariant of 3-dimensional manifolds is one of the applications of 
Floer homology. The other application is to symplectic geometry 
especially to Arnold's conjecture on the periodic orbits of 
a periodic Hamiltonian system.
Conley-Zehnder \cite{CZ} proved it in the case of (symplectic) torus 
$T^{2n}$, as follows.
Let $H : T^{2n} \times S^1 \to \R$ be a smooth function.
For $t\in S^1$, we put $H_t(x) = H(x,t)$ and let $\frak X_{H_t}$ be the 
Hamiltonian vector field associated to the function $H_t$ with respect 
to a certain symplectic structure $T^{2n}$ (with constant coefficient).
We consider the set $\mathcal{PER}_H$ of solutions of the equation
\begin{equation}\label{hameq}
\frac{d}{dt} \gamma = \frak X_{H_t} \circ \gamma
\end{equation}
where $\gamma : S^1  \to T^{2n}$. 
Conley-Zehnder proved that the order of $\mathcal{PER}_H$ is 
not smaller than $2^{2n}$, the Betti-number of $T^{2n}$,
under a certain non-degeneracy condition.
A solution of Equation (\ref{hameq}) can be regarded as a 
critical point of the action functional $\mathcal A_H$ defined by 
\begin{equation}\label{action}
\mathcal A_H(\gamma) 
= -\int_{D^2} u^* \omega + \int_{S^1} H(\gamma(t),t) dt.
\end{equation}
Here we consider only the loops $\gamma : S^1 \to T^{2n}$
which are homotopic to the constant map 
and $u : D^2 \to T^{2n}$ is a map such that
$u\vert_{\partial D^2} = \gamma$.
The integral of the symplectic form $\omega$ which 
is the first term of the right hand side
is independent of the choice of $u$, because of  Stokes' 
theorem, since $\pi_2(T^{2n}) = 0$.
The fact that a periodic solution of  
Hamilton equation (\ref{hameq}) is 
a critical point of the action functional $\mathcal A_H$ is 
a classical fact (maybe discovered by Hamilton himself).
However it had been difficult to use `Morse theory' of action functional $\mathcal A_H$
to study periodic solutions of  
Hamilton equation (\ref{hameq}).
In fact the properties of $\mathcal A_H$ are far from many of the functionals 
studied in geometric analysis. 
In the case when the functional $\frak F$ satisfies a condition 
that $\{ x \mid \frak F(x) \le c \}$ is `compact' 
in a certain weak sense (Palais-Smale's condition C is 
a typical way to formulate it), then 
one can show that the set of critical points of $ \frak F$ 
is related to the topology of the configuration space.
However in the case of the action functional $\mathcal A_H$ 
such `compactness' does not hold in any reasonable sense.
In fact the set of critical points of $\mathcal A_H$ is
related to the topology of $T^{2n}$ but not to the topology 
of the loop space $\Omega(T^{2n})$.
\par
Conley-Zehnder \cite{CZ} used a finite dimensional 
approximation of the loop space $\Omega(T^{2n})$ by 
Fourier expansion and used an  appropriate finite dimensional 
approximation of the action functional $\mathcal A_H$
to study the set of critical points of the action functional $\mathcal A_H$.
\par
Note that in gauge theory there exists a functional, 
the Chern-Simons functional: 
\begin{equation}\label{chernsimons}
\frak{cs}(A) = \int_M {\rm Tr}(A \wedge dA + \frac{2}{3} A \wedge A \wedge A)
\end{equation}
on the set of gauge equivalence classes of 
connections on a 3-dimensional manifold $M$, such that its critical point set 
coincides with the moduli space $R(M)$ of  flat connections on $M$.
Floer homology studies two functionals (\ref{action}) and (\ref{chernsimons})
in a similar way.
\par\medskip
Witten's paper \cite{Wi} had also an important impact to the discovery 
of Floer homology.
Witten explained how  Morse theory can be regarded as a 
(supersymmetric topological) field theory.
Using a Morse function $f : M \to \R$, Witten 
deformed a Laplace operator $\Delta$ on $p$-forms 
to 
$$
\Delta_t = d_t \circ d_t^* + d_t^*\circ d_t,
\qquad d_t = e^{-tf} d e^{tf}
$$
and found that:
\begin{enumerate}
\item 
For large $t$ the set of small eigen-spaces of $\Delta_t$ on $p$-forms can 
be identified with the vector space whose basis is 
identified with the set of critical points of $f$ with 
Morse index $p$.
\item
The restriction of $d_t$ to the set of small eigen-spaces of $\Delta_t$ 
defines a cochain complex which is isomorphic to $(C(M;f),d)$ where:
\begin{enumerate}
\item As a vector space $C(M;f)$ has a basis $\{[p] \mid p \in {\rm Crit}f\}$
where ${\rm Crit}$ is the set of critical points.
\item The matrix coefficient $\langle d[p],[q]\rangle$ is the 
number counted with sign of the integral curves of the gradient vector 
field ${\rm grad}f$ joining $p$ and $q$.
\end{enumerate}
\end{enumerate}
The chain complex defined by (2) above is called the Witten complex.
(Actually very similar constructions had been known in the classical 
works by Morse, Smale, Milnor etc. The importance of Witten's work 
is explaining its relation to quantum field theory, supersymmetry, and etc.)
\par
The important point of the Witten complex in Floer theory is the following:
It uses the moduli space of the solutions of the equation
\begin{equation}\label{gradienteq}
\frac{d \ell}{d\tau}(\tau) = {\rm grad}_{\ell(\tau)} f
\end{equation}
together with the asymptotic boundary conditions
\begin{equation}\label{asyboundary}
\lim_{\tau \to -\infty} \ell(\tau) = p, \quad 
\lim_{\tau \to +\infty} \ell(\tau) = q.
\end{equation}
Floer studied infinite dimensional versions where 
the Morse function $f$ is replaced by either the 
action functional $\mathcal A_H$ or the Chern-Simons functional $\frak{cs}$.
The equation (\ref{gradienteq}) then becomes 
\begin{equation}\label{Floersequ}
\frac{\partial u}{\partial \tau}
= J \left(
\frac{\partial u}{\partial t} - X_{H_t}
\right)
\end{equation}
or
\begin{equation}\label{GTFLoereq}
\frac{\partial A}{\partial \tau} = *_M F_A.
\end{equation}
Here $u : S^1 \times \R \to X$ is a map 
to a symplectic manifold $X$, $J$ is a compatible almost 
complex structure and
$A$ is a connection of a trivial $SU(2)$ bundle on $M \times \R$.
We take a gauge (temporal gauge) such that
$A$ has no $d\tau$ component, where $\tau$ is the coordinate of $\R$.
\par
Studying (\ref{Floersequ}) or (\ref{GTFLoereq}) 
with initial condition 
$u(0,t) =$ given, or $A\vert_{\tau=0}$ = given, 
is difficult.
Actually it is known that for almost all (smooth) initial values 
they do not have solutions.\footnote{Let us consider the case $X = \C$.
The initial value $u(0,t)$ has Fourier expansion 
$u(0,t) = \sum_{k \in \Z}a_k e^{2\pi\sqrt{-1}kt}$.
The solution $u(\tau,t)$ should be $u(\tau,t) = \sum_{k \in \Z}a_k e^{2\pi k(\tau+\sqrt{-1}t)}$.
For this series to converge for $\vert \tau\vert < \delta$ it is necessary 
$\vert a_k\vert < C e^{-2\pi\delta k}$. This condition is much more restrictive than 
the condition for $\sum_{k \in \Z}a_k e^{2 \pi\sqrt{-1}kt}$ to converge to a smooth map.}
In other words the gradient flow of $\mathcal A_H$ or 
$\frak{cs}$ is not well-defined.
\par
On the other hand, if we put an asymptotic boundary condition 
similar to (\ref{asyboundary}), the equations 
(\ref{Floersequ}) and (\ref{GTFLoereq}) behave 
nicely. Namely:
\begin{enumerate}
\item 
Its `weak solution' is automatically smooth.
\item
The moduli spaces of its solutions are finite dimensional.
\item
The moduli spaces of its solutions have nice compactifications,
which are similar to those of finite dimensional Morse theory.
\end{enumerate}
This is based on the fact that the equation (\ref{Floersequ}) 
is a variant of Gromov's pseudo-holomorphic curve equation (\ref{phceq2})
and (\ref{GTFLoereq}) is a particular case of the ASD-equation (\ref{ASDeq}).
\par
Let $X$ be a compact symplectic manifold
and $H: X \times S^1 \to \R$  a smooth function.
We denote by $\mathcal{PER}_H$ the set of solutions $\gamma : S^1 \to X$ of the equation 
(\ref{hameq}) which are homotopic to zero.
We assume that elements of $\mathcal{PER}_H$ satisfy an appropriate 
non-degeneracy condition.
We put
$$
CF(X;H) = \bigoplus_{\gamma \in \mathcal{PER}_H} \mathbb F [\gamma],
$$
where $\mathbb F$ is the coefficient ring which is explained later.

\begin{thm}\label{perHF}
There exists a boundary operator
$
d: CF(X;H) \to CF(X;H)
$
such that $d\circ d = 0$.
The Floer homology 
$$
HF(X;H): = \frac{{\rm Ker}d}{{\rm Im}d}
$$
is isomorphic to the ordinary homology $H(X;\mathbb F)$ with $\mathbb F$ 
coefficients.
\end{thm}
\begin{cor}\label{cor32}
In the situation when  elements of $\mathcal{PER}_H$ are all 
non-degenerate, we have
$$
\# \mathcal{PER}_H \ge {\rm rank} H(X;\mathbb F).
$$
\end{cor}

Floer \cite{fl3} proved Theorem \ref{perHF} in the case 
when $X$ is monotone.
Here a symplectic manifold $(X,\omega)$ is said to be monotone if 
there exists a positive number $c$ such that
\begin{equation}\label{monotoneabs}
c \int_{S^2} u^* \omega = u_*([S^2]) \cap c_1(X)
\end{equation}
for all $u : S^2 \to X$.
In that case $\mathbb F = \Z$.
This assumption is relaxed by Hofer-Salamon \cite{HS} and Ono \cite{On}
to the semi-positivity.
Here $(X,\omega)$ is said to be semi-positive 
if there does not exist $u : S^2 \to X$ such that
$$
\int_{S^2} u^* \omega > 0, \qquad 0 > u_*([S^2]) \cap c_1(X) \ge 6 - 2n.
$$
In this case, $\mathbb F$ is a Novikov ring (with $\Z$ as a ground ring).
There are several variants of the definition of a Novikov ring.\footnote{
This ring itself is known before Novikov. Novikov \cite{nov} first pointed out that 
to study Morse theory of closed one form (which is not necessary exact) 
we need to use this ring. 
In the case $\alpha \mapsto \int_{\alpha}\omega$ is non-zero 
on $\pi_2(X)$, the action functional (\ref{action}) is not single valued.
So Floer theory (of, say, periodic Hamiltonian system) should be regarded as a 
Morse theory of closed 1 form.}
A version which is called (the universal) Novikov ring (with the ground ring $R$)
is the set of all formal sums
\begin{equation}\label{eleNov}
\sum_{i=0}^{\infty} a_i T^{\lambda_i}
\end{equation}
where $a_i \in R$ and $\lambda_i \in \R_{\ge 0}$ with 
$\lim_{i\to \infty} \lambda_i = + \infty$
 (\cite{fooobook}).
The universal Novikov ring with $R$ as the  ground ring is written as 
$\Lambda_0^R$. 
\par
In the case when $\mathbb F$ is the universal Novikov ring with the ground ring $\Q$,
Theorem \ref{perHF} is proved by 
Fukaya-Ono \cite{FO}, Liu -Tian \cite{LT}, Ruan \cite{Ru}.
\par
Note that Theorem \ref{perHF} for $\mathbb F$ to be the universal Novikov ring with 
the ground ring $R$ implies Corollary \ref{cor32} with $\mathbb F = R$.
In the case when $\mathbb F$ is a finite field 
Corollary \ref{cor32} is proved in a recent paper by
Abouzaid-Blumberg \cite{AB}.
See also \cite{BX}.
\par
All of those proofs use Morse theory of the functional $\mathcal A_H$ 
and the equation (\ref{Floersequ}) to define 
Floer homology.
The difference between the methods of papers mentioned above lies on the 
way to overcome various difficulties appearing in the infinite 
dimensional situations. We do not discuss it here.
\par
To prove that Floer homology $HF(X;H)$ is isomorphic to the 
ordinary homology there are three different methods established in the literature.
\begin{enumerate}
\item 
\begin{enumerate}
\item We relax the condition that the periodic orbits of $\mathfrak X_H$ are 
non-degenerate, so that the case $H=0$ will be included. 
\item We show the Floer homology $HF(X;H)$ is independent of $H$ in that generality.
\item We prove that in case $H=0$ Floer homology $HF(X;0)$ is isomorphic 
to the ordinary homology.
\end{enumerate}
\item
We study the case when $H : X \times S^1 \to \R$ is independent of $S^1$ factor
and so is a function on $X$. We furthermore require that $H$ is a Morse function 
and its $C^2$-norm is sufficiently small. Then we show that 
the boundary operator $d$ to define the Floer homology $HF(X;H)$ is equal to the 
boundary operator of the Witten complex of $CF(X;H)$ (the one of finite 
dimensional Morse theory).
\item
\begin{enumerate}
\item
We study two Lagrangian submanifolds in 
$(X\times X,-\pi_1^*\omega+\pi_2^*\omega)$.
One is the diagonal $\Delta = \{(x,x) \mid x \in X\}$ and the other 
is the graph ${\rm Gra}_{\varphi_H^1} : \{(x,\varphi_H^1(x) \mid x \in X\}$
of $\varphi_H^1$.
Here $\varphi_H^t : X \to X$ is defined by 
\begin{equation}\label{varphiH}
\varphi^0(x) = x, 
\qquad
\frac{d\varphi_H^t(x)}{dt} = \frak X_{H_t}(\varphi_H^t(x)).
\end{equation}
\item
We show that the Lagrangian Floer homology\footnote{See Section \ref{sec:lagHF}.} $HF(\Delta,{\rm Gra}_{\varphi_H^1})$
is well-defined, isomorphic to $HF(X;H)$, independent of $H$ and is isomorphic to $H(X)$.
\end{enumerate}
\end{enumerate}
To prove that $HF(X;H)$ is isomorphic to the 
ordinary homology in the case when $\mathbb F$ is the universal 
Novikov ring with the ground ring $\Q$, 
the method (1) is used in \cite{LT},\cite{Ru}, 
the method (2) is used in \cite{FO}.
The method (3) is worked out later in \cite{fooobook} and \cite{fooo:inv}.
Actually there are two methods which can be used to show that 
$HF(\Delta,{\rm Gra}_{\varphi_H^1})$ is isomorphic to $H(X)$.
One uses the fact that $H(\Delta) \to H(X\times X)$ is injective.
The other uses the anti-holomorphic involution $X\times X \to X\times X$ 
for which $\Delta$ is the fixed point set.
The first method is used in \cite{fooobook}, \cite{fooo:inv}.
The second method works under a certain assumption on $X$ 
when the ground ring is $\Z_2$.
\par\medskip
In Yang-Mills gauge theory (Donaldson-Floer theory) 
Floer homology (instanton homology) is defined 
in one of the following two cases. \footnote{This condition is 
called admissibility in certain references.}
\begin{enumerate}
\item[(AIF1)]
$M$ is a 3-dimensional closed manifold such that $H(M;\Z) \cong H(S^3;\Z)$
and
$\mathcal E_M \to M$ is the trivial $SU(2)$ bundle.
\item[(AIF2)]
$M$ is a 3-dimensional oriented closed manifold.
$\mathcal E_M \to M$ is a principal $SO(3)$ bundle. There exists a 2-dimensional submanifold 
$\Sigma \subset M$ such that the restriction of $\mathcal E_M$ to $\Sigma$ is 
non-trivial.
\end{enumerate}
The Chern-Simons functional (\ref{chernsimons}) can also be defined 
in the case (AIF2) such that its gradient flow equation is (\ref{GTFLoereq}).
\par
We consider the set $R(M)$ of gauge equivalence classes of flat connections 
 of $\mathcal E_M \to M$. 
In the case (AIF2) all the elements $[a]$ of $R(M)$ are irreducible,
that is, the bundle automorphism of $\mathcal E_M$ preserving $a$ is trivial.
In the case (AIF1) all the elements $[a]$ of $R(M)$ except $[a] = [0]$ are irreducible.
The reducible connections correspond to the singularity of the set of gauge 
equivalence classes of connections. (AIF1),(AIF2) are used to go around the 
trouble which the singularity of the set of gauge 
equivalence classes causes.
\par
Let $\mathcal B(M;\mathcal E_M)$ be the set of gauge 
equivalence classes of connections on $\mathcal E_M$.
We can define an appropriate function $h: \mathcal B(M;\mathcal E_M) \to \R$,
such that the set of  solutions of the perturbed equation 
\begin{equation}\label{perturbcurv}
*_M F_a + D_ah = 0
\end{equation}
is isolated. Here $*_M$ is the Hodge $*$ operator and $D_ah$ is the 
derivative of $h$ at $a$. The two terms in (\ref{perturbcurv}) are sections of 
$\Lambda^1 \otimes ad\mathcal E_M$.
Here $\mathcal E_M$ is the $su(2) = so(3)$ bundle 
induced by the adjoint representation from the principal 
bundle $\mathcal E_M$.
We also required that the linearized operator 
$$
\mathcal D_a: = *_M d_a + {\rm Hess}_a h
$$
is invertible for solutions $a$ of (\ref{perturbcurv}).
Let $R(M;h)$ be the set of gauge equivalence classes of solutions of 
(\ref{perturbcurv}). In case (AIF1) we remove the trivial connection 
from $R(M;h)$. We put
$$
CF(M,\mathcal E_M;h)
= \bigoplus_{[a] \in R(M;h)} \Z [a].
$$
Floer defined a boundary operator 
$d: CF(M,\mathcal E_M;h) \to CF(M,\mathcal E_M;h)$
by
\begin{equation}\label{Flbondarydef}
d[a] = \sum_{a'} \langle da,a'\rangle [a']
\end{equation}
where the matrix element $\langle da,a'\rangle$ is the number 
counted with sign of solutions of the equation
\begin{equation}\label{GTFLoereqpert}
\frac{\partial A}{\partial \tau} = * F_{A\vert_{M\times \{\tau\}}} + D_{A\vert_{M\times \{\tau\}}}h
\end{equation}
with asymptotic boundary conditions:
\begin{equation}
\lim_{\tau\to-\infty}[A\vert_{M\times \{\tau\}}] = [a],
\qquad
\lim_{\tau\to+\infty}[A\vert_{M\times \{\tau\}}] = [a'].
\end{equation}
\begin{thm}{\rm(Floer)}\label{gaugeFloer}
$d \circ d = 0$. The cohomology, called instanton (Floer) homology
$$
I(M;\mathcal E_M): = \frac{{\rm Ker}d}{{\rm Im}d}
$$
is independent of $h$ and is an invariant of 
a 3-dimensional manifold $M$ equipped with $\mathcal E_M$.
\end{thm}
Floer proved Theorem \ref{gaugeFloer} in \cite{fl1} in the case (AIF1)
and in \cite{fl5} in the case (AIF2).

\section{Topological field theory.}
\label{HFtoposec}

An important development of Floer homology in gauge theory
is the discovery of its relation to 4-dimensional Donaldson 
invariant. (This is due to Donaldson and Floer and is explained 
in \cite{Don5}.)
Let $X$ be a 4-dimensional manifold with boundary $M = \partial X$ and 
$\mathcal E_X \to X$ is a, say, $SU(2)$ bundle over $X$.
Suppose that the restriction of $\mathcal E_X$ to $M$ is trivial.
Then, under a certain hypothesis, one can define a 
relative Donaldson invariant as follows.
\par
Let $[a]$ be a gauge equivalence class of a flat connection $a$ on $M$.
We take a Riemannian metric on $\overset{\circ}X: = 
X \setminus \partial X$ such that 
$\overset{\circ}X$ minus a compact set is isometric to 
$M \times (0,\infty)$. We consider the moduli space of 
connections $A$ on $\mathcal E_X$ which solves\footnote{Here $F_A$ and $*_X$ 
denote the curvature of the connection $A$
and the Hodge star operator of $X$, respectively.}:
\begin{equation}\label{ASD2}
F_A + *_X F_A = 0
\end{equation}
and satisfies the asymptotic boundary condition
\begin{equation}\label{asympto11}
\lim_{\tau \to \infty} [A\vert_{M \times \{\tau\}}] = [a].
\end{equation}
We also require that the energy $\mathcal E(A): = \int_X \Vert F_A\Vert^2$ is finite.
The moduli space $\mathcal M(X;a;E)$ of such connections with given $E = \mathcal E(A)$ 
becomes a finite dimensional space and has a nice compactification.
In a simplest case when the virtual dimension is zero, it gives an element
\begin{equation}\label{relativeinv}
\sum_{a,E_a} \# \mathcal M(X;a;E_a) [a] \in CF(M).
\end{equation}
Here the sum is taken over $a,E_a$ such that the 
virtual dimension of $\mathcal M(X;a;E_a)$ is zero.
\par
Using the fact that (\ref{GTFLoereq})  and (\ref{ASD2}) 
coincide on $M \times (0,\infty)$, we can show that 
(\ref{relativeinv}) is a cycle with respect to the 
boundary operator (\ref{Flbondarydef}) and so obtain a 
relative invariant
in the (instanton) Floer homology $I(M;{\rm trivial})$.
In case the (virtual) dimension of $\mathcal M(X;a;E)$ is positive 
we cut $\mathcal M(X;a;E_a)$ using  homology classes of the 
space of connections on $X$ (typically obtained from homology classes of 
$X$) in the same way as the case of Donaldson invariant (\ref{forDona}) 
and obtain a relative invariant.
\par
This construction becomes a prototype of the definition of 
topological field theory (\cite{witten1, At}), which might 
be formulated as follows.\footnote{
The description below is not intended to formulate a precise axiom. It is rather 
an informal guideline how  such a theory will be built. More precise and 
systematic formulation of topological field theory is now being built.}
\begin{enumerate}
\item[(TF1)]
To a closed oriented $n$-dimensional manifold $X^n$ it associates a number $Z_X$.
\item[(TF2)]
To a closed oriented $(n-1)$-dimensional 
manifold $M^{n-1}$ it associates a vector space $HF(M)$
with inner product, such that $HF(-M)= HF(M)^*$ 
(where $V^*$ denotes the dual vector space of $V$)
and $HF(M_1 \sqcup M_2) = HF(M_1) \otimes HF(M_2)$.
\item[(TF3)] 
Let $X$ be an oriented $n$-dimensional manifold such that 
$X$ minus a compact set is the union of $M_- \times (-\infty,0)$ 
and $M_+ \times (0,+\infty)$.
Then it associates a linear map:
$$
Z_X : HF(M_-) \to HF(M_+).
$$
\item[(TF4)] 
Let $M_1,M_2,M_3$ be closed oriented $(n-1)$-dimensional manifolds and 
$X_{ij}$ be an oriented $n$-dimensional manifold for $(ij) = (12)$ or $(23)$
such that
$X_{ij}$ minus a compact set is the union of $-M_i \times (-\infty,0)$ 
and $M_j \times (0,+\infty)$.
\par
We glue $X_{12}$ and $X_{23}$ along $M_2$ and obtain 
$X_{13}$ such that $X_{13}$ minus a compact set is the union of $-M_1 \times (-\infty,0)$ 
and $M_3 \times (0,+\infty)$.
Then we have
$$
Z_{X_{13}} = Z_{X_{23}} \circ Z_{X_{12}}: HF(M_1) \to HF(M_3).
$$
\end{enumerate}
We remark that Donaldson-Floer theory actually does {\it not} 
satisfy this axiom itself.
In fact the instanton Floer homology is defined only under a certain 
assumption on 3-dimensional manifolds and Donaldson invariant in general uses a 
certain auxiliary data 
(such as a homology class of a 4-dimensional manifold) 
and in a certain case ($b_2^+ = 1$) it depends 
on the `chamber'. Moreover it is not defined in a certain case ($b_2^+ = 0$).
It seems that such delicate `unstable' phenomenon is {\it the} reason 
why this theory is so nontrivial and powerful.
The `axiomatic' understanding of Donaldson-Floer theory seems 
to be a subject yet to be studied and clarified in the future.
\par
In the early 1990's various people expected that 
relative Donaldson invariant (Donaldson-Floer theory) 
will be a tool to calculate Donaldson invariant, 
via decomposing 4-dimensional manifolds into pieces.
However the mathematical study of gauge theory  
developed in a different way.
The major tool to compute Donaldson invariant
turns out to be Kronheimber-Mrowka's structure theorem (\cite{KM}) 
and 
its relation to Seiberg-Witten invariant (\cite{witten3}).
\par
In the early 1990's there was also an attempt to expand the 
topological field theories to those on 
$n$-$(n-1)$-$(n-2)$ dimensional theory.
It may be formulated as follows.
\begin{enumerate}
\item[(TF5)]
Let $N$ be an $(n-2)$-dimensional closed oriented manifold.
To $N$ the topological field theory associates a category $\mathscr C(N)$.
For two objects $c,c'$ of $\mathscr C(N)$, the set of morphisms 
$\mathscr C(N)(c,c')$ is a vector space with an inner product.
The category $\mathscr C(-N)$ associated to $-N$ 
(Here $-N$ is the manifold $N$ with the opposite orientation.) is the opposite category 
$\mathscr C(N)^{\rm op}$. The set of objects of $\mathscr C(N)^{\rm op}$ 
is identified with the set of objects of  $\mathscr C(N)$.
For two objects $c,c'$ 
$$
\mathscr C(N)^{\rm op}(c,c') = \mathscr C(N)(c',c).
$$
\item[(TF6)]
Let $M$ be an oriented $(n-1)$-dimensional manifold such that 
$M$ minus a compact set is the union of $N_- \times (-\infty,0)$ 
and $N_+ \times (0,+\infty)$.
Then the topological field theory associates a functor:
$$
HF_M : \mathscr C(N_-) \to \mathscr C(N_+).
$$
\item[(TF7)] 
Let $N_1,N_2,N_3$ be closed oriented $(n-2)$-dimensional manifolds and 
$M_{ij}$  an oriented $(n-1)$-dimensional manifold for $(ij) = (12)$ or $(23)$
such that
$M_{ij}$ minus a compact set is the union of $-N_i \times (-\infty,0)$ 
and $N_j \times (0,+\infty)$.
\par
We glue $M_{12}$ and $M_{23}$ along $N_2$ and obtain 
$M_{13}$ such that $M_{13}$ minus a compact set is the union of $-N_1 \times (-\infty,0)$ 
and $N_3 \times (0,+\infty)$.
Then we have
$$
HF_{M_{13}} = HF_{M_{23}} \circ HF_{M_{12}}: \mathscr C(N_1) \to \mathscr C(N_3).
$$
See \cite[Definition 8.5]{takagi} for the case $N_1 = \emptyset$ and/or $N_3 = \emptyset$.
\end{enumerate}
See \cite[Theorem 3.2]{fu1} for the formulation in the case 
when $X$ is an $n$-dimensional manifold with corners 
such that $\partial X = M_- \cup M_+$
and $M_- \cap M_+ = N$.

In the case of Donaldson-Floer theory (Yang-Mills gauge theory), 
Donaldson proposed a candidate of the 
category $\mathscr C(\Sigma)$ to be associated to a 2-dimensional manifold 
in the year 1992\footnote{during a conference 
at University Warwick. At the same conference 
Y. Ruan explained his work \cite{Ru1} to define 
an invariant of a symplectic manifold using the fundamental class of 
the moduli space of pseudo-holomorphic curves. This idea 
was not  explicit before. (It was implicit in Gromov's work.)} as follows. (Here we write $\Sigma$ in place of $N=N^{4-2}$.)
\par
Let $(M,\mathcal E_M)$ be a pair of 3-dimensional manifold 
with boundary and an $SU(2)$ or $SO(3)$ bundle on it.
We consider one of the following two situations:

\begin{enumerate}
\item[(AIFB1)]
$\mathcal E_M$ is a trivial $SU(2)$ bundle.
\item [(AIFB2)]
$\mathcal E_M$ is an $SO(3)$ bundle.
The second Stiefel-Whitney class of the restriction 
of $\mathcal E_M$ to the boundary $\Sigma = \partial M$ 
is the fundamental class $[\Sigma] \in H^2(\Sigma;\Z_2)$.
\end{enumerate}
Note that (AIFB2) implies that the number of connected components 
of $\Sigma$ is even.
\par
We denote by $\mathcal E_{\Sigma}$ the restriction of $\mathcal E_{M}$
to $\Sigma = \partial M$.
Let $R(\Sigma;\mathcal E_{\Sigma})$ be the space of 
gauge equivalence classes of the flat connections of $\mathcal E_{\Sigma}$.
\par
In case (AIFB2), the space $R(\Sigma;\mathcal E_{\Sigma})$ is a smooth manifold 
and in case (AIFB1), the space $R(\Sigma;\mathcal E_{\Sigma})$ 
has a singularity. In both cases its dimension is $6g-6$ where 
$g$ is the genus of $\Sigma$. In the disconnected case 
$R(\Sigma;\mathcal E_{\Sigma})$ is the direct product of 
the spaces $R(\Sigma_a;\mathcal E_{\Sigma_a})$ for connected components $\Sigma_a$ of
$\Sigma$.
\par
(The regular part of)  $R(\Sigma;\mathcal E_{\Sigma})$ has a symplectic 
structure \cite{go}.
In fact the tangent space at $[a]$ of $R(\Sigma;\mathcal E_{\Sigma})$
is identified with the first cohomology $H^1(\Sigma;ad(a))$ of the 
flat $su(2) = so(3)$ bundle associated to the principal bundle 
$\mathcal E_{\Sigma}$ by the adjoint representation.
The cup product defines an anti-symmetric form on $H^1(\Sigma;ad(a))$
which we can check to be a symplectic form.
\par
Let $R(M;\mathcal E_M)$ be the set of 
gauge equivalence classes of  flat connections of $\mathcal E_{M}$.
The restriction of a connection defines a map
$$
{\rm Res} : R(M;\mathcal E_M) \to R(\Sigma;\mathcal E_{\Sigma}).
$$
By Stokes' theorem we can show
$$
{\rm Res}^*\omega = 0,
$$
where $\omega$ is the symplectic form on $R(\Sigma;\mathcal E_{\Sigma})$.
Let $\mathcal B(\Sigma;\mathcal E_{\Sigma})$ be the space of 
gauge equivalence classes of connections of $\mathcal E_{\Sigma}$.
In a similar way as 
(\ref{perturbcurv}) we can find an appropriate 
perturbation $h: \mathcal B(\Sigma;\mathcal E_{\Sigma}) \to \R$, 
such that $h(A)$ depends only on a restriction of $A$ to a complement 
of a neighborhood of $\partial M$, such that the following holds.
Let $R(M;\mathcal E_M;h)$ be the space of gauge equivalence classes of solutions of the equation:
\begin{equation}\label{GTFLoereqpertpartial}
\frac{\partial A}{\partial \tau} = * F_A + D_{A\vert_{M\times \{\tau\}}}h.
\end{equation}
\begin{prop}{\rm (Herald \cite{He}. See also \cite{DFL2}.)}
For a generic choice of $h$ the space $R(M;\mathcal E_M;h)$
has dimension $\frac{1}{2}\dim R(\Sigma;\mathcal E_{\Sigma};h)$.
The map
$$
{\rm Res} : R(M;\mathcal E_M;h) \to R(\Sigma;\mathcal E_{\Sigma})
$$
becomes a Lagrangian immersion outside the set of singular points of 
$R(\Sigma;\mathcal E_{\Sigma})$.
\end{prop}
In the case (AIFB1), the space $R(M;\mathcal E_M;h)$ contains a reducible connection,
where it becomes singular.
In the case (AIFB2), the space $R(M;\mathcal E_M;h)$ does not contain a 
reducible connection and is a smooth manifold.
In the latter case, ${\rm Res}$ becomes a Lagrangian immersion.
\par
The candidate proposed by Donaldson for $\mathscr C(\Sigma)$ 
is one whose object is a Lagrangian submanifold of 
$R(\Sigma;\mathcal E_{\Sigma})$ and 
the morphisms are Lagrangian Floer homology. (See Section \ref{sec:lagHF}).
\par
This proposal is related to several results and conjectures which appeared 
around the same time (early 1990's).
\par
Let us first consider the case when the 3-dimensional manifold is the 
handle body $H_g$ and the case (AIFB1).
The map
$$
{\rm Res} : R(H_g,\mathcal E_{H_g}) \to R(\Sigma_g;\mathcal E_{\Sigma_g})
$$
for the trivial bundle $\mathcal E_{H_g}$ is a Lagrangian embedding 
outside the set of reducible connections.
\par
Let $M = H_g^1 \cup_{\Sigma_g} H_g^2$ be a Heegaard decomposition of a 
homology 3-sphere $M$.
The instanton (Floer) homology $I(M)$ is defined by using the trivial 
$SU(2)$ bundle. (Theorem \ref{gaugeFloer}.)

\begin{conj}{\rm (Atiyah-Floer conjecture \cite{At})}\label{AFconj}
The instanton (Floer) homology $I(M)$ is isomorphic to the 
Lagrangian Floer homology $HF(R(H^1_g;\mathcal E_{H_g}),R(H^1_g;\mathcal E_{H_g}))$.
\end{conj}

Actually the statement itself has a difficulty. In fact 
since $ R(\Sigma_g;\mathcal E_{\Sigma_g})$ is singular 
the definition of Lagrangian Floer homology in Conjecture 
\ref{AFconj} is not yet established. 
\par
Note that Casson's definition of Casson invariant uses 
Heegaard decomposition and Taubes' version is based on 
gauge theory. Therefore Conjecture \ref{AFconj} can be 
regarded as a `categorification' of Taubes' theorem that two 
definitions coincide.
\par
We like to mention that other than those we describe in this 
article, there are various approaches to Conjecture \ref{AFconj} 
by various mathematicians, 
such as \cite{Yo,LLW,Weh,MWo,dun}.
\par\medskip
Since the case (AIFB1) has a difficulty, we discuss the case (AIFB2).
We consider $\Sigma_g$ a genus $g$ 2-dimensional oriented manifold 
and an $SO(3)$ bundle $\mathcal E_{\Sigma_g}$ on it 
such that $w_2(\mathcal E_{\Sigma_g}) = [\Sigma_g]$.
We put $M = \Sigma_g \times [0,1]$ and consider the 
$SO(3)$ bundle $\mathcal E_M$ induced from   $\mathcal E_{\Sigma_g}$.
In this case  $R(\Sigma;\mathcal E_{\Sigma})$ is a smooth 
symplectic manifold.
Note that $\partial M$ is the disjoint union of two copies of $\Sigma_g$.
Therefore 
$R(\partial M;\mathcal E_{\partial M}) 
= -R(\Sigma_g;\mathcal E_{\Sigma_g}) \times R(\Sigma_g;\mathcal E_{\Sigma_g})$.
Here we put the minus sign to the first factor. 
It means that the symplectic form of the first factor is $-\omega$ and 
the one of the second factor is $\omega$.
In fact the  two copies of $\Sigma_g$ in $\partial M$  have opposite 
induced orientation 
and the symplectic structure on $R(\Sigma_g;\mathcal E_{\Sigma_g})$ 
changes the sign if we change the orientation of $\Sigma_g$.
The map
$$
{\rm Res} : R(M;\mathcal E_{M}) \to -R(\Sigma_g;\mathcal E_{\Sigma_g}) \times R(\Sigma_g;\mathcal E_{\Sigma_g})
$$
is the diagonal embedding.
In this situation an analogue of Conjecture \ref{AFconj} is proved by Dostoglou-Salamon
\cite{DS} as follows.
We consider
$M_1 = M_2 = \Sigma_g \times [0,1]$.
Then
$
\partial M_1 = -\partial M_2 
$
is a disjoint union of two copies of $\Sigma_g$, 
which we write $\Sigma_g^1 
\sqcup -\Sigma^2_g$. We take a diffeomorphism 
$F: \partial M_1 \cong -\partial M_2$ as follows.
$F = {\rm identity}$ on $\Sigma_g^1$ and 
$F = \varphi$ on $\Sigma_g^2$, where 
$\varphi$ is a certain orientation preserving diffeomorphism.
Then 
$$
M = M_1 \cup_{F} M_2
$$
is a mapping cylinder of $\varphi$ and is 
an $\Sigma_g$ bundle over $S^1$.
The bundle $\mathcal E_{\Sigma_g}$ 
induces $\mathcal E_{M}$ on $M$.
\par
The diffeomorphism $\varphi$ induces a symplectic diffeomorphism
$\varphi_* : R(\Sigma_g;\mathcal E_{\Sigma_g}) \to R(\Sigma_g;\mathcal E_{\Sigma_g})$.
We consider two Lagrangian submanifolds of 
$-R(\Sigma_g;\mathcal E_{\Sigma_g}) \times R(\Sigma_g;\mathcal E_{\Sigma_g})$:  
one is
the diagonal
$$
\Delta = \{(x,x) \mid x \in R(\Sigma_g;\mathcal E_{\Sigma_g})\}
$$
the other is
$$
{\rm Gra}(\varphi_*) = 
\{(x,\varphi_*(x)) \mid x \in R(\Sigma_g;\mathcal E_{\Sigma_g})\}.
$$
Since $\pi_2(-R(\Sigma_g;\mathcal E_{\Sigma_g}) \times R(\Sigma_g;\mathcal E_{\Sigma_g}),\Delta)$ and 
$\pi_2(-R(\Sigma_g;\mathcal E_{\Sigma_g}) \times R(\Sigma_g;\mathcal E_{\Sigma_g}),{\rm Gra}(\varphi_*))$
are $0$, Lagrangian Floer homology 
$$
HF(\Delta,{\rm Gra}(\varphi_*))
$$
is defined. (It is a $\Z_4$ graded $\Z$ module.)
\begin{thm}{\rm (\cite{DS})}\label{Dosala}
$I(M;\mathcal E_{M}) \cong HF(\Delta,{\rm Gra}(\varphi_*))$.
\end{thm}
This theorem can be regarded as a special case of (TF7).
\begin{rem}
Note that, in \cite{DS}, Theorem \ref{Dosala} is stated in a different way.
For a symplectic diffeomorphism $\varphi : X \to X$ one can define 
an analogue $HF(X;\varphi)$ of the Floer homology of periodic Hamiltonian 
system. Namely  $HF(X;\varphi)$ is a cohomology group of 
a chain complex $CF(X;\varphi)$ whose generator is a fixed point of 
$\varphi$. In the case when $X$ is a monotone symplectic manifold 
and $c_1(X)$ is divisible by $N$ the Floer homology $HF(X;\varphi)$ 
is a $\Z$ module with period $2N$. 
Dostoglou and Salamon proved $I(M;\mathcal E_{M}) \cong HF(X;\varphi_*)$.
\end{rem}
\par
In the case $\Sigma = T^2$ the next result written in 
Braam-Donaldson \cite{BD}\footnote{Braam-Donaldson 
attributes it to Floer \cite{fl45}.} is regarded as another  special case of (TF.7).
We consider a nontrivial $SO(3)$ bundle $\mathcal E_{T^2}$ 
on the two-dimensional torus $T^2$. 
It is easy to see that $R(T^2;\mathcal E_{T^2})$
the space of flat connections on $T^2$ consists of a single 
point. 
Let $M$ be a 3-dimensional manifold whose 
boundary is a disjoint union of $T^2$'s.
Suppose that $\mathcal E_M$ is an $SO(3)$ bundle on $M$ 
such that its second Stiefel--Whitney class $w_2(\mathcal E_M)$
restricts to the fundamental class of $\partial M$.
It implies that the number of boundary components of $\partial M$
is even. We take an orientation reversing involution $\tau$ of $\partial M$ which 
induces a free $\Z_2$ action on
$\pi_0(\partial M)$.
Then we glue $T^2 \subset \partial M$ with $\tau(T^2) \subset \partial M$
for each connected component and obtain 
a closed $3$ manifold, which we denote by $M_{\tau}$.
The $SO(3)$ bundle $\mathcal E_M$ induces an 
$SO(3)$ bundle $\mathcal E_{M_{\tau}}$
on $M_{\tau}$ in an obvious way.
\begin{thm}
\label{exision}{\rm (Floer, Braam-Donaldson)}
The instanton Floer homology $I(M_{\tau};\mathcal E_{M_{\tau}})$
is independent of the choice of $\tau$.
\end{thm}
We may regard this result as a  special case of (TF7)
as follows.
The space $R(\partial M;\mathcal E_{\partial M})$ 
is one point.
`Relative invariant' in this case is a 
chain homotopy type of a chain complex,
$CF(M;\mathcal E_{\partial M})$.
The gluing axiom (TF7) 
claims that for any $\tau$, 
the instanton Floer homology $I(M_{\tau};\mathcal E_{M_{\tau}})$,
is isomorphic to the cohomology of $CF(M;\mathcal E_{\partial M})$.
\par
There are two other results which are the cases when the 
bundle $\mathcal E_{\Sigma}$ on a 2-dimensional submanifold is trivial.
One is the case when $\Sigma = S^2$.  
This case corresponds to the study of the Floer homology of the connected sum 
$M_1 \# M_2$, and is studied in \cite{fu15,WL}, where it is proved that 
there exists a spectral sequence which relates $I(M_1)$, $I(M_2)$ and 
$I(M_1 \# M_2)$.
Note that $R(S^2,{\rm trivial})$ is one point and the isotropy group 
of the action of the gauge transformation group on this point is $SU(2)$ or 
$SO(3)$.  
\par
The other  is the case when $\Sigma = T^2$ and $\mathcal E_{\Sigma}$ is 
trivial. Floer studies this case and obtained an 
exact triangle which relate three instanton Floer homologies 
corresponding three different ways to identify $T^2$ with $\partial(D^2 \times S^1)$.
It is called a Floer's exact triangle. See \cite{fl5,BD}.
Note that $R(T^2;{\rm trivial})$ is $T^2/\Z_2$ and the isotropy group of the 
action of the gauge transformation group at the generic point is $U(1)$.
\par
From those three cases where $\Sigma = S^2$ or $T^2$, we find that 
for gluing axiom (TF7) to hold we need to modify Donaldson's 
proposal and include more general objects than 
Lagrangian submanifolds of $R(\Sigma;\mathcal E_{\Sigma})$ as 
objects of the category $\mathscr C(\Sigma;\mathcal E_{\Sigma})$.
In fact, in the situation of Theorem \ref{exision},
the object of $\mathscr C(T^2,{\rm nontrivial})$ is a chain complex.
So we need a kind of mixture of Lagrangian submanifold 
and chain complex. 
See Section \ref{sec:Ainfinity} for a way to obtain 
such a category.
In the case when $\Sigma = S^2$ or $T^2$ with trivial bundle, 
the way to relate the connected sum formula or Floer's 
triangle to (TF7) is not yet understood.
The difficulty is the fact that the isotropy group of the generic 
point of $R(\Sigma;{\rm trivial})$ has positive dimension in those cases.

\section{Lagrangian Floer theory.}
\label{sec:lagHF}

Among various Floer theories, Lagrangian Floer theory 
is the first  Floer studied \cite{fl2}. 
However actually the foundation of Lagrangian Floer theory is 
more delicate than other Floer theories such as 
Floer homologies of periodic Hamiltonian systems which we explained in Section \ref{sec:Floer}.
\par
Let $(X,\omega)$ be a compact symplectic manifold
and $L_i \subset X$  an embedded Lagrangian submanifold
for $i=0,1$.
We consider the space
$$
\Omega(L_0,L_1) = \{\gamma ; [0,1] \to X \mid \gamma(0) \in L_0, 
\gamma(1) \in L_1\}
$$
of arcs joining $L_0$ to $L_1$. 
To each connected component $\Omega(L_0,L_1)_o$ of $\Omega(L_0,L_1)$ we fix 
a base point $\gamma_o \in \Omega(L_0,L_1)_o$.
For $\gamma \in \Omega(L_0,L_1)_o$ we take a path 
joining $\gamma_o$ to $\gamma$. Such a path may be regarded as a 
map $u: \R \times [0,1] \to X$ such that:
\begin{enumerate}
\item[(path1)]
$u(\tau,0) \in L_0$, $u(\tau,1) \in L_1$.
\item[(path2)]
$\lim_{\tau\to-\infty} u(\tau,t)=\gamma_o(t)$.
\item[(path3)]
$\lim_{\tau\to+\infty} u(\tau,t)=\gamma(t)$.
\end{enumerate}
We define the action functional 
$\mathcal A$ by
\begin{equation}
\mathcal A(\gamma) = 
\int_{\R \times [0,1]} u^*\omega \in \R.
\end{equation}
Condition (path1) and Stokes' theorem imply 
that $\mathcal A(\gamma)$ depends only on the homotopy 
class of the path $u$ joining $\gamma_o$ to $\gamma$. 
It may depend on the homotopy class of $u$.
So $\mathcal A$ is a function on an appropriate 
covering space of $\Omega(L_0,L_1)$.
Its derivative however is well-defined.
Floer homology of Lagrangian submanifolds uses 
the gradient vector field of $\mathcal A$.
We take an almost complex structure $J$ of $X$ 
such that
$g(V,W) = \omega(V,JW)$ becomes a Riemannian metric.
We use it to define an $L^2$ norm of the 
section of $\gamma^*TX$ for $\gamma \in \Omega(L_0,L_1)$.
The space of sections of $\gamma^*TX$ 
is the tangent space  $T_{\gamma}\Omega(L_0,L_1)$ 
so $g$ defines a Riemannian metric on $\Omega(L_0,L_1)$.
The gradient vector field of $\mathcal A$ with respect 
to this metric is described as follows.
We consider an arc $\ell: (a,b) \to \Omega(L_0,L_1)_o$,
which can be identified with a map $u_{\ell} : (a,b) \times [0,1] \to X$
satisfying (path1). 
Then one can show that $\ell$ is a gradient line of $\mathcal A$
if and only if it satisfies (\ref{Floersequ}) for $H=0$, that is,
\begin{equation}\label{Floersequ2}
\frac{\partial u_{\ell}}{\partial \tau}
= J \left(
\frac{\partial u_{\ell}}{\partial t}
\right).
\end{equation}
It implies that the critical point set of $\mathcal A$ is identified with 
the intersection $L_0 \cap L_1$.
Thus a naive idea to define Lagrangian Floer homology is as follows.
We define:
\begin{equation}\label{Flchaincx}
CF(L_0,L_1;\mathbb F): = \bigoplus_{p \in L_0 \cap L_1}\mathbb F[p].
\end{equation}
For $p,q \in L_0 \cap L_1$, 
the matrix coefficient $\langle dp,q\rangle$ of the boundary operator 
$d: CF(L_0,L_1;\mathbb F) \to CF(L_0,L_1;\mathbb F)$ 
is the number (up to the shift of $\R$-direction) of solutions of the
equation (\ref{Floersequ2}) such that (path1) and 
the following two more conditions are satisfied.\footnote{More 
precisely we count only the component whose virtual dimension 
is $0$.}

\begin{enumerate}
\item[(path2)']
$\lim_{\tau\to-\infty} u(\tau,t)=p$.
\item[(path3)']
$\lim_{\tau\to+\infty} u(\tau,t)=q$.
\end{enumerate}
Floer established this theory under a rather restrictive assumption.
Let $H : X\times [0,1] \to \R$ be a smooth function.
We define $\varphi_H^t$ by (\ref{varphiH}).
A map $\varphi: X \to X$  is said to be a 
Hamiltonian diffeomorphism if $\varphi = \varphi_H^1$ for 
a certain $H$. A Hamiltonian diffeomorphism preserves the symplectic 
structure.

\begin{thm}{\rm(Floer \cite{fl2})}
Suppose that $X = T^*M$ (a cotangent bundle of a compact manifold $M$), 
$L_0 \subset X$ is the zero section, and $L_1 = \varphi(L_0)$
for a certain Hamiltonian diffeomorphism $\varphi$.
\par
Then we can define 
$d: CF(L_0,L_1;\Z_2) \to CF(L_0,L_1;\Z_2)$
such that $d\circ d =0$.
Moreover the Floer homology
$
HF(L_0,L_1;\Z_2): = \frac{{\rm Ker}d}{{\rm Im}d}
$
is isomorphic to the ordinary homology of $M$.
\end{thm}

The well-defined-ness of Floer homology can be proved in the so called 
exact case in the same way as \cite{fl2}.
Here a Lagrangian submanifold $L$ is said to be exact if 
for any $u:(D^2,\partial D^2) \to (X,L)$ the equality
$\int_{D^2}u^*\omega = 0$ holds.

Oh \cite{Oh} generalized Floer's result to monotone 
Lagrangian submanifolds as follows.
Let $L \subset X$ be a Lagrangian submanifold
and $u : (D^2,\partial D^2) \to (X,L)$  a 
continuous map.
Since $D^2$ is contractible $u$ determines a trivialization 
of $u\vert_{S^1}^*TX$. (Here $S^1 = \partial D^2$.)
For $z \in S^1$ the tangent space $T_{u(z)}L$
is a Lagrangian linear subspace of $T_{u(z)}X$.
By the trivialization $u$ determines a loop 
of the space of Lagrangian linear subspaces of a fixed 
symplectic vector space.
The fundamental group of the Lagrangian Grassmannian is known to 
be $\Z$ and so the above construction defines a map 
$\mu: \pi_2(X;L) \to \Z$, which is called the Maslov index.
Maslov index controls the (virtual) dimension of the 
moduli space of pseudo-holomorphic disks.
\par
A Lagrangian submanifold $L \subset X$ is said to be monotone 
if there exists a positive number $c$ such that
\begin{equation}\label{monotonerel}
c \int_{S^2} u^* \omega = \mu([u])
\end{equation}
for all the maps $u : (D^2,\partial D^2) \to (X,L)$.
Note that this condition is similar to (\ref{monotoneabs}).
In fact Maslov index can be regarded as a relative version 
of Chern number. Moreover if there exists a monotone Lagrangian submanifold 
in $X$ then $X$ is known to be monotone.
\par
The minimal Maslov number is the smallest positive number 
which is $\mu(\beta)$ for some $\beta \in \pi_2(X;L)$.
(If $\mu(\beta)$ is never positive and $L$ is monotone, minimal Maslov number is $\infty$ by definition.)

\begin{thm}{\rm(Oh \cite{Oh})}\label{LagFlOh}
Let
$L_0,L_1 \subset X$
be monotone Lagrangian submanifolds.
\par
Then we can define 
$d: CF(L_0,L_1;\Z_2) \to CF(L_0,L_1;\Z_2)$
such that $d\circ d =0$ in one of the following two cases: 
\begin{enumerate}
\item[(i)]
The minimal Maslov numbers of $L_0$ and $L_1$ are not strictly greater than $2$.
\item [(ii)] $L_1 = \varphi(L_0)$ for a certain Hamiltonian diffeomorphism $\varphi$,
and  the minimal Maslov numbers of $L_i$  are not smaller that $2$.
\end{enumerate}
\par
Moreover the Floer homology
$
HF(L_0,L_1;\Z_2): = \frac{{\rm Ker}d}{{\rm Im}d}
$
has the following properties.
\begin{enumerate}
\item 
If $\varphi : X \to X$ is a Hamiltonian diffeomorphism then
$$
HF(L_0,L_1;\Z_2) \cong HF(\varphi(L_0),L_1;\Z_2).
$$
\item
If $L_0 = L_1 = L$, there exists a spectral sequence whose $E_2$ page 
is $H(L;\Z_2)$ and which converges to $HF(L,L;\Z_2)$.
\end{enumerate}
\end{thm}

In the case when $L_i$ are (relatively) spin, we can work over $\Z$ coefficient 
instead of $\Z_2$ coefficient in Theorem \ref{LagFlOh}.
This fact is established in \cite[Chapter 2]{fooobook}, \cite[Chapter 8]{fooobook2}.

It was known already to Floer that beyond monotone case Floer homology 
of Lagrangian submanifolds may not be defined.
Namely $d\circ d =0$ may not hold.

The way to define and study Lagrangian Floer theory 
in the general case is established in \cite{fooobook,fooobook2}. 
A Lagrangian submanifold $L$ is said to be relatively spin 
if there exists $st \in H^2(X;\Z_2)$ which restricts to 
the second Stiefel-Whitney class of $L$. We call $st$ a background 
class and say $L$ is $st$-relatively spin if the back ground class is $st$.
We consider the universal Novikov ring $\Lambda_0$ with $R = \Q$ (or $R = \R$) as  
the ground ring  of Floer homology.
Let $\Lambda_+$ be its ideal consisting of (\ref{eleNov}) with $\lambda_i > 0$.

\begin{thm}{\rm (\cite{fooobook,fooobook2})}\label{FOOOFloer}
Let $st \in H^2(X;\Z_2)$.
For any $st$-relatively spin Lagrangian submanifold 
we can define a subset $\widetilde{\mathcal{MC}}(L) \subseteq H^{\rm odd}(L;\Lambda_+)$
with the following properties.
\begin{enumerate}
\item
There is a map $\mathcal Q:  H^{\rm odd}(L;\Lambda_+) \to 
 H^{\rm even}(L;\Lambda_+)$ of the form
$$
\mathcal Q(b) = \sum T^{\lambda_i} \mathcal Q_i(b)
$$
where $\mathcal Q_i$ is a formal power series with $\R$ coefficient 
and $\lambda_i > 0$, $\lim_{i\to \infty} \lambda_i = \infty$.
$\widetilde{\mathcal{MC}}(L)$ is the zero set of $\mathcal Q$.
The image of $\mathcal Q$ is in the kernel of the Gysin homomorphism 
$i_{!} : H^*(L) \to H^*(X)$.
\item
Let $L_i$ be $st$-relatively spin and
$b_i \in \widetilde{\mathcal{MC}}(L_i)$ for $i=0,1$. We assume that $L_0$ is 
transversal to $L_1$. Then we can define  a boundary operator
$$
d^{b_0,b_1} : CF(L_0,L_1;\Lambda_0) \to CF(L_0,L_1;\Lambda_0),
$$
where $CF(L_0,L_1;\Lambda_0)$ is as in (\ref{Flchaincx}).
It satisfies $d^{b_0,b_1} \circ d^{b_0,b_1} = 0$.
Hence Floer homology
$$
HF((L_0,b_0),(L_1,b_1);\Lambda_0): = \frac{{\rm Ker}d^{b_0,b_1}}{{\rm Im}d^{b_0,b_1}}
$$
is defined.
\item
If $\varphi : X \to X$ is a symplectic diffeomorphism 
then there exists a map $\varphi_*: H^{\rm odd}(L;\Lambda_+)
\to H^{\rm odd}(\varphi(L);\Lambda_+)$
such that
$\varphi_*(\widetilde{\mathcal{MC}}(L))
= \widetilde{\mathcal{MC}}(\varphi(L))$. Moreover 
$$
HF((\varphi(L_0),\varphi_*(b_0)),(\varphi(L_1),\varphi_*(b_1);\Lambda_0) \cong
HF((L_0,b_0),(L_1,b_1);\Lambda_0).
$$
The map $\varphi_*$ is written as
$$
\varphi_* = \varphi_{\#} + \sum_i T^{\lambda_i}  \varphi_i
$$
where $\varphi_{\#} = H^{\rm odd}(L;\R)
\to H^{\rm odd}(\varphi(L);\R)$ is a linear map induced by the diffeomorphism $\varphi$.
The map $\varphi_i$ is a formal power series\footnote{Namely 
it becomes a formal power series with $\R$ coefficient when we fix 
a basis of $H^{\rm odd}(L;\R)$ and of $H^{\rm odd}(\varphi(L);\R)$. Therefore it induces a 
map $H^{\rm odd}(L;\Lambda_+)
\to H^{\rm odd}(\varphi(L);\Lambda_+)$.} $H^{\rm odd}(L;\R)
\to H^{\rm odd}(\varphi(L);\R)$ and $\lambda_i \in \R_+$ with $\lim \lambda_i = +\infty$.
\item
If $\varphi : X \to X$ is a Hamiltonian diffeomorphism 
then 
$$
HF((\varphi(L_0),\varphi_*(b_0)),(L_1,b_1;\Lambda_0) 
\otimes \Lambda \cong
HF((L_0,b_0),(L_1,b_1);\Lambda_0) \otimes \Lambda .
$$
Here $\Lambda$ is the field of fractions of $\Lambda_0$.
\item
Suppose $L_0 = L_1 = L$, $b_0 = b_1 = b$. 
Then there exists a spectral sequence whose $E^2$ page is 
the ordinary cohomology $H(L;\Lambda_0)$ and which 
converges to $HF((L,b),(L,b);\Lambda_0)$.
The image of the differential is contained in the 
subquotient of the kernel of the Gysin homomorphism 
$i_{!} : H^*(L) \to H^*(X)$.\footnote{In particular if 
$i_{!} : H^*(L) \to H^*(X)$ is injective then 
$E_2 \cong E_{\infty}$.}
\end{enumerate}
\end{thm}
\begin{rem}\label{rem44}
We may consider $b \in H^{\rm odd}(L;\Lambda_0)$
actually. Since $\mathcal Q_i(b)$ does not make 
sense in this generality for a formal power series $\mathcal Q_i$
we need to state it in a bit more careful way as follows.
We consider
$$
\frac{H^1(L;\Lambda_0)}{2\pi i H^1(L;\Z)} 
\times
\prod_{k>0} H^{2k+1}(L;\Lambda_0).
$$
Taking a basis of the free parts of $H^1(L;\Z)$ and of
$H^{2k+1}(L;\Z)$, its element is written by coordinates 
$y^1_j = \exp(x^1_j)$ and $x^{2k+1}_j$.
Here $x^1_j$ and $x^{2k+1}_j$ are 
coordinates of $H^{\rm odd}(L;\Lambda_0)$ 
corresponding to the $j$-th basis of $H^1(L;\Z)$ and of 
$H^{2k+1}(L;\Z)$, respectively.
Then we can write $Q$ as
$$
Q(b) = \sum_i T^{\lambda_i} Q_i((y^1_j),(x^{2k+1}_j)).
$$
Here $Q_i$ is a {\it polynomial} of $y^1_j = \exp(x^1_j)$ and $x^{2k+1}_j$
and $\lim_{i\to \infty} \lambda_i = +\infty$.\footnote{We can prove this fact 
by using a `disk analogue of divisor axiom'. See \cite[Lemma 13.1]{fukaya:cyc} and 
\cite{Yu}.}
Therefore, the equation $Q(b) = 0$ makes sense for 
$b \in \frac{H^1(L;\Lambda_0)}{2\pi i H^1(L;\Z)} 
\times
\prod_{k>0} H^{2k+1}(L;\Lambda_0)$.
\end{rem}
An element of $\widetilde{\mathcal{MC}}(L)$ is called 
a bounding cochain or a Maurer-Cartan element of $L$ and plays an important 
role. The equation $\mathcal Q(b) = 0$ is actually a Maurer-Cartan equation.
See the next section.
\par
At first sight the existence of an extra parameter $b$ is unsatisfactory.
However actually it expands the applicability of  Lagrangian Floer 
theory.
In fact the Floer homology $HF((L,b),(L,b);\Lambda)$ becomes trivial 
very frequently and, in various examples, it becomes non-zero only at a 
very special value of $b$. So this extra freedom allows us wider possibility 
to obtain a non-trivial Floer homology.
(See for example \cite{toric1}.)
\par
Theorem \ref{FOOOFloer} is generalized by Akaho-Joyce \cite{AJ} 
to the case of immersed Lagrangian submanifolds as follows.
Let $L = (\tilde L,i_L)$ be an immersed Lagrangian submanifold of $X$.
Namely $\tilde L$ is an $n = \dim X/2$ dimensional closed manifold 
and $i_L : \tilde L \to X$ is an immersion such that $i_L^*\omega = 0$.
We assume that self-intersection is transversal.
So 
$$
{\rm SW}(L): = \{(p,q) \in \tilde L\times \tilde L \mid p\ne q ,\,\,\, i_L(p) = i_L(q)\}
$$
is a finite set. (Note that $\#SW(L)$ is twice of the number of self-intersections.)
\begin{equation}\label{AJXCF}
CF(L;R) = H(L;R) \oplus \bigoplus_{(p,q) \in {\rm SW}(L)} R[p,q].
\end{equation}
We can define the degree $d(p,q)$ of $(p,q) \in {\rm SW}(L)$
such that $d(p,q) + d(q,p) = n$.

\begin{thm}{\rm (\cite{AJ})}\label{AJFloer}
Theorem \ref{FOOOFloer} (1)(2)(3)(4) hold for immersed Lagrangian submanifolds $L_i$ 
when we replace $H(L;R)$ by $CF(L;R)$.
\end{thm}

\section{$A_{\infty}$ algebras and $A_{\infty}$ categories.}
\label{sec:Ainfinity}

The language of $A_{\infty}$ algebras and categories is inevitable 
to understand  Lagrangian Floer theory appearing in Theorems \ref{FOOOFloer},\ref{AJFloer}.
\par
\begin{figure}[h]
\centering
\includegraphics[scale=0.6]{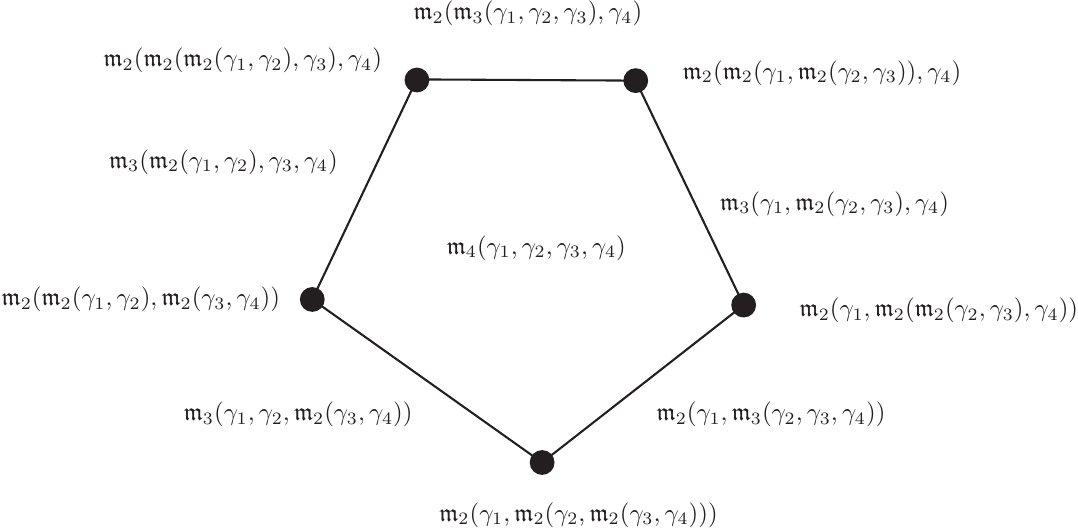}
\caption{Stasheff 2-gon}
\label{Figure1}
\end{figure}

The notions of  $A_{\infty}$ spaces and algebras were invented by 
Stasheff in his study \cite{St} of loop spaces.
The product operation of loops $\gamma: S^1 \to X$ is not
strictly associative because of the parametrization problem.
(Here the product $\gamma_1 \circ \gamma_2$ of loops is defined by
$$
(\gamma_1 \circ \gamma_2)(t)
= 
\begin{cases}
\gamma_1(2t)  &\text{if $t \in [0,1/2]$} \\
\gamma_2(2t-1)  &\text{if $t \in [1/2,1]$}.)
\end{cases}
$$
However there is a canonical homotopy 
between two compositions $(\gamma_1 \circ \gamma_2)\circ \gamma_3$
and $\gamma_1 \circ (\gamma_2 \circ \gamma_3)$.
Namely there is a $[0,1]$ parametrized family of loops 
$\frak m_3(\gamma_1,\gamma_2,\gamma_3): S^1 \times [0,1] \to X$
such that 
$$
\frak m_3(\gamma_1,\gamma_2,\gamma_3)(t;0)
= ((\gamma_1 \circ \gamma_2)\circ \gamma_3)(t),
\quad
\frak m_3(\gamma_1,\gamma_2,\gamma_3)(t;1)
= (\gamma_1 \circ (\gamma_2 \circ \gamma_3))(t).
$$
The `homotopy associativity' of the product of loops is 
actually stronger. There is a two parameter family 
of loops $\frak m_4(\gamma_1,\gamma_2,\gamma_3,\gamma_4)$
which bounds 
the union of
$$
\aligned
&\frak m_3(\frak m_2(\gamma_1,\gamma_2),\gamma_3, \gamma_4),
\quad
\frak m_2(\frak m_3(\gamma_1,\gamma_2,\gamma_3),\gamma_4),
\quad
\frak m_3(\gamma_1,\frak m_2(\gamma_2,\gamma_3), \gamma_4) \\
&\frak m_2(\gamma_1,\frak m_3(\gamma_2,\gamma_3, \gamma_4)),
\quad 
\frak m_3(\gamma_1,\gamma_2,\frak m_2(\gamma_3, \gamma_4)).
\endaligned
$$
Here we write $\frak m_2(\gamma,\gamma')$ in place of 
$\gamma \circ \gamma'$.
See Figure \ref{Figure1}.
We can continue and obtain a $(k-2)$-parameter family of 
loops $\frak m_k(\gamma_1,\dots,\gamma_k)$ 
whose boundary parametrizes the union of the compositions of 
$\frak m_{k_1}$ and $\frak m_{k_2}$ with $k_1 + k_2 = k+1$.
This is the definition of an $A_{\infty}$ space.
Its algebraic analogue is an $A_{\infty}$ algebra, 
which is defined as follows.
\begin{defn}{\rm (Stasheff)}\label{AinfSta}
An $A_{\infty}$ algebra is a graded $\mathbb F$ module $C$ 
together with operations 
$$
\frak m_k: \underbrace{C \otimes \dots \otimes C}_k \to C
$$
of degree $2-k$ for $k=1,2,\dots$ which satisfies the following $A_{\infty}$ relation.
\begin{equation}\label{Ainfintyform}
0 = \sum_{k_1+k_2 = k+1}\sum_{i=1}^{k+1-k_2}
(-1)^* \frak m_{k_1}(x_1,\dots,\frak m_{k_2}(x_i,\dots,x_{i+k_2-1}) \dots, x_k).
\end{equation}
Here $* = \deg x_1 + \dots + \deg x_{i-1} + i-1$.\footnote{This is the sign 
convention of \cite{fooobook} which is slightly different from \cite{St}.}
\end{defn}

In  Lagrangian Floer theory, we use curved and filtered $A_{\infty}$ 
algebra (and category) which is different from Definition \ref{AinfSta}
at the following points.
\begin{enumerate}
\item The operation $\frak m_0: \mathbb F \to C$ can be non-zero. 
\item The coefficient ring $\mathbb F$ is the universal Novikov ring $\Lambda_0$
and the operations are assumed to preserve the filtration.
Moreover $\frak m_0(1)$ is assumed to be in $C \otimes_{\Lambda_0} \Lambda_+$.
\end{enumerate}
Let us elaborate on those points.
In the case when $\frak m_0 = 0$, (\ref{Ainfintyform}) implies 
$\frak m_1\circ \frak m_1 = 0$. 
However in the case when $\frak m_0 \ne 0$ it implies
$$
(\frak m_1\circ \frak m_1)(x) 
= (-1)^{\deg x}\frak m_2(x,\frak m_0(1)) - \frak m_2(\frak m_0(1),x).
$$
Therefore $\frak m_0$ is an obstruction for Lagrangian Floer homology 
to be well-defined.
\par
The universal Novikov ring has a filtration $\frak F^{\lambda}\Lambda_0$
which consists of elements 
$\sum_{i=0}^{\infty} a_i T^{\lambda_i}$ such that 
$\lambda_i \ge \lambda$ for all $i$ with $a_i \ne 0$.
The module $C$ over $\Lambda_0$ is assumed to be a 
completion of the free $\Lambda_0$ module
$\overline C \otimes_R \Lambda_0$.
In other words it consists of (infinite) sums 
$\sum_{i=0}^{\infty} c_i T^{\lambda_i}$ 
with $\lim_{i\to\infty}\lambda_i = +\infty$, where $c_i \in \overline C$
and $\overline C$ is a free $R$ module.
Then we can define a filtration on $C$ in a similar way as the filtration on $\Lambda_0$.
We call such $C$ a completed free filtered $\Lambda_0$ module.
We then require the operations to preserve the filtration.
We call such $(C,\{\frak m_k\})$ a (curved) filtered $A_{\infty}$
algebra.
The topology induced by the filtration is called the $T$-adic topology.

The relation of an $A_{\infty}$ structure to 
a perturbation of the structure and to mathematical physics 
had been known before early 1990's.
(See for example \cite{GS}.)
We can use it in Lagrangian Floer theory as follows. 
Let $(C,\{\frak m_k\})$ be a (curved) filtered $A_{\infty}$
algebra.
For $b \in C^{\rm odd}$ with $b \in \frak F^{\lambda}C$, $\lambda > 0$, 
we consider the Maurer-Cartan equation:
\begin{equation}\label{MCequation}
\sum_{k=0}^{\infty} \frak m_k(b,\dots,b) = 0.
\end{equation}
Note that the term for $k=0$ of left hand side is $\frak m_0(1)$.
Since $\frak m_k(b,\dots,b) \in \frak F^{k\lambda}C$
the left hand side converges in $T$-adic topology.
\par
Suppose $b$ solves (\ref{MCequation}). We define the 
operations $\frak m_k^b$ by the following formula:
\begin{equation}\label{defbyb}
\frak m_k^b(x_1,\dots,x_k)
= 
\sum_{\ell_0,\dots,\ell_k = 0}^{\infty}
\frak m_{k+\sum \ell_i}(b^{\ell_0},x_1,b^{\ell_1},\dots,b^{\ell_{k-1}},x_k,b^{\ell_k}).
\end{equation}
The right hand side converges in $T$-adic topology.
It is easy to show that $\{\frak m_k^b\}$ satisfies the
$A_{\infty}$-relation (\ref{Ainfintyform}).
Moreover the Maurer-Cartan equation (\ref{MCequation})
implies $\frak m_0^b = 0$.
In particular we have
$$
\frak m_1^b \circ \frak m_1^b = 0.
$$
Thus we can eliminate the `curvature' $\frak m_0$ by 
using the solution of  Maurer-Cartan equation (\ref{MCequation}).
\par
Filtered $A_{\infty}$ algebras appear in Lagrangian Floer theory 
as follows.

\begin{thm}{\rm (\cite{fooobook,fooobook2,AJ})}\label{LagfilAINF}
Let $L$ be a relatively spin Lagrangian submanifold of a symplectic 
manifold. Then we can associate a structure of a filtered $A_{\infty}$
algebra to $H(L;\Lambda_0)$.
\par
The same holds for an immersed Lagrangian submanifold if we replace  $H(L;\Lambda_0)$ by 
$C(L;\Lambda_0)$.
\end{thm}
We can define a map $H^{\rm odd}(L;\Lambda_+) 
\to H^{\rm ev}(L;\Lambda_+)$  by
$$
b \mapsto \sum_{k=0}^{\infty} \frak m_k(b,\dots,b).
$$
This is the map $\mathcal Q$ in Theorem \ref{FOOOFloer}.
\par
Bondal-Kapranov \cite{bondkap} defined the notion of a DG-category.
We can modify it to define a filtered $A_{\infty}$ category as follows.

\begin{defn}
A curved filtered $A_{\infty}$ category $\mathscr C$ 
is the following:
\begin{enumerate}
\item The set of objects $\frak{OB}(\mathscr C)$ is given.
\item For $c,c' \in \frak{OB}(\mathscr C)$ a graded completed 
free $\Lambda_0$ module $\mathscr C(c,c')$ is given. 
This is the set of morphisms.
\item We put
$$
B_k\mathscr C(c,c') 
= \bigoplus_{c_0,\dots,c_k} \bigotimes_{i=1}^k \mathscr C(c_{i-1},c_i),
$$
where the direct sum is taken over $c_0,\dots,c_k \in \frak{OB}(\mathscr C)$
such that $c_0 = c$, $c_k = c'$.
We also put
$$
B_0\mathscr C(c,c')  = 
\begin{cases}
0   & c\ne c', \\
\Lambda_0   & c = c'.
\end{cases}
$$
The $\Lambda_0$ module homomorphisms
$$
\frak m_k : B_k\mathscr C(c,c') \to \mathscr C(c,c')
$$
are given for $k=0,1,2,\dots$. It is called the structure operations.
The structure operations are required to preserve filtrations.
\item
The $A_{\infty}$ relation (\ref{Ainfintyform}) is satisfied.
\end{enumerate}
\end{defn}
A filtered $A_{\infty}$ category is said to be strict if 
$\frak m_0 = 0$. 
To a curved filtered $A_{\infty}$ category $\mathscr C$ we can associate 
a strict filtered $A_{\infty}$ category $\mathscr C_{s}$ as follows.
We remark that for $c \in \frak{OB}(\mathscr C)$
the restriction of $\frak m_k$ defines a structure of 
a curved filtered $A_{\infty}$ algebra on $\mathscr C(c,c)$.
A bounding cochain of $c$ is by definition an element $b$ 
of $\mathscr C^{\rm odd}(c,c)$ such that $b \in \frak F^{\lambda}\mathscr C(c,c)$
for $\lambda > 0$ and $b$ satisfies (\ref{MCequation}).
\par
An object of $\mathscr C_{s}$ is a pair $(c,b)$ where 
$c \in \frak{OB}(\mathscr C)$ and $b$ is its bounding cochain.
The module of morphisms 
$\mathscr C_{s}((c,b),(c',b'))$ is $\mathscr C(c,c')$.
The structure operations of $\mathscr C_{s}$ is defined 
by modifying the structure operations of $\mathscr C$ 
in the same way as (\ref{defbyb}).
\par
A strict $A_{\infty}$ category such that $\frak m_k = 0$
for $k\ne 1,2$ is called a DG-category.
\par
We can define a notion of a (filtered) $A_{\infty}$ functor 
between two (filtered) $A_{\infty}$ categories.
(See \cite[Section 7]{fu4}.)
\par
In Lagrangian Floer theory a filtered $A_{\infty}$ category 
appears in the following way.
Let $(X,\omega)$ be a symplectic manifold.
We fix a background class $st \in H^2(X;\Z_2)$.
We consider a finite set $\mathbb L$ of $st$-relatively 
spin (immersed) Lagrangian submanifolds such that
for $L,L' \in \mathbb L$, $L$ is transversal to $L'$.
\begin{thm}\label{LagAinfcat}
There exists a curved filtered $A_{\infty}$ category 
$\frak{Fuk}(X;\mathbb L)$ the set of whose objects is $\mathbb L$.
For $L \in \mathbb L$ the curved filtered $A_{\infty}$ algebra
$\frak{Fuk}(X;\mathbb L)(L,L)$ is one in Theorem \ref{LagfilAINF}.
\end{thm}
We can prove Theorem \ref{LagAinfcat} from Theorem \ref{LagfilAINF} as follows.
We consider the disjoint union of the elements of $\mathbb L$
and regard it as a single immersed Lagrangian submanifold $\hat L$.
We apply Theorem \ref{LagfilAINF} (Akaho-Joyce's immersed case) to $\hat L$
to obtain a curved filtered $A_{\infty}$ algebra. 
The structure operations of this curved filtered $A_{\infty}$ algebra 
induce structure operations of $\frak{Fuk}(X;\mathbb L)$.
See also \cite{AFOOO,cor}.
We denote by $\frak{Fukst}(X;\mathbb L)$ the strict category 
associated to $\frak{Fuk}(X;\mathbb L)$.
Its object is a pair $(L,b)$ where $L \in \mathbb L$ and $b$ is its 
bounding cochain.
\par
Let $(L_i,b_i)$ be an object of $\frak{Fukst}(X;\mathbb L)$ for $i=0,1$.
The $\frak m_1$ operator of $\frak{Fukst}(X;\mathbb L)$
on $\frak{Fukst}(X;\mathbb L)((L_0,b_0),(L_1,b_1))$ is 
by definition
$$
x \mapsto \sum_{k,\ell=0}^{\infty} \frak m_{k+\ell+1}(b_0^k,x,b_1^{\ell}).
$$
(See (\ref{defbyb}).) Here $\frak m_{k+\ell+1}$ in the right hand side 
is the structure operation of $\frak{Fuk}(X;\mathbb L)$.
This is the boundary operator $d^{b_0,b_1}$ in Theorem \ref{FOOOFloer}.
\par
As was mentioned at the end of  Section \ref{HFtoposec}, we need a bit more general 
object than those of $\frak{Fukst}(X;\mathbb L)$ for various purposes.
One is establishing topological field theory picture in Donaldson-Floer 
theory. The others are for homological mirror symmetry 
(Section \ref{secHMS}) and for Lagrangian 
correspondence (Section \ref{sec:COR+GAUGE}).
One way to do so is to use the notion of an $A_{\infty}$ module over 
an $A_{\infty}$ category\footnote{In \cite{fu2} etc. the author used 
an $A_{\infty}$ functor from an $A_{\infty}$ category $\mathscr C$ to the 
DG category $\mathcal{CH}$ whose objects are chain complexes. 
These two notions are the same as explained in \cite[Subsection 5.1]{cor}.} and 
to use Yoneda embedding.

\begin{defn}
Let $\mathscr C$ be a filtered $A_{\infty}$ category.
A right $A_{\infty}$ module $\mathscr D$ over $\mathscr C$ associates 
a chain complex $\mathscr D(c)$ to each $c \in \frak{OB}(\mathscr C)$
and a map
\begin{equation}\label{fraknismap}
\frak n_k: \mathscr D(c) \otimes_{\Lambda_0} B_k\mathscr C(c,c')
\to \mathscr D(c')
\end{equation}
for $k=0,1,\dots$ and $c,c' \in \frak{OB}(\mathscr C)$ with 
the following properties.
\begin{enumerate}
\item
$\frak n_0$ is the boundary operator of $\mathscr D(c)$.
\item
The following $A_{\infty}$ relation is satisfied.
\begin{equation}\label{rmodAinf}
\aligned
0=&\sum_{k_1+k_2 = k} \frak n_{k_2}(\frak n_{k_1}(y;x_1,\dots,x_{k_1});
x_{x_1+1},\dots,x_{k}) \\
&+ 
\sum_{k_1+k_2=k+1}\sum_{i=1}^{k_1} (-1)^*
\frak n_{k_1}(y;x_1,\dots,\frak m_{k_2}(x_i,\dots,x_{i+k_2-1}),\dots,x_k)
\endaligned
\end{equation}
where $* = \deg y + \deg x_1 + \dots + \deg x_{i-1} + i$.\footnote{
There is an alternative sign convention where $+i$ is replaced by $+i-1$.
See \cite[Subsection 5.1]{cor}.}
\end{enumerate}
\end{defn}
We can define the notion of a right $A_{\infty}$ module homomorphism
between two right $A_{\infty}$ modules and obtain 
a DG-category whose object is a right  $A_{\infty}$ module
over $\mathscr C$. We denote it by $\frak{RMOD}(\mathscr C)$.
\par
We can define the notion of an $A_{\infty}$ functor. 
There exists an $A_{\infty}$ functor from 
$\mathscr C$ to $\frak{RMOD}(\mathscr C)$ which is an 
$A_{\infty}$ analogue of the Yoneda embedding, as follows.
\par
Let $c \in \frak{OB}(\mathscr C)$.
We associate a right $A_{\infty}$ module $\mathscr D_c$ 
by
$$
\mathscr D_c(c') = \mathscr C(c,c')
$$
and
$$
\frak n_k(y;x_1,\dots,x_k) = \frak m_{k+1}(y,x_1,\dots,x_k)
$$
Here the left hand side is the right module structure 
and the right hand side is the structure operation of $\mathscr C$.
The equality (\ref{rmodAinf}) is a consequence of (\ref{Ainfintyform}).

This is the way how the objects are sent by the Yoneda functor.
We can define the morphism part by using the structure operations 
$\frak m$ of $\mathscr C$. See \cite[Section 9]{fu4}.
\par
An $A_{\infty}$ analogue of Yoneda's lemma is Theorem \ref{Yoneda}.
\begin{defn}
A strict unit of an $A_{\infty}$ category $\mathscr C$ assigns 
${\bf e}_c \in \mathscr C(c,c)$ (of degree $0$) to each object $c$ such that:
\begin{enumerate}
\item $\frak m_k(\dots,{\bf e}_c,\dots) = 0$ unless $k=2$.
\item $\frak m_2({\bf e}_c,x) = (-1)^* \frak m_2(x,{\bf e}_c) = x$,
where $* = \deg x$.
\end{enumerate}
An $A_{\infty}$ category is said to be unital if it has a strict unit.
\end{defn}
\begin{thm}\label{Yoneda}
If $\mathscr C$ is a strict and unital filtered $A_{\infty}$ category, 
then the Yoneda embedding 
$\mathscr C \to \frak{RMOD}(\mathscr C)$ is a homotopy 
equivalence to a full subcategory.
\end{thm}
See \cite{fu4,cor} for the homotopy equivalence of filtered $A_{\infty}$
category. See \cite[Section 9]{fu4} for the proof of Theorem \ref{Yoneda}.
\par
By Theorem \ref{Yoneda} we may regard an object of $\mathscr C$ as 
a right $\mathscr C$ module. In other words a right 
$\mathscr C$ module is regarded as an `extended' object of $\mathscr C$.
Based on this observation the following is proposed in \cite{fu1,fu2}.

\begin{conj}{\rm (\cite{fu1,fu2})}\label{conjk222}
Let $M$ be a 3-dimensional manifold with boundary $\Sigma$ and 
$\mathcal E_M$  an $SO(3)$ bundle such that 
 $w_2(\mathcal E_{\Sigma}) = [\Sigma]$.
(Here $\mathcal E_{\Sigma}$ is the restriction of $\mathcal E_M$
to $\Sigma$.)
\par
Then we can associate a right filtered $A_{\infty}$ module 
$HF_M$ on  $\frak{Fukst}(R(\Sigma;\mathcal E_{\Sigma});\mathbb L)$
for any finite set $\mathbb L$ of Lagrangian submanifolds of 
$R(\Sigma;\mathcal E_{\Sigma})$.
\par
Furthermore the following holds.
Let $M_i$ be a 3-dimensional manifold for $i=1,2$ with $
\partial M_i = \Sigma$ and 
$\mathcal E_{M_i}$  an $SO(3)$ bundle such that 
$w_2(\mathcal E_{\Sigma}) = [\Sigma]$ where 
$\mathcal E_{\Sigma}$ is the restriction of $\mathcal E_{M_i}$
to $\Sigma$. 
Let $M$ be a closed $3$-dimensional manifold obtained 
by gluing $-M_1$ and $M_2$ along $\Sigma$.
An $SO(3)$ bundle $\mathcal E_M$ on $M$ is obtained by 
gluing $\mathcal E_{M_1}$ and $\mathcal E_{M_2}$.
Then we can choose $\mathbb L$ such that the isomorphism
$$
I(M;\mathcal E_M) \cong H(Hom(HF_{M_1},HF_{M_2}))
$$
holds. Here the left hand side is the instanton Floer homology 
and the right hand side is the cohomology of the morphism 
complex in the DG-category of 
right filtered $A_{\infty}$ modules.
\end{conj}
Together with A. Daemi and M.  Lipyanskiy 
the author is on the way of proving this conjecture.
(We will discuss it more in Section \ref{sec:COR+GAUGE}.)
\par
Note that the right $A_{\infty}$ module $HF_M$ is 
supposed to associate a cohomology
to an object $(L,b)$  of $\frak{Fukst}(R(\Sigma;\mathcal E_{\Sigma});\mathbb L)$,
that is, a pair of a Lagrangian submanifold $L \in \mathbb L$ 
of $R(\Sigma;\mathcal E_{\Sigma})$ and its bounding cochain $b$.
Let us write it as $HR(M;(L,b))$.
An idea in \cite{fu1,fu2} for the construction of such cohomology 
theory $HR(M;(L,b))$ is to study gauge theory of 
$M \times \R$, which has a boundary $\Sigma \times \R$ and 
ends $M \times \{\pm \infty\}$,
use $L$ to set an appropriate boundary condition 
on  $\Sigma \times \R$, and try to imitate the construction 
of instanton Floer homology, which is the case $\partial M = \emptyset$.
\par
This part of the idea is later realized in the case when $L$ is monotone and 
$b =0$ by Salamon-Wehrheim in \cite{Sawe}.
Namely they construct such a Floer homology $HR(M;L)$ of $M$ with 
coefficient $L$, a monotone Lagrangian submanifold of 
$R(\Sigma;\mathcal E_{\Sigma})$.
\par
To prove the first part of Conjecture \ref{conjk222} 
we also need to define a right module structure, 
the structure map (\ref{fraknismap}).
See \cite{takagi} on this point.

\section{Homological Mirror symmetry.}
\label{secHMS}

Kontsevich \cite{Ko} used $A_{\infty}$ categories to formulate 
his famous homological mirror symmetry conjecture.
Actually he used derived category of $A_{\infty}$ categories 
which we review below.
\par
Let $\mathscr C$ be a strict (filtered) $A_{\infty}$ category.
We assume for simplicity that for an object $c$ of $\mathscr C$
its degree shift\footnote{See the beginning of 
Subsection \ref{subsec:gene2}.} $c[k]$ is also exists as an object of  $\mathscr C$.
We consider a finite sequence $c_1,\dots,c_n$ of objects of 
$\mathscr C$ and let $x_{ij} \in \mathscr C(c_i,c_j)$
for $i< j$.
We call $\mathcal C = (\{c_i\},\{x_{ij}\})$ a twisted complex 
if for all $i<j$ the next equation is satisfied.
\begin{equation}\label{Twistedcompx}
\sum \frak m_k(x_{\ell_0\ell_1},\dots,x_{\ell_{k-1}\ell_k}) = 0
\end{equation}
where the sum is taken over the set of all
$\ell_0,\dots,\ell_k$ with $k=1,2,\dots$, $\ell_0 = i$, $\ell_{k} = j$, 
$\ell_0 < \dots <\ell_k$.
Generalizing a similar notion in the case of DG category
(which is due to Bondal-Kapranov) 
Kontsevich introduced the notion of a twisted
complex and showed that there is a triangulated 
category $\mathbb D(\mathscr C)$ whose object is a 
twisted complex.
A triangulated category is an additive category
(that is, the category the set of whose morphisms forms an abelian 
group) so that for each morphism $f: c \to c'$ 
there is a mapping cone which satisfies  
appropriate axioms. (See \cite{Ha}.)
\par
A morphism between two twisted complexes  
$\mathcal C^{(m)} = (\{c_i^{(m)}\},\{x^{(m)}_{ij}\})$
$m=1,2$ 
is a tuple $(y_{ij})$ such that $y_{ij} \in \mathscr C(c_i^{(1)},c_j^{(2)})$.
The differential
is defined by $d(y_{ij}) = (z_{ij})$
where
$$
z_{ij} = \sum \frak m_{k+n+1}(x^{(1)}_{\ell_0\ell_1},\dots,x^{(1)}_{\ell_{k-1}\ell_k},y_{\ell_k m_0},
x^{(2)}_{m_0m_1},\dots,x^{(2)}_{m_{n-1}m_n}).
$$
Here the sum is taken over all $\ell_0,\dots,\ell_k$, 
$m_0,\dots,m_n$ with $\ell_0=i$, $m_n = j$.
\par
The equation (\ref{Twistedcompx}) and the $A_{\infty}$ relation 
imply that $d\circ d = 0$.
A closed morphism from $\mathcal C^{(1)}$ to $\mathcal C^{(2)}$
is an element $(y_{ij})$ such that $d(y_{ij}) = 0$.
For a closed morphism we can define its mapping cone. See \cite[Section 6]{fu4}.
We can then define a triangulated category $\mathbb D(\mathscr C)$
as the localization of this category 
so that a closed morphism which induces an 
isomorphism on cohomologies becomes an isomorphism.
\par
For a complex manifold $X$ we can define an abelian category 
$\mathscr{SH}(X)$ of its coherent 
sheaves.  Its derived category $\mathbb D(\mathscr{SH}(X))$
is defined roughly as follows. (See \cite{Ha}.)
We consider the additive category whose object is a 
chain complex of coherent sheaves. 
We can define mapping cone of such chain complex and 
also define the notion of a weak equivalence and take the
localization.
We thus obtain $\mathbb D(\mathscr{SH}(X))$.
\par
The homological mirror symmetry conjecture by 
Kontsevich is stated as:
\begin{equation}
\mathbb D(\frak{Fukst}(X)) 
\cong 
\mathbb D(\mathscr{SH}(X^{\vee})),
\end{equation}
where $X^{\vee}$ is a mirror of $X$.
\par
We can formulate the homological mirror symmetry conjecture
using the terminology of filtered $A_{\infty}$
category rather than using that of derived category 
as follows.
Note that the filtered $A_{\infty}$ category 
$\frak{Fukst}(X)$ is linear over the 
universal Novikov ring $\Lambda_0$, 
which is similar to the formal 
power series ring $\C[[T]]$.
The ring $\C[[T]]$ is regarded as `the 
set of functions' of a formal neighborhood 
of $0$ in $\C$.
We regard the mirror of $X$ as a 
scheme\footnote{or a formal scheme or a rigid analytic space, etc.} over $\C[[T]]$, that is a 
formal deformation $X^{\vee} \to {\rm Sp}(\C[[T]])$.
We require that it is a formalization of a 
maximal degenerate family 
$\pi: X^{\vee}_{\epsilon} \to D^2(\epsilon)$ of, 
say, Calabi-Yau manifolds. 
Namely we require:
\begin{enumerate}
\item[(md1)]
For $q \in D^2(\epsilon)\setminus \{0\}$
the fiber $M_q: = \pi^{-1}(q)$ is a Calabi-Yau 
manifold. (A K\"ahler manifold whose 
canonical bundle is trivial.)
\item[(md2)]
The fiber $M_0: = \pi^{-1}(0)$ of $0$ is a normal 
crossing divisor.
\item[(md3)]
There exists a point $p$ in $M_0$ where 
the intersection of a neighborhood $U$ of $p$ 
in $X^{\vee}$ with $M_0$ is bi-holomorphic 
to  
$$
\{(z_0,\dots,z_{n+1}) \mid z_i \in D^2, z_0\dots z_{n+1} = 0\}.
$$
\end{enumerate}
\begin{rem}
It seems that the third condition is equivalent to the following 
condition. 
\begin{enumerate}
\item[(*)] Let $q_i \in D^2\setminus\{0\}$ be a sequence converging to $0$ 
and we use a metric on $M_{q_i} = \pi^{-1}(q)$ which is Ricci-flat 
and has diameter $1$. Then 
the topological dimension of the Gromov-Hausdorff limit 
$\lim_{i\to \infty} M_{q_i}$
is $n$, the complex dimension of $M_{q_i}$.
\end{enumerate}
In fact when (md3) is satisfied, it is expected that the Riemannian manifold $M_{q_i}$ 
is mostly occupied in a small neighborhood of the points  $p$ 
in $X^{\vee}$ as in (md3). (Namely the complement of such a neighborhood 
has small diameter.) 
\par
On the other hand, other than the case of dimension 2 (\cite{GW}) or tori,
the relation between (*) and (md3) is not established yet.
(See \cite{LiY} for a related recent result.)
Gross-Wilson \cite{GW} and 
Kontsevich-Soibelman \cite{KS} conjectured that condition (*) 
is equivalent to a definition of the maximal degenerate point 
in complex geometry, that is,  the monodromy $\rho : H_N(M_q) \to H_N(M_q)$
satisfies $(\rho - {\rm id})^{n-1} \ne 0$ and $(\rho - {\rm id})^{n} = 0$.
The author does not know the present status of this conjecture.
\end{rem}
\par
We consider the $n$-fold (branched) covering 
${\rm Sp}(\C[[T^{1/n}]]) \to {\rm Sp}(\C[[T]])$,
which induces $X^{\vee}_n \to {\rm Sp}(\C[[T^{1/n}]])$ via pull back.
We consider the category of coherent sheaves on 
$\mathscr{SH}(X^{\vee}_n)$
which is $\C[[T^{1/n}]]$ linear. 
Take the inductive limit
$$
\mathscr{SH}(X^{\vee}_{\infty}) := \varinjlim \mathscr{SH}(X^{\vee}_n)
$$
which is linear over
$$
\Lambda_{0,\Q}^{\C}: =
\left\{\sum a_i T^{\lambda_i} \in \Lambda_{0}^{\C} \mid \lambda_i \in \Q_{\ge 0}\right\}.
$$
We then obtain a DG-category (linear over $\Lambda_{0,\Q}^{\C}$) 
whose object is a 
chain complex of objects of $\mathscr{SH}(X^{\vee}_{\infty})$.
We denote it by 
${\mathbb{DG}}(\mathscr{SH}(X^{\vee}_{\infty}))$.
\par
In the symplectic side we modify 
$\frak{Fukst}(X)$ slightly so that it becomes 
$\Lambda_{0,\Q}^{\C}$ linear in place of being
$\Lambda_{0}^{\C}$ linear as follows.
\begin{defn}
Let $(X,\omega)$ be a symplectic manifold.
We assume that the cohomology class $[\omega]$ is in $H^2(X;\Z)$ and 
take a complex line bundle $\frak L$ 
such that $c_1(\frak L) = [\omega]$.
We take a Hermitian connection $\nabla$ of 
$\frak L$ such that $F_{\nabla} = 2\pi i\, \omega$.
A pair $(\frak L,\nabla)$ is called a 
pre-quantum bundle.
\par
Let $L$ be a Lagrangian submanifold of $X$. 
Since $\omega\vert_L =0$, the restriction of 
$(\frak L,\nabla)$ to $L$ is flat.
Let ${\rm Hol}_{\nabla} : \pi_1(L) \to U(1)$
be the holonomy representation of this flat bundle.
We say $L$ is rational 
if the image of ${\rm Hol}_{\nabla}$ is a finite group.
\end{defn}
By modifying the construction of $\frak{Fukst}(X)$
we can define a filtered $A_{\infty}$ category 
$\frak{Fukst}_{\Q}(X)$ such that:

\begin{enumerate}
\item 
An object of $\frak{Fukst}_{\Q}(X)$ is a pair 
of a rational Lagrangian submanifold $L$ 
of $X$ and its bounding cochain $b$ defined over 
$\Lambda_{0,\Q}^{\C}$ (see the explanation below).
\item
$\frak{Fukst}_{\Q}(X)$ is $\Lambda_{0,\Q}^{\C}$ linear.
\end{enumerate}
Suppose that $L$ is a  rational Lagrangian submanifold.
We consider the filtered $A_{\infty}$ algebra
$CF(L,L)$. The exponent $\lambda$  in the weight $T^{\lambda}$
appearing in the structure operations of the $A_{\infty}$ structure of $CF(L,L)$ 
is the symplectic area of a disk $(D^2,\partial D^2) \to (X,L)$.
Using the rationality it is easy to see that 
$\lambda \in \Q_{>0}$. Therefore the Maurer-Cartan 
equation (\ref{MCequation}) is defined over $\Lambda_{0,\Q}^{\C}$.
Thus the notion of bounding cochain $b$ defined over 
$\Lambda_{0,\Q}^{\C}$ makes sense.
\par
We can then modify the definition of $A_{\infty}$
operations slightly so that (2) is satisfied.
See \cite[Proposition 2.2]{fu5}.

\begin{conj}
There exists a filtered $A_{\infty}$ functor
$\frak{Fukst}_{\Q}(X) \to {\mathbb{DG}}(\mathscr{SH}(X^{\vee}_{\infty}))$
which induces an isomorphism of their derived categories 
in certain cases.
\end{conj}

\subsection{Elliptic curves and tori}
\label{subsec:ellipticcurve}

\begin{figure}[h]
\centering
\includegraphics[scale=0.3]{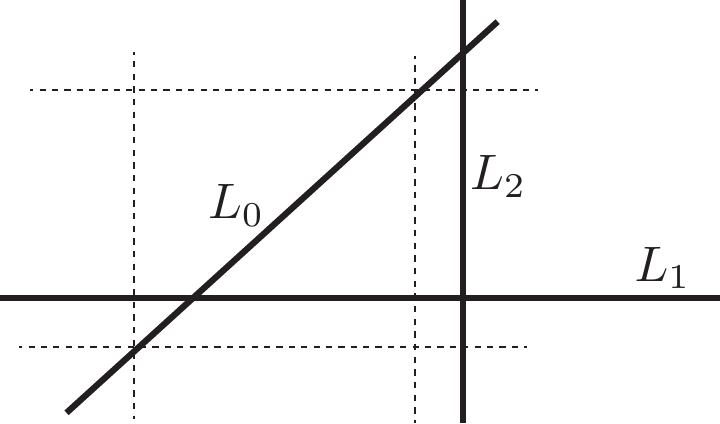}
\caption{Universal cover of an elliptic curve and its three Lagrangians}
\label{Figure15}
\end{figure}

The first example of homological mirror symmetry is the case 
of  elliptic curves and is discovered by Kontsevich \cite{Ko}.
He considered three Lagrangian submanifolds 
$L_0$, $L_1$, $L_2$ in $T^2$ (Figure \ref{Figure15}).

In the mirror, $L_0$, $L_1$ and $L_2$  correspond
to the structure sheaf, the ample line bundle $\mathscr O(1)$ and 
a point (sky scraper sheaf) $p \in T^2$, respectively.
\par
Kontsevich calculated the operation
\begin{equation}\label{m2elli}
\frak m_2: HF(L_0;L_1) \otimes HF(L_1;L_2) \to HF(L_0;L_2).
\end{equation}
In our situation $HF(L_0;L_1)$, $HF(L_1;L_2)$, $HF(L_0;L_2)$
are all rank one and so (\ref{m2elli}) is a number.
We move $L_2$ without changing its direction.
We also put a flat $U(1)$ connection.
These two operations give a family parametrized by one complex number,
which becomes a coordinate of the (universal cover of the) mirror elliptic curve.
\par
By definition, the map (\ref{m2elli}) is obtained by counting 
holomorphic triangles bounding $L_0$, $L_1$, $L_2$ together with an appropriate weight.
In our 1-dimensional case, we can calculate it explicitely 
and there is exactly one such a triangle in each homotopy 
class, which corresponds one to one to the natural numbers $k$.
The symplectic area of the triangle is 
$$\frac{(k+c)^2}{2}
$$
where $c$ is a parameter of $L_2$. Thus Kontsevich 
obtained the theta series:
$$
\sum_k\exp \left( -\frac{(k+z)^2}{2} \right)
$$
where $z$ is the coordinate which parametrizes a pair of $L_2$ and a 
flat $U(1)$ connection on it.
\par
In the mirror (\ref{m2elli}) should be
$$
H_{\overline\partial}(T^2;\mathscr O(1)) \otimes_{\C} \C \to \mathscr O(1)_p 
$$
where $p$ is the point corresponding to $z$.
Moving $p$, the family of the above operations should be the value of 
a global section of $\mathscr O(1)$ at $p$,
which is exactly the theta function.
This interesting calculation provides a nice evidence 
of homological mirror symmetry.
\par
Later Polishchuk and Zaslow \cite{PZ} studied the case of an elliptic curve 
in more detail and proved the homological mirror symmetry 
(in the cohomology level) in that case.
Their results are partially generalized in \cite{fu45}
to higher dimension, using the idea of family Floer homology
(See Section \ref{subsecfamily}).
\par
Abouzaid-Smith \cite{AS}  uses the fact that $\frak{Fukst}(T^2 \times T^2)$
is related to the category of filtered $A_{\infty}$
functors $\frak{Fukst}(T^2) \to  \frak{Fukst}(T^2)$
(See Section \ref{sec:COR+GAUGE}) to prove the homological mirror symmetry 
for a certain $T^4$.

\subsection{Generator of the Category 1: 
Decomposition of the monodromy.}

One difficulty in studying homological mirror symmetry is 
that categories are not easy to study or describe.
In fact, usually category is too huge to `calculate'.
A way to `calculate' an $A_{\infty}$ category 
is to find its generator and to 
compute the Floer homology groups between generators 
together with the structure operations.
This is the way how homological mirror symmetry 
have been proved in many cases.
We will discuss two of the methods to 
find a generator of an $A_{\infty}$ category
$\frak{Fukst}(X)$ of symplectic side.
\par
The first one is used in \cite{Sei3}.\footnote{
This paper is published rather recently.
However it appeared in the arXive in the year 2003, which is much 
earlier than its publication.}
Seidel's paper \cite{Sei3} studies the case of the quartic surface.
However his method is applied in greater generality.
Let us consider a family of Calabi-Yau manifolds:
\begin{equation}\label{fibration}
\pi: \frak X \to S^2.
\end{equation}
We consider its restriction 
$$
\pi: \pi^{-1}(D^2(\epsilon)) \to D^2(\epsilon)
$$
to a neighborhood of $0$ and assume 
that it is a maximal degeneration family, that is,
it satisfies (md1), (md2), (md3) above.
Take $q \in D^2(\epsilon) \setminus \{0\}$.
The monodromy around the singular fiber $M_0 = 
\pi^{-1}(0)$ defines a symplectic 
diffeomorphism $\Phi : M_q \to M_q$, where 
$M_q = \pi^{-1}(q)$.
Seidel's work is based on the following three points.
\begin{enumerate}
\item[(Sei1)]
One can define a $\Z$-grading 
on the Floer homology for a certain class of Lagrangian 
submanifolds in $M_q$.
\item[(Sei2)]
For $L$ in such a class and $d \in \Z$ the elements of
Floer homology 
$HF(L;\Phi^k(L))$ have degree $<d$ for sufficiently large $k$. 
\item[(Sei3)]
$\Phi$ is Hamiltonian isotopic to a 
composition of finitely many Dehn twists around Lagrangian 
spheres in $M_q$.
\end{enumerate}
\par
\begin{figure}[h]
\centering
\includegraphics[scale=0.3]{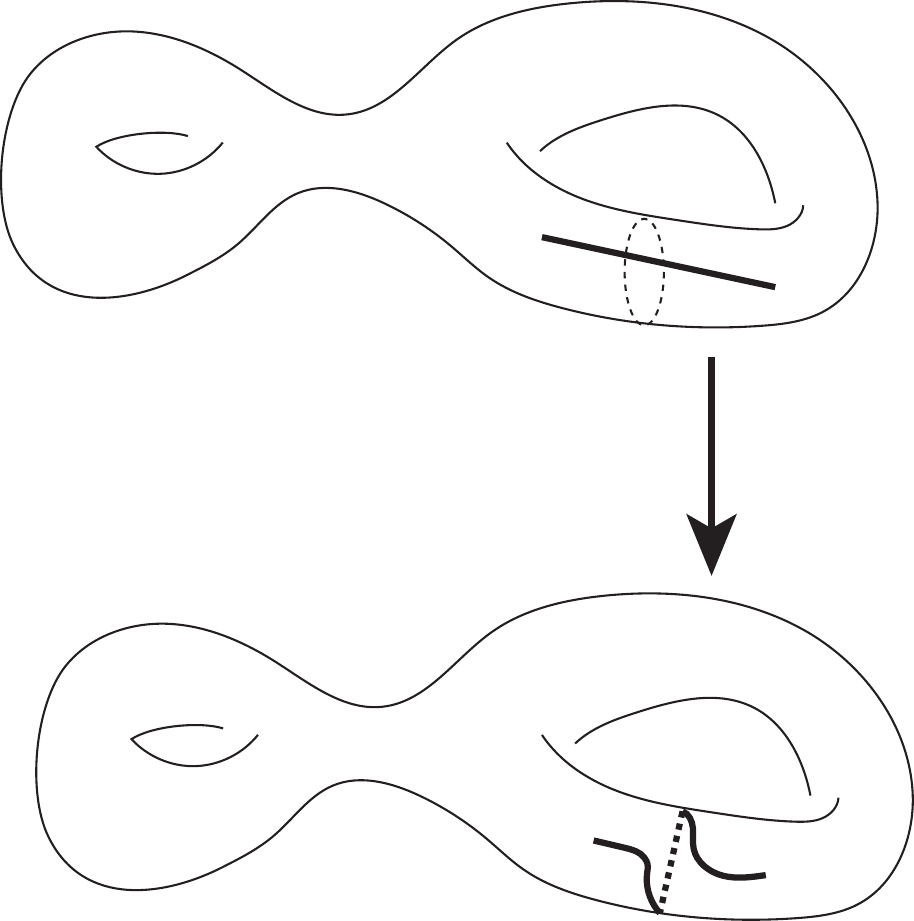}
\caption{Dehn twist}
\label{Figure2}
\end{figure}

A brief explanation of them are in order. We first 
explain (Sei3).
For a 2-dimensional manifold $\Sigma$ and 
a closed loop $\gamma \subset \Sigma$ we can 
associate a diffeomorphism, the Dehn twist,
$D_{\gamma} : \Sigma \to \Sigma$, which is
$[C] \mapsto [C] - ([\gamma]\cap[C])[\gamma]$
in homology.
(Figure \ref{Figure2}.)
Its higher dimensional analogue 
is a symplectic diffeomorphism, the Dehn twist, 
$\varphi_S: X \to X$
associated to a Lagrangian sphere $S^n$ in $2n$-dimensional 
symplectic manifold $X$. 
In ($n$-dimensional) homology, it is again
$[C] \mapsto [C] - ([S^n]\cap[C])[S^n]$.
Seidel used the fact that 
the monodromy of $M_q$ is decomposed 
into a composition of Dehn twists 
obtained by various Lagrangian spheres in $M_q$.
The reason of this decomposition is as follows.
\par
We consider a fibration (\ref{fibration}),
which we call a Lefschetz fibration.
We study the critical value of $\pi$.
One critical value is $0 \in S^2$.
By Assumption (md3), 
the fiber $M_0$ is much degenerate.
On the contrary, we can require that the other 
singular fibers $M_{q_i}$ are mildly degenerate.
Namely we assume $M_{q_i}$ is an immersed 
submanifold with one transversal self-intersection point.
It is a consequence of classical Picard-Lefschetz
theory that in such a case there exists 
a Lagrangian sphere $S_i^n$ and the monodromy 
around $q_i$ (that is $M_q \to M_q$ where $q$ is close to 
$q_i$) is a Dehn twist $\varphi_{S^n_i}$.
Thus the monodromy $\varphi$ around $0$ is 
decomposed into the composition of Dehn twists 
$\varphi_{S^n_i}$.
\par
In \cite{Sei15} Seidel calculated how 
the Dehn twist $\varphi_S$ acts on 
the $A_{\infty}$ category $\frak{Fukst}(X)$.
It is described by the following 
long exact sequence:
\begin{equation}\label{SLES}
\to HF(\varphi_S(L),L') 
\to HF(L,L') \to HF(L,S) \otimes HF(S,L') \to
\end{equation}
We consider the `difference' between two 
objects $L$ and $\varphi_S(L)$.
Via Yoneda embedding we regard them as right modules 
over $\frak{Fukst}(X)$.
To a Lagrangian submanifold  $L'$ (which plays the role of 
a `test function') the difference becomes 
$$
L' \mapsto HF(L,S) \otimes HF(S,L').
$$
So it is contained in the subcategory generated by 
$L' \to HF(S,L')$, that is, $S$ itself.
In this way, together with (Sei2)
(Sei3), the exact sequence (\ref{SLES}) shows that 
the set of  Lagrangian spheres $S^n_{q_i}$  
generates  $\frak{Fukst}(X)$.
\begin{rem}
We remark that (\ref{SLES}) is related to 
Floer's Dehn surgery triangle \cite{fl5} 
via Atiyah-Floer conjecture.
\end{rem}
\begin{rem}
We also remark that the above mentioned method is 
a part of the important project\footnote{which is on going 
and developing for a long time.} to calculate Floer homology 
of Lagrangian submanifolds based on 
symplectic Lefschetz fibration \cite{Don6}.
\end{rem}
\par
Note that in (Sei2) $\Z$-grading of Lagrangian 
Floer homology plays an essential role. 
We can define such $\Z$-grading 
in the case of a Lagrangian submanifold with vanishing Maslov index as follows
(\cite{Sei1}).
For a positive integer $n$ we consider $\C^n$ with standard symplectic form 
$\sum dx_i \wedge dy_i$ (where $z_i = x_i + \sqrt{-1}y_i$ is the standard 
coordinate.)
We denote by $\mathcal{LAG}(\C^n)$ the set of all oriented $n$-dimensional 
real linear subspaces $V$ of $\C^n$ such that $\sum dx_i \wedge dy_i$ is $0$ on $V$.
The manifold $\mathcal{LAG}(\C^n)$ is called the oriented Lagrangian Grassmannian.
Let $X$ be a symplectic manifold. For each $p \in X$ 
the tangent space $T_pX$ is identified with $\C^n$. (The identification 
respects the symplectic 
forms.) We collect $\mathcal{LAG}(T_pX)$ for all $p \in X$ and denote it by 
$\mathcal{LAG}(TX)$. The space $\mathcal{LAG}(TX)$ has a structure 
of a smooth manifold. There is a smooth map $\mathcal{LAG}(TX) \to X$ which sends 
$\mathcal{LAG}(T_pX)$ to $p$. This map gives a structure of a smooth fiber bundle
on $\mathcal{LAG}(TX)$ over $X$.
It is well known that $\pi_1(\mathcal{LAG}(T_pX)) = \Z$.
Let $K_X$ be the complex line bundle,
the $n$-th exterior power of the tangent bundle
and $SK_X$  its unit circle bundle. 
There is a bundle map $\mathcal{LAG}(TX) \to SK_X$
which associates  $[v^*_1 \wedge \dots \wedge v^*_n]$ to $V \in \mathcal{LAG}(TX)$, 
where $(v_1,\dots,v_n)$ is an oriented basis of $V$ and $v^*_i$ is its dual basis.
The next diagram commutes.
\begin{equation}
\begin{CD}
\pi_1(\mathcal{LAG}(T_pX)) 
@ >>>
\pi_1(\mathcal{LAG}(TX))
@ >>>
\pi_1(X)
\\
@ VVV @ VVV @ VVV\\
\pi_1(S^1)   
@ >>>
\pi_1(SK_X)
@ >>>
\pi_1(X)
\end{CD}
\end{equation}
Moreover the left vertical arrow is an isomorphism.
\par
We assume that $X$ is Calabi-Yau.
Therefore $K_X$ is the trivial bundle.
It follows from the diagram that there exists 
$\pi_1(\mathcal{LAG}(TX)) \to \Z$ which is an isomorphism 
on $\pi_1(\mathcal{LAG}(T_pX))$.
Therefore there exists a $\Z$ fold cover 
$\tilde{\mathcal{LAG}}(TX)$ of $\mathcal{LAG}(TX)$ 
which restricts to the universal cover on each $\mathcal{LAG}(T_pX)$.
\par
Let $L$ be an oriented Lagrangian submanifold of $X$.
The association $p \mapsto T_pL$ defines a section of 
$\mathcal{LAG}(TX)\vert_L$, which we denote by $s_L$.
We can show that there exists a section 
$\tilde s_L$ of  $\tilde{\mathcal{LAG}}(TX)\vert_L$
which lifts $s_L$ if and only if the Maslov index of $L$ is zero.
We call the lift $\tilde s_L$ the grading of $L$.
The graded Lagrangian submanifold is a pair $(L,\tilde s_L)$.
See \cite{Sei1}. It is proved there that 
for a pair of graded Lagrangian submanifolds 
$(L_i,\tilde s_{L_i})$, $i=1,2$,  we can define a 
$\Z$ grading of elements $p \in L_1\cap L_2$ which 
induces a $\Z$ grading of Floer homology.
\par
We do not discuss the proof of (Se2) here.
\par
The method of \cite{Sei3} is expanded by various authors 
and produced various examples where homological mirror symmetry 
is proved. One important recent development is \cite{sheri1}.\footnote{
In fact \cite{sheri1} uses the argument we will explain in the next subsection 
rather than the one in this subsection to find a generator. However 
\cite{sheri1}'s argument can be regarded as a generalization of \cite{Sei3}.}

\subsection{Generator of the Category 2: 
Open-closed map.}
\label{subsec:gene2}

Before explaining the second method, we recall the notion 
of a split generation of an $A_{\infty}$ category.
Let $\mathscr C$ be an $A_{\infty}$ category.
We enhance the set of objects of $\mathscr C$ in the 
following three ways.
\begin{enumerate}
\item[(spg1)]
If $c$ is an object of $\mathscr C$ and $k \in \Z$ 
we include an object $c[k]$, its degree shift, such that
$$
\mathscr C^d(c[k],c') = \mathscr C^{d-k}(c,c')
\qquad
\mathscr C^d(c',c[k]) = \mathscr C^{d+k}(c',c).
$$
The (higher) composition of morphisms between 
shifted objects is the same as those before shifted except 
the sign. The sign is determined for example 
if $xs \in  \mathscr C(c[1],c')$ which is identified 
with $x \in \mathscr C(c,c')$ then 
$\frak m_2(xs,y) = (-1)^{\deg y+1}\frak m_2(x,y)s$.
Here $s$ is an operator which shifts the degree by $1$.
\item[(spg2)]
If $c,c'$ are objects and $x \in \mathscr C(c,c')$ is a 
closed morphism, we include the mapping cone of $x: c \to c'$ 
as an object. (We can do so in the $A_{\infty}$ category of 
twisted complexes.)
\item[(spg3)]
We include `direct summand' of an object as a new object.
We can do so by using the notion of an idempotent. See 
\cite[Section 4]{Sei2}.
\end{enumerate}
\par
For any $A_{\infty}$ category $\mathscr C$ 
we can repeatedly add objects by one of the above 
three operations. We denote by $\mathscr C^+$
the $A_{\infty}$ category obtained in this way.

\begin{defn}
Let $X$ be a symplectic manifold and $\mathbb L$  a 
finite set of relatively spin Lagrangian submanifolds.
We say $\mathbb L$ is a split generator if the following 
holds. 
\par
Let $L'$ be an arbitrary relatively spin Lagrangian submanifold.
We put $\mathbb L' = \mathbb L \cup \{L'\}$. 
We obtain a filtered $A_{\infty}$ category
$\frak{Fukst}(\mathbb L)$ (resp. $\frak{Fukst}(\mathbb L')$)
whose objects are pairs $(L,b)$ such that $L 
\in \mathbb L$ (resp. $L 
\in \mathbb L'$) and $b$ is a bounding cochain of $L$ (resp. $L'$).
\par
We change the coefficient ring from $\Lambda_0$ to $\Lambda$ 
to obtain filtered $A_{\infty}$ 
categories $\frak{Fukst}(\mathbb L)_{\Lambda}$ and $\frak{Fukst}(\mathbb L')_{\Lambda}$.
\par
Then by adding objects as in (spg1), (spg2), (spg3)
we obtain fitered $A_{\infty}$ 
categories  $\frak{Fukst}(\mathbb L)_{\Lambda}^+$ and 
$\frak{Fukst}(\mathbb L')_{\Lambda}^+$.
\par
Now we require that the canonical embedding
\begin{equation}\label{geneiso}
\frak{Fukst}(\mathbb L)_{\Lambda}^+ \to \frak{Fukst}(\mathbb L')_{\Lambda}^+
\end{equation}
is a homotopy equivalece of categories.
\end{defn}

The second method to find a generator has an origin 
in the following proposal (due to Kontsevich\footnote{The author heard this idea in early 2000's from 
Kontsevich in his E-mail, as an idea to prove the 
homological mirror symmetry for torus.}).
Let $(X,\omega)$ be a symplectic manifold.
We consider the product $X\times X$ together with 
its symplectic form $-\pi_1^*\omega + \pi_2^*\omega$,
which we denote by $-X \times X$.
The diagonal $\Delta = \{(x,x) 
\in -X \times X \mid x \in X\}$ is a Lagrangian submanifold of $-X \times X$
and $(\Delta,0)$ is an object of $\frak{Fukst}(-X\times X)$.
\par
Suppose that $\mathbb L$ is a finite set of relatively 
spin Lagrangian submanifolds of $X$.
We put
$$
\mathbb L \times \mathbb L
= \{L\times L' \subset -X \times X \mid L,L' \in \mathbb L\}.
$$
We also put $(\mathbb L \times \mathbb L)' 
= (\mathbb L \times \mathbb L) \cup \{\Delta\}.$
We define:
\begin{equation}\label{geneiso2}
\frak{Fukst}(\mathbb L\times \mathbb L)_{\Lambda}^+ \to \frak{Fukst}((\mathbb L\times \mathbb L)')_{\Lambda}^+
\end{equation}
in the same way as (\ref{geneiso}).
\par
The proposal claims that if (\ref{geneiso2}) is a homotopy  
equivalence then $\mathbb L$ is a split generator.
\par
One can justify this proposal by using Lagrangian correspondence 
as follows.
As we will explain in Section \ref{sec:COR+GAUGE}, 
a Lagrangian submanifold $\frak L$ of $-X \times X$ together 
with its bounding cochain $\frak b$ defines a filtered $A_{\infty}$
functor 
$$
\mathcal{W}_{(\frak L,\frak b)} : \frak{Fukst}(X) \to \frak{Fukst}(X).
$$
In the case $(\frak L,\frak b) =  (\Delta,0)$ the 
functor $\mathcal{W}_{(\Delta,0)}$ is the identity functor.
On the other hand, it is easy to see that, if $\frak L \in \mathbb L \times \mathbb L$
and $\frak b$ is obtained 
from a pair of bounding cochains of elements of $\mathbb L$
(see \cite[Section 16]{cor}), 
then the image of the functor  $\mathcal{W}_{(\frak L,\frak b)}$
is contained in $\frak{Fukst}(\mathbb L)$.
Therefore if (\ref{geneiso2}) is a homotopy equivalence then (\ref{geneiso})
is a homotopy equivalence.
\par
The condition that (\ref{geneiso2}) is a homotopy equivalence 
could be rewritten by using the open-closed map $\frak p$ as follows.
For an $A_{\infty}$ category $\mathscr C$ we can define its Hochshild 
(co)homology $HH(\mathscr C)$ and cyclic homology $HC(\mathscr C)$.
In the case when $X$ is compact\footnote{When 
$X$ is non-compact the quantum cohomology $HQ^*(X;\Lambda_0)$ 
should be replaced by the symplectic homology.
See \cite{Sei12,Ga}.} there are maps
$$
\frak q_{\mathbb L}: HQ^*(X;\Lambda_0) \to HH^*(\frak{Fukst}(\mathbb L)),
\quad
\frak p_{\mathbb L}: HC_*(\frak{Fukst}(\mathbb L)) \to H_*(X;\Lambda_0).
$$
(See \cite{fl255,Ko,Sei12,Al,fooobook,BC}.) The map
$\frak q_{\mathbb L}$ is a ring homomorphism
from the quantum cohomology $HQ^*(X;\Lambda_0)$ to 
the Hochshild cohomology.
The domain of $\frak p_{\mathbb L}$ can be taken also as 
the Hochshild homology $HH_*(\frak{Fukst}(\mathbb L))$. Then 
it is dual to $\frak q_{\mathbb L}$. (See \cite{AFOOO,Ga}.)
There is a heuristic argument that 
the condition (\ref{geneiso2}) being a homotopy equivalence 
becomes the following:
\begin{enumerate}
\item[($\star$)]
The unit $1_X \in H_*(X;\Lambda)$ 
is contained in the image of $\frak p_{\mathbb L}\otimes \Lambda$.
\end{enumerate}
Let us explain the relation between the condition ($\star$) on the open-closed map $\frak p_{\mathbb L}$ 
and the condition that (\ref{geneiso2}) is a homotopy equivalence.
Let $(L \times L,b\times b)$ be an object of  
$\frak{Fukst}(\mathbb L\times \mathbb L)$.
Here $L \in \mathbb L$ and $b$ is a bounding cochain of 
$L$. It induces a bounding cochain $b \times b$ 
of $L \times L$.\footnote{See \cite[Section 16]{cor}.}
We simplify the situation and assume that
there exist closed morphisms 
$$
\frak x : (L \times L,b\times b) \to (\Delta,0),
\quad 
\frak y : (\Delta,0)  \to (L \times L,b\times b)
$$
such that
\begin{equation}\label{idintheimage}
\frak m_2^{-X \times X}(\frak x,\frak y) = 1_{\Delta}.
\end{equation}
Here $\frak m_2^{-X \times X}$ is the structure operations of 
$\frak{Fukst}(-X \times X;(\mathbb L \times \mathbb L) \cup \{\Delta\})$
and $1_{\Delta} = HF(\Delta) \cong H(X)$ is the unit.
\par
The equality (\ref{idintheimage}) is a simplified version of the condition that (\ref{geneiso2}) 
becomes a homotopy equivalence of categories.
\par
We assume furthermore that 
$\frak x$ (resp. $\frak y$) is represented as 
$\frak x = {\rm PD}[p] \in H^n(L) \cong H_0(L)$, (resp. $\frak y = {\rm PD}[L]
\in H^0(L) \cong H_n(L)$).
Here $p \in L = (L \times L) \cap \Delta$.
In this (over)simplified case, 
the product in the left hand side is defined by the moduli space 
as in Figure \ref{Figure4} below.
\begin{figure}[h]
\centering
\includegraphics[scale=0.3]{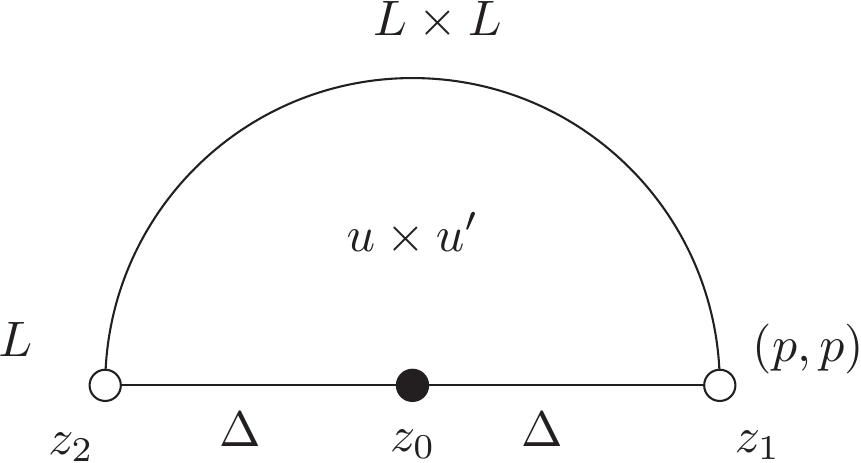}
\caption{A pseudo-holomorphic map to $-X \times X$.}
\label{Figure4}
\end{figure}
We (conformally) identify the semi-circle in the figure with a triangle.
Three vertices of the triangle are identified with the points $z_0,z_1,z_2$ 
in the figure, respectively, and $z_0$ is the vertex corresponding to the output.
The direct product $u \times u'$ is a map from the semi-circle to $-X \times X$ which is 
required to be pseudo-holomorphic.
The part of the boundary of the semi-circle which is the intersection 
of the semi-circle with $S^1$, is required to be mapped to $L \times L$.
The other part of the boundary is required to be mapped to the diagonal.
The points $z_1,z_2$ are required to be mapped to 
$(p,p)$ and $L = (L\times L) \cap \Delta$,
respectively.
(The condition for $z_2$ is actually automatic in this case.)
We consider the moduli spaces of such maps and use the evaluation map
at $z_0$. The sum of such output (with appropriate weight) is 
$\frak m_2^{-X \times X}(\frak x,\frak y)$ by definition.\footnote{
This is the case when the bounding cochain $b$ is zero. If $b$ is 
nonzero we need to consider more marked points on the circle part of 
the boundary.}
It is a chain of the diagonal $\cong X$.
\par
\begin{figure}[h]
\centering
\includegraphics[scale=0.3]{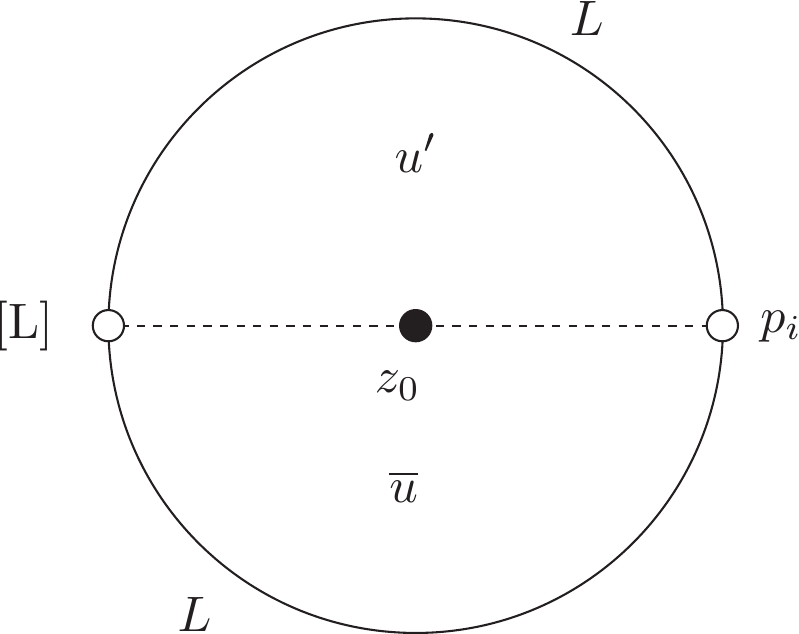}
\caption{A pseudo-holomorphic map to $X$.}
\label{Figure5}
\end{figure}

Now using the fact that the almlost complex structure we use 
on $-X \times X$ is $-J_X \oplus J_X$ and $u\times u'$ is 
pseudo-holomorphic, we can use reflexion principle 
to obtain the map in the Figure \ref{Figure5}.
Here $\overline u$ is defined by $\overline u(z) = u(
\overline z)$.
The moduli space of pseudo-holomorphic curves 
as in Figure \ref{Figure5} is the one we use to define 
the open-closed map
$
\frak p_{\mathbb L}([p]).
$
Thus (\ref{idintheimage}) is equivalent to 
$
\frak p_{\mathbb L}([p]) = 1_{X}.
$
It implies ($\star$).
\par
\begin{figure}[h]
\centering
\includegraphics[scale=0.3]{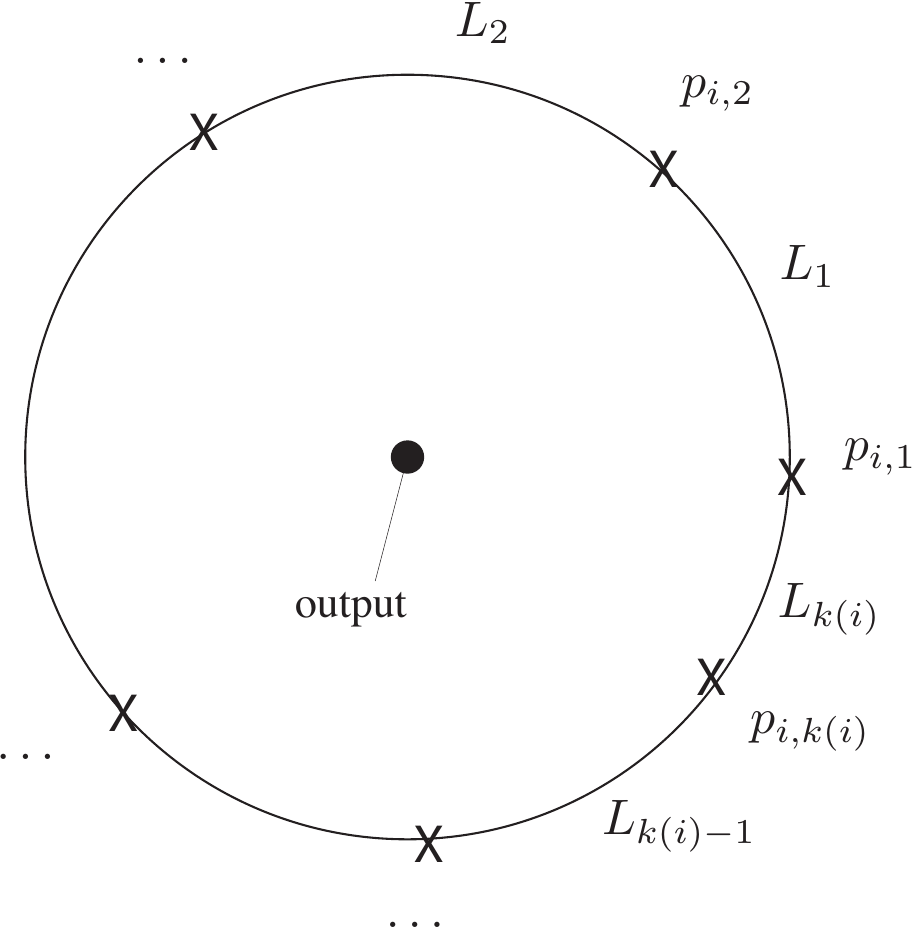}
\caption{A pseudo-holomorphic map to define the open-closed map.}
\label{Figure52}
\end{figure}
\begin{figure}[h]
\centering
\includegraphics[scale=0.5]{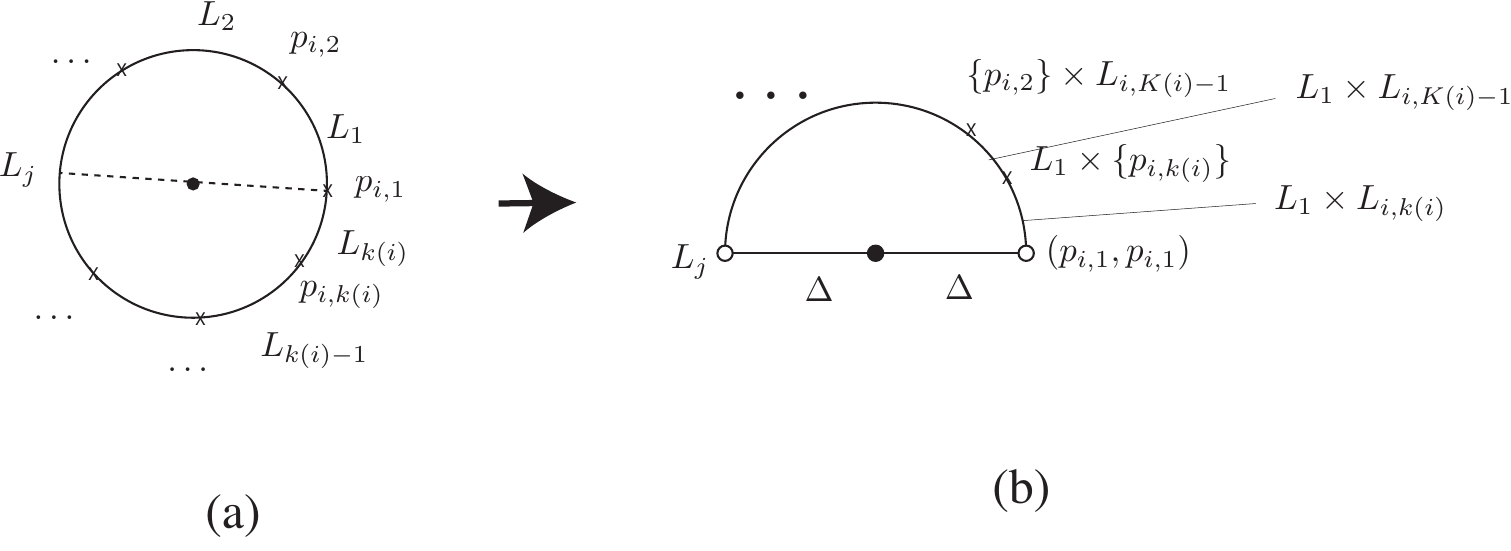}
\caption{A pseudo-holomorphic polygon obtained by a reflexion.}
\label{Figure53}
\end{figure}

In the more general case where
$$
\sum_i \frak p_{\mathbb L}([p_{i,1}]\otimes \dots \otimes [p_{i,k(i)}]) = 1_{X}
$$
the left hand side is obtained from the moduli space depicted by Figure \ref{Figure52}.
It corresponds to the polygon in Figure \ref{Figure53} (b) on $-X\times X$.
This polygon is obtained from the map (\ref{Figure52}) 
by reflexion at the dotted line of Figure \ref{Figure53} (a).

\par
Abouzaid \cite{Ab2} proved 
that under the condition ($\star$) (where $H_*(X;\Lambda)$ 
is replaced by the symplectic homology \cite{fh,BO}) 
$\mathbb L$ is a split generator, 
in a certain non-compact situation.
The compact case is on the way being written (\cite{AFOOO}).
The proof is based on the Cardy relation \cite{CL}
which identifies inner products of 
cyclic (or Hochshild) homology (See \cite{fooo:Ast,Sh,AFOOO}.) and the Poincar\'e duality via 
the open-closed map $\frak p$.\footnote{In the non-compact case, 
which Abouzaid established in  \cite{Ab2}, there is no Poincar\'e duality
on the ambient space. Abouzaid avoid using Poincar\'e duality by a skillful argument which 
we do not discuss here.}
It does not directly use the above explained idea to relate 
the condition ($\star$) to (\ref{geneiso2}).\footnote
{The author believes that there is an alternative proof using 
the idea to relate open-closed map to (\ref{geneiso2}), directly.
Such a proof is not worked out yet.}
\par
As we mentioned already one of the origins of this generating 
criteria is Kontsevich's proposal which seems to be 
related to a similar idea in the study of 
algebraic cycles in algebraic geometry 
(via mirror symmetry).
\par
Another origin is in symplectic geometry, 
such as those in \cite{Al,BC}.
It can be stated as follows.
Let $L,L'$ be two Lagrangian submanifolds of $X$.
Suppose
\begin{equation}\label{PDnonzero}
\langle \frak p(1_L),\frak p(1_{L'}) \rangle \ne 0.
\end{equation}
where $\langle \cdot \rangle : H_*(X) \otimes H_*(X) \to 
\Lambda$ is the Poincar\'e pairing, $\frak p$ is the 
open-closed map and $1_L$, $1_{L'}$ are the fundamental 
classes of $L$ and $L'$. 
Then, for any Hamiltonian diffeomorphism $\varphi : X \to X$, 
we have
\begin{equation}\label{ineq621}
\varphi(L) \cap L' \ne \emptyset.
\end{equation}
Note that the leading order term of the left hand side of (\ref{PDnonzero}) (that is, 
the term which does not contain $T$) is the intersection 
number $L\cdot L'$. So (\ref{ineq621}) 
is an enhancement of the obvious fact that $L\cdot L' \ne 0$
implies $
\varphi(L) \cap L' \ne \emptyset
$.

\subsection{Family Floer homology.}
\label{subsecfamily}

Another approach to homological mirror symmetry is by 
using family Floer homolology.
This approach is related to Strominger-Yau-Zaslow's proposal 
\cite{SYZ} to construct the mirror manifold via a dual torus fibration.
We first briefly review their proposal.
\par
We consider a symplectic manifold $X$ together with 
a map $\pi : X \to B$, such that for a `generic' point $\frak b$ of $B$ 
the fiber $L_{\frak b}: = \pi^{-1}({\frak b})$ is a Lagrangian torus.
Let $B_0 \subset B$ be the set of such `generic' points
and put $X_0 = \pi^{-1}(B_0)$.
The tangent space $T_{\frak b}B_0$ is identified with the first cohomology group 
$H^1(L_{\frak b};\R)$ of the fiber. The group $H^1(L_b;\R)$ 
contains a lattice $H^1(L_{\frak b};\Z)$.
The cotangent fiber $T^*_{\frak b}B_0$ is identified with 
$H^1(L_{\frak b};\R)$ which contains $H^1(L_{\frak b};\Z)$.
Frequently it is also assumed that $\pi : X \to B$ 
has a section $s: B \to X$ such that its image is a Lagrangian submanifold.
In that case $X_0$ is identified with 
the quotient of $T^*_{\frak b}B_0$ by the lattice $\Gamma: = \cup_{{\frak b}\in B_0} H^1(L_{\frak b};\Z)$.
(Here $s(\frak b)$ becomes $0 \in H^1(L_{\frak b};\R)/H^1(L_{\frak b};\Z)$
by this identification.)
Note that the total space of the 
cotangent bundle has a symplectic structure and the above 
identification respects the symplectic structures.
We consider the fiber-wise dual $TB_0$ and its dual lattice $\Gamma^{\vee}$.
One can show that the symplectic structure on $T^*B_0/\Gamma$ 
induces a complex structure on $TB_0/\Gamma^{\vee}$.
We put $X_0^{\vee}: = TB_0/\Gamma^{\vee}$.
The complex manifold $X_0^{\vee}$ is called a semi-flat mirror to $X_0$.
\par
The symplectic structure on $X_0$ extends to $X$ (by definition).
However the complex structure on $X_0^{\vee}$ in general does 
not extend to its compactification.
It is conjectured that there is a `quantum correction' 
to the complex structure of $X_0^{\vee}$ 
(which is determined by a certain `count' of holomorphic disks 
bounding a Lagrangian fibers $L_{\frak b}$)
so that, after correction, the complex manifold 
$X_0^{\vee}$ is compactified to a compact Calabi-Yau 
manifold $X^{\vee}$ and the map $X_0^{\vee}\to B_0$ extends to 
$X^{\vee}\to B$.
See \cite{fu52,KS2,GrS} for example.
When a mirror manifold is obtained in this way, 
$\pi : X \to B$ is called a SYZ-fibration.
\par
The family Floer homology program is a proposal to 
use SYZ picture of mirror symmetry to 
produce homological mirror functor 
$\frak{Fukst}(X) \to {\mathbb{DG}}(\mathscr{SH}(X^{\vee}))$.\footnote{
This idea was first communicated by M. Kontsevich to the author 
in 1997 during the author's stay in IHES.}
We consider a Calabi-Yau manifold $X$ together with an 
SYZ fibration $\pi: X \to B$. We consider a symplectic structure 
of $X$ only, forgetting its complex structure.
(We need an almost complex structure $J$ of $X$ to define the notion 
of a pseudo-holomorphic curve. However $J$ may not be integrable.) 
Suppose that we obtain its (SYZ) mirror  $X^{\vee}\to B$.
$X^{\vee}$ is a Calabi-Yau manifold. We regard it as a 
complex manifold, forgetting its symplectic structure 
(that is, the K\"ahler form).
\par
We consider a Lagrangian submanifold $L$ of $X$.
We assume, for simplicity, that $L$ is transversal to the 
fibers $L_{\frak b}:= \pi^{-1}({\frak b})$.
\par
A point of $X_0^{\vee}$ is a pair of an element $\frak b \in B_0$ 
and a cohomology class $a \in H^1(T_{\frak b};\R)/H^1(T_{\frak b};\Z)$.
The class $a$ determines uniquely a flat $U(1)$ bundle 
$\mathcal L_a$ on $L_{\frak b}$.
Let $\mathcal E_L$ be a vector bundle of $X^{\vee}$ 
which is the homological mirror to $L$.
In this situation, the homological mirror symmetry conjecture
`implies':
\begin{equation}\label{HMS1}
(\mathcal E_L)_{(\frak b,a)} = HF((L_{\frak b},a);L).
\end{equation}
Here the left hand side is the fiber of the vector 
bundle and is a $\C$-vector space.
The right hand side is the Floer homology 
of $L$ and $L_{\frak b}$. The role of $a$ can be explained 
in two different ways.
\begin{enumerate}
\item[(rega1)]
We regard the right hand side of (\ref{HMS1}) as the Floer homology with local coefficient.
Note that the boundary operator of  Lagrangian Floer homology 
is by definition a signed and weighted count of holomorphic 
strips $u$ with an appropriate boundary condition.
(See (\ref{Floersequ2}) and (path2)' (path3)' in Section \ref{sec:lagHF}.)
We put the weight $T^{\int u^*\omega}$ usually.
(The variable $T$ is a formal parameter, the Novikov parameter).
Here we take the weight $e^{-\int u^*\omega} {\rm Hol}_a(\partial u)$
instead, where $e$ is the Napier's number and 
${\rm Hol}_a(\partial u)$ is the holonomy of the flat line bundle 
$\mathcal L_a$ along the loop $u\vert_{\partial D^2}$.
(${\rm Hol}_a(\partial u) \in U(1) = \{z \in\C\mid\vert z\vert = 1\}$.)
\item[(rega2)] 
We regard $a \in H^1(L_{\frak b};\R) \subset H^1(L_{\frak b};\Lambda_0)$
and regard it as a bounding cochain as explained in Remark \ref{rem44}.
Note that $H^1(L_{\frak b};\R)/H^1(L_{\frak b};\Z)$
which appears in Remark \ref{rem44} is the moduli space of 
flat $U(1)$ bundle on $L$.
\end{enumerate}

We can then try to use (\ref{HMS1}) as the 
{\it definition} of the vector bundle $\mathcal E_L$.
Namely we conjecture that the $(\frak b,a) \in X_0^{\vee}$-parametrized family 
of vector spaces $HF((L_{\frak b},a);L)$ has a structure of 
a holomorphic vector bundle on $X_0^{\vee}$ that can be extended to $X^{\vee}$.

Note that the Lagrangian submanifold $S = s(B)$ (the image of the Lagrangian 
section $s$) has the property
$HF((L_{\frak b},a);S) = \C$. So we conjecture that it becomes the trivial 
line bundle (that is, the structure sheaf) of $X^{\vee}$.

The product structure $\frak m_2$ defines a map
$$
HF((L_{\frak b},a);S) \otimes HF(S;L)  \to HF((L_{\frak b},a);L).
$$
Since $HF((L_{\frak b},a);S) \cong \C$ we find that 
an element of  $HF(S;L)$ determines a section of 
$\mathcal E_L$. A way to define a holomorphic structure 
on the family $HF((L_{\frak b},a);L)$ is requiring this section 
to be holomorphic. 
\par
The isomorphism $HF^0(S;L) \cong H^0_{\overline{\partial}}(\mathcal E_L)$
is a part of the claim that $L \mapsto \mathcal E_L$ 
defines a fully faithful embedding of 
$A_{\infty}$ categories.
\par
In \cite{fu455}, the homological 
mirror symmetry of symplectic and complex tori 
are proved by this method, in the case when  Lagrangian 
submanifolds $L$ are flat and are transversal to the fibers.
A certain discussion in the case when there is a singular fiber 
and non-flat $L$ can be found in \cite{fu52}.
\par
A few years after \cite{fu455}, Kontsevich-Soibelman \cite{KS}
proposed to use rigid analytic geometry 
to study homological mirror symmetry and family 
Floer homology.
Note that when we use the weight
$e^{-\int u^*\omega} {\rm Hol}_a(\partial u)$
to define  Floer's boundary operator 
(and also a similar weight for the structure operations 
of an $A_{\infty}$ category), it is in general 
difficult to prove that the series defining the boundary operator
etc. converges. 
In the case of flat Lagrangian submanifolds in  
symplectic tori, the series defining structure operations
are certain variants of the theta series and actually converges 
very rapidly.\footnote{As mentioned in Subsection \ref{subsec:ellipticcurve}, 
the fact that theta series appears in the structure operations 
defining $\frak{Fukst}(X)$ was first observed by
Kontsevich in the case when $X$ is an elliptic curve.}   
\par
However beyond the case of tori and flat Lagrangian submanifolds, 
there is no method known to obtain 
an estimate of the number of holomorphic strips or 
disks and to show the convergence of the series with weight 
$e^{-\int u^*\omega} {\rm Hol}_a(\partial u)$.
\par
In symplectic Floer theory, the usual method 
is  introducing the universal Novikov ring 
and defining the boundary operator and the structure operations 
as a `formal power series'.
Kontsevich-Soibelman \cite{KS} observed that 
using the universal Novikov ring in the mirror (complex) side 
corresponds to studying the mirror manifold 
as a rigid analytic space.
In a series of papers \cite{toric1,toric2,fooo:Ast}
we worked out the idea using rigid analytic family 
of Floer homologies to study the case of toric manifolds.
Abouzaid \cite{Ab3,Ab4,Ab5} realized this program 
and proved a version of homological mirror symmetry 
in the case when the SYZ-fibration $X \to B$ has 
no singular fiber and the fibers never bound holomorphic 
disks.
J. Tu \cite{Tu} and H. Yuan \cite{Yu} partially relax these conditions 
and include the case where there exist singular fibers.
J. Solomon \cite{So} showed that for the fibers of SYZ-fibration, 
any $b \in H^{\rm odd}(T^{2n};\Lambda_0)$ satisfies 
the Maurer-Cartan equation (\ref{MCequation}). 

\subsection{Matrix factorizations and weak bounding cochains.}

Mirror symmetry is generalized to a manifold $M$ which is not 
necessary Calabi-Yau.
In such a case the mirror of $M$ is not a manifold 
but is a pair $(M^{\vee},W)$ of a manifold $M^{\vee}$ 
and a function $W$ on it.
Such a mirror symmetry is studied both 
in the case when $M$ is a symplectic manifold and the case when $M$ is a complex manifold.
\par\medskip
We discuss in this subsection  the case when $M$ is a symplectic manifold.
The other case is discussed in the next subsection.
In this case $(M^{\vee},W)$ is a pair of a  
manifold $M^{\vee}$, and a function $W$ on it.
There are cases both $W$ is a holomorphic function 
and $W$ is a rigid analytic function.
For simplicity we explain the case when $M^{\vee}$ is a 
complex manifold and $W$ is a holomorphic function on it.
The function $W$ is called a Landau-Ginzburg potential.
\par
In the Calabi-Yau case, the mirror to the quantum cohomology 
$HQ(M)$ is the $\overline\partial$-cohomology 
$H_{\overline\partial}(M^{\vee},\Lambda^*TM^{\vee})$. 
Here $\Lambda^*TM^{\vee}$ is an exterior power of the (complex) 
tangent bundle of $M^{\vee}$. 
This cohomology controls the `extended' deformation 
of the complex structure. In fact, 
Kodaira-Spencer theory says $H^1_{\overline\partial}(M^{\vee},TM^{\vee})$
controls the deformation of complex structures of $M^{\vee}$.
\par
In the case when $M$ is not Calabi-Yau, 
its mirror is a pair $(M^{\vee},W)$. 
In this case the mirror to the quantum cohomology ring 
$HQ(M)$ is a vector space which controls the deformation 
of the pair $(M^{\vee},W)$.
In many cases, $W$ has isolated critical points.
In such a case, only a neighborhood of 
the critical point set is used to study the mirror symmetry.
So $M^{\vee}$ is a local object (a neighborhood of finitely many 
points) and has no deformation.
The deformation of $W$ is controlled by the Jacobi ring.
Let $p \in M^{\vee}$ be a point such that $d_pW = 0$.
We take a complex coordinate $z_1,\dots,z_n$ of $M^{\vee}$
in a neighborhood of $p$. 
The Jacobi ring $Jac_p(W)$ is the quotient:
\begin{equation}
Jac_p(W)
: =
\frac{\mathscr O_p}{\left(\frac{\partial W}{\partial z_i}:i=1,\dots,n\right)}.
\end{equation}
Here $\mathscr O_p$ is the ring of germs of holomorphic functions 
of $M^{\vee}$ at $p$. 
The denominator is its ideal generated by the germs of 
partial derivatives $\frac{\partial W}{\partial z_i}$.
In fact the deformation of $W$ (in its first order approximation) can be written as 
the form
$$
t \mapsto W_t: = W + t F
$$ 
where $F \in \mathscr O_p$. Two such deformations $W_t$ and $W'_t$
are regarded as equivalent if there is a family of 
bi-holomorphic maps $\varphi_t$ with $\varphi_t(0)=0$
such that  $W'_t \sim W'_t\circ \varphi_t$.
It is easy to see that $W + t F$ is equivalent to 
$W + t F'$ if $F - F'$ is in the ideal 
$(\frac{\partial W}{\partial z_i}:i=1,\dots,n)$ in its first order approximation.
In the case when $p$ is an isolated critical point, 
the Jacobi ring $Jac_p(W)$ is finite dimensional.
\par
K. Saito \cite{Sai83} defined an inner-product (higher residue pairing) 
and  various other structures on the (family of)  Jacobi 
rings\footnote{Saito's work is much earlier than the discovery 
of quantum cohomology.} which is very close to the structure 
found (by Dubrovin \cite{dub}) on quantum cohomologies.
Such a structure is called a Frobenius manifold 
structure.  
These two structures are expected to be identified via 
mirror symmetry.
It seems that among the mathematicians, Givental \cite{Gi} first 
mentioned this conjecture.  The author does not know how 
this story was developed among physicists.
\par
The homological mirror symmetry between 
a symplectic manifold $M$ and a pair $(M^{\vee},W)$
is expected to be a relation between a filtered $A_{\infty}$
category $\frak{Fukst}(M)$ and 
a category of matrix factorizations of $(M^{\vee},W)$.
Let us briefly review the notion of a matrix factorization.
The notion of a matrix factorization is introduced by 
Eisenbud \cite{ei} a long time ago.
It is a $\Z_2$ graded module $C = C_0 \oplus C_1$ 
over $\mathscr O_p$
with degree one operator $d: C \to C$ such that
\begin{equation}\label{matfact}
d\circ d = W \cdot {\rm Id}.
\end{equation}
When $(C,d), (C',d')$ are  matrix factorizations 
we consider the set of $\mathscr O_p$
module homomorphisms $C \to C'$ and
denote it by $Hom(C,C')$.
It is a $\Z_2$ graded $\mathscr O_p$ module.
We define $\delta: Hom(C,C') \to Hom(C,C')$ by
\begin{equation}\label{matfactdiff}
\delta(\varphi) = d' \circ \varphi - (-1)^{\deg \varphi}\varphi \circ d.
\end{equation} 
It is easy to show $\delta \circ \delta = 0$.
Therefore there exists a DG-categry ${\rm Mat}(W;p)$
whose object is a matrix factorization and the
module of morphisms from $(C,d)$ to $(C',d')$ is $Hom(C,C')$,
with boundary operator $\delta$ and the obvious composition.
\par
The fact that the matrix factorizations define the D-brane category 
of the pair $(M^{\vee},W)$ is pointed out by physicists
(\cite{KL,HW}) and also in \cite{Or} in early 2000's.
\par
In Lagrangian Floer theory a similar 
structure was know by Floer and Oh \cite{Oh}
in 1990's.
Let us consider a pair $L_1$, $L_2$ of monotone Lagrangian submanifolds
intersecting transversally.\footnote{Hereafter we assume all the
Lagrangian submanifolds involved are oriented and relatively spin.}
Theorem \ref{LagFlOh} by Oh says if the minimal Maslov number of 
$L_1$ and $L_2$ are strictly larger than $2$ then 
Floer's boundary operator $d: CF(L_1,L_2) \to CF(L_1,L_2)$
satisfies 
$d\circ d = 0$. Here $CF(L_1,L_2)$ is the free abelian group 
whose basis is identified with the intersection points.
\par
In the case when the minimal maslov number is $2$, the equality 
$d\circ d = 0$ may not hold.
Suppose that the
minimal maslov number of a monotone Lagrangian submanifold 
$L$ is $2$.
We consider $\beta \in \pi_2(X;L)$ whose Maslov index $\mu(\beta)$
is $2$.  Let $\mathcal M_1(\beta)$ 
is the moduli space of pseudo-holomorphic maps
$
u: (D^2,\partial D^2) \to (X,L)
$
whose homotopy classes are $\beta$.
We identify $u$ and $u'$ if there exists 
a biholomorphic map $v: D^2 \to D^2$ such that
$u\circ v = u'$ and $v(1) = 1$.
Using the fact that there exists no holomorphic disk whose Maslov index is 
strictly smaller than $2$, we can show that $\mathcal M_1(\beta)$ 
is an $n$-dimensional smooth manifold without boundary 
for a generic compatible almost complex structure.
The map which sends $[u] \in \mathcal M_1(\beta)$ to 
$u(1) \in L$ is a smooth map $: \mathcal M_1(\beta) \to L$
between $n$-dimensional oriented manifolds 
without boundary and so its mapping degree $n_{\beta}$
is well-defined.  We define
\begin{equation}\label{cLLL}
c_L: = \sum_{\beta \in \pi_2(X;L),\,\,\, \mu(\beta)=2} n_{\beta}.
\end{equation}
In the case when $(L_1,L_2)$ is a pair of monotone Lagrangian submanifolds
intersecting transversally and with the minimal Maslov number $2$, 
Oh \cite{Oh} proved the following equality:
\begin{equation}\label{Ohmaslov2}
d\circ d = (c_{L_2} - c_{L_1})\cdot {\rm Id}.
\end{equation}
The equality (\ref{Ohmaslov2}) is the same as (\ref{matfact}) except 
$c_{L_2} - c_{L_1}$ is not a function but is 
a number (integer).
\par
We can generalize this story without assuming monotonicity as follows.
Let $L$ be a relatively spin Lagrangian submanifold of $X$.
We obtain a structure of a unital, curved filtered $A_{\infty}$ algebra 
on $H(L;\Lambda_0)$.\footnote{In case $L$ is immersed we can 
use $CF(L;\Lambda_0)$ as in (\ref{AJXCF}) and the story goes 
in the same way.}
We say $b \in H^1(L;\Lambda_0)$ is a weak bounding cochain 
if it satisfies:
\begin{equation}
\sum_{k=0}^{\infty} \frak m_k(b,\dots,b) = c \,\,{\bf e}_L
\end{equation}
for some constant $c \in \Lambda_0$.\footnote{In the case when $b \notin H(L;\Lambda_+)$
but $b \in H(L;\Lambda_0)$
there is an issue of convergence of the right hand side. However we can  still 
define the notion of weak bounding cochain in the same way as
Remark \ref{rem44}.}
Let $\widetilde{\mathcal M}_{\rm weak}(L)$ be the set of 
all weak bounding cochains of $L$.\footnote{Actually 
it is more natural to introduce an equivalence relation, gauge equivalence 
(See \cite[Section 4.3]{fooobook}.), between weak bounding cochains, and introduce 
a space ${\mathcal M}_{\rm weak}(L)$ consisting 
of gauge equivalence classes of weak bounding cochains.
The potential function then becomes a function on ${\mathcal M}_{\rm weak}(L)$.
It seems likely that ${\mathcal M}_{\rm weak}(L)$ becomes a 
certain version of an (Artin) stack in the rigid analytic category 
and $\frak{PO}$ is a rigid analytic function on it.
(This statement is not proved in the general case in the literature.)}
We define a potential function
$\frak{PO}: \widetilde{\mathcal M}_{\rm weak}(L) \to \Lambda_0$
by 
\begin{equation}\label{defpotential}
\sum_{k=0}^{\infty} \frak m_k(b,\dots,b) = \frak{PO}(b) \cdot {\bf e}_L.
\end{equation}
Let $(L_1,L_2)$ be a transversal pair of 
relatively spin Lagrangian submanifolds
and $b_i$ a weak bounding cochain of $L_i$ for $i=1,2$.
We define $d^{b_1,b_2} : CF(L_1,L_2;\Lambda_0) \to CF(L_1,L_2;\Lambda_0)$\footnote{$CF(L_1,L_2;\Lambda_0)$ is the free $\Lambda_0$ module whose basis 
is the set of intersection points $L_1 \cap L_2$.} by
\begin{equation}
d^{b_1,b_2}(x): = 
\sum_{k=0}^{\infty}\sum_{\ell=0}^{\infty}\frak m_{k+1+\ell}
(b_1^{k},x,b_2^{\ell}).
\end{equation}
Here $\frak m_*$ are the restrictions of the structure operations 
of the curved $A_{\infty}$ category whose objects are $L_1$ and $L_2$.
Using the $A_{\infty}$ relation and (\ref{defpotential})
we can show (\cite[Proposition 3.7.17]{fooobook}):
\begin{equation}\label{weakmatrix}
d^{b_1,b_2} \circ d^{b_1,b_2} = (\frak{PO}(b_2) - \frak{PO}(b_1))\cdot {\rm Id}.
\end{equation}
This equality shows that $\frak{PO}$ can be regarded as a Landau-Ginzburg 
potential and $(CF(L_1,L_2;\Lambda_0),d^{b_1,b_2})$ is identified with the morphism 
spaces between two matrix factorizations.
\begin{rem}
Let $(C,d)$, $(C',d')$ be matrix factorizations with respect to the 
Landau-Ginzburg potentials $W$, $W'$, respectively.
We consider $Hom(C,C')$ and define the boundary operator 
$\delta$ by the same formula as (\ref{matfactdiff}).
Then we have
$\delta\circ\delta = (W'-W)\cdot {\rm Id}$.
This formula coincides with (\ref{weakmatrix}).
\end{rem}
\begin{rem}
In 1990's Givental gave a comment to a talk by the author 
that the Landau-Ginzburg potential of the mirror may be obtained 
as a FOOO's obstruction class $\frak m_0$. 
The author does not know how 
this story was developed among physicists.
It seems that one of its origin is \cite{witten25}.
See also \cite{HIV}.
This fact is used in an important paper 
by Hori-Vafa \cite{HV}.
In \cite{HV} this fact is used in the case of certain toric 
manifolds. Cho-Oh \cite{CO} studied that case and calculated  (\ref{cLLL}) 
in several important toric manifolds. 
After an important work by Cho \cite{Ch}, FOOO \cite{toric1,toric2,fooo:Ast} 
further studied the case of toric manifolds.
See also \cite{Au1}.
\end{rem}
The formula (\ref{weakmatrix}) implies that if 
$\frak{PO}(b_1) = \frak{PO}(b_2)$
the Floer homology group $HF((L_1,b_1),(L_2,b_2);\Lambda_0)$
is defined as the homology group of the boundary operator $d^{b_1,b_2}$.
The Floer homology $HF((L_1,b_1),(L_2,b_2);\Lambda)$
is expected to be the mirror to the cohomology of the 
morphism complex of the category of
matrix factorizations of $W = \frak{PO}$.\footnote{
We remark that  $\frak{PO}$ is not a holomorphic function 
but is a rigid analytic function. So we need to 
use a rigid analytic analogue of the theory of 
a matrix factorization. Such a theory is not developed much yet.}
We remark that the matrix factorization category of $W$ at $p$
is trivial unless $p$ is a critical point of $W$.
\par
Thus $HF((L_1,b_1),(L_2,b_2);\Lambda)$ is expected to be zero 
unless $b_1$, $b_2$ are critical points of $\frak{PO}$.\footnote{
This fact is rigorously proved in many cases. 
However the author does not know the proof which works in complete 
generality. One issue is since $\widetilde{\mathcal M}_{\rm weak}(L)$ 
the domain of $\frak{PO}$ is not a usual smooth complex 
manifold, it is not so obvious what we mean 
by a critical point of $\frak{PO}$.}
We may also expect that the set of critical values of $\frak{PO}$ 
is a finite set.\footnote{This fact is  rigorously proved also in many cases
but not in complete generality.}
We define:
$$
\frak C: = \{ c \in \Lambda_0 \mid \exists 
(L,b),\,\, b \in \widetilde{\mathcal M}_{\rm weak}(L), 
\,\, \frak{PO}(b) = c,\,\,\, HF((L,b),(L,b);\Lambda) \ne 0\}.
$$
Then, for each $c \in \frak C$, we can define a 
strict and unital filtered $A_{\infty}$ 
category $\frak{Fukst}(X;c)$
so that its object is a pair $(L,b)$ of a Lagrangian submanifold $L$ 
and its weak bounding cochain $b$ with  $\frak{PO}(b) = c$.
The structure operations of $\frak{Fukst}(X;c)$ can 
be defined in the same way as the case of $\frak{Fukst}(X)$, 
except $\frak m_0$. Note that $\frak m_0^b(1) = \frak{PO}(b){\bf e}_{(L,b)} \ne 0$.
However we redefine $\frak m_0^b$ by putting  $\frak m_0^b(1) =0$.
In the same way as $\frak m_1^b\circ \frak m_1^b = 0$ (which 
follows from (\ref{weakmatrix})) we can check $A_{\infty}$ formula 
when we redefine $\frak m_0^b(1)$ to be $0$.
\par
Thus, in non-Calabi-Yau cases, we have several strict 
filtered $A_{\infty}$ 
categories associated to $X$.
In the case of  Landau-Ginzburg B-model $(M^{\vee},W)$, 
there are several matrix factorization categories.
Namely for each critical value of $W$ we can associate
a category of matrix factorizations.
The mirror symmetry is expected to identify them together with decompositions.
\par
Other than the toric case (\cite{fooo:Ast})\footnote{
To show that the result of \cite{fooo:Ast} 
implies this conjecture we need to show 
that the finitely many Lagrangian tori which are $T^n$ orbit 
and have nontrivial 
Floer homology, generates the Fukaya category 
of the toric manifold. This is 
a consequence of the result discussed in 
Subsection \ref{subsec:gene2}, but is still on the way being written.
In the case when the toric manifold this follows from \cite{el}.},
homological mirror symmetry of this type is proved for a certain toric degeneration \cite{NNU1,NNU2}
and monotone hypersurfaces of $\C P^n$ \cite{sheri2} and etc.
\par
The decompositions of the category are expected to correspond 
to the decompositions of the quantum cohomology as follows.

\begin{conj}
The quantum cohomology ring $HQ(M,\Lambda)$ is
decomposed into the direct product of the rings
\begin{equation}\label{form622}
HQ(M,\Lambda): = \prod_{c \in \frak C} HQ(M,\Lambda;c).
\end{equation} 
The closed open map 
$\frak q: HQ(M,\Lambda) \to HH(\frak{Fukst}(X))$
is decomposed into a direct product of the ring homomorphisms:
\begin{equation}
\frak q_c: HQ(M,\Lambda;c) \to HH(\frak{Fukst}(X);c).
\end{equation}
We have an equality 
\begin{equation}
c_1 \cup_q x = cx 
\end{equation}
for $x \in HQ(M,\Lambda;c)$. Here $c_1$ is the first Chern class 
and $\cup_q$ is the quantum cup product.
\end{conj}
This conjecture is known to be true in many cases including 
the case of toric manifolds. 
It seems unlikely that it is proved in complete generality at the stage 
when this article is written.
\par
The author heard this conjecture in a lecture by Kontsevich
in the year 2006.
The author does not know how the story developed among physicists. 
\par
Note that in  Landau-Ginzburg B-model $(M^{\vee},W)$
the Jacobi ring $Jac(M^{\vee},W)$ is decomposed as:
\begin{equation}\label{Jacdecomp}
Jac(W) = \prod_{c \in {\rm Crit}(W)} Jac_c(W).
\end{equation}
Here ${\rm Crit}(W)$ is the set of critical values of $W$.
If the set of critical points is isolated 
$J_c(W)$ is the direct product of $J_p(W)$ for critical points $p$ 
with $W(p) = c$. If the critical point set is not isolated 
the story seems to be more involved.
\par
It is expected that the two decompositions (\ref{form622}) and (\ref{Jacdecomp})
coincide via mirror symmetry.

\subsection{Seidel's directed $A_{\infty}$ category.}

In this subsection, we discuss the mirror symmetry between 
a complex manifold $M$ and a pair $(M^{\vee},W)$.
Here $M$ is not Calabi-Yau.
Since the present author is not an expert of this 
part of  homological mirror symmetry, this subsection is rather brief
compared to the importance of the topic.
When $\mathscr F_i$, $i=1,2$ are coherent sheaves of $M$,
the Serre duality says
\begin{equation}
H^k(Hom(\mathscr F_1,\mathscr F_2))
\cong 
H^{n-k}(Hom(\mathscr F_2,\mathscr F_1)\otimes K_M),
\end{equation}
where $n$ is the complex dimension of $M$ and $K_M$ is the 
canonical bundle, that is, the $n$-th exterior power of 
complex cotangent bundle of $M$.
In the Calabi-Yau case $K_M$ is trivial 
so 
$H^k(Hom(\mathscr F_1,\mathscr F_2))
\cong 
H^{n-k}(Hom(\mathscr F_2,\mathscr F_1))$.
Note that the Lagrangian Floer homology 
of compact Lagrangian submanifolds $L_i$ 
satisfies
$HF((L_1,b_1),(L_2,b_2))) \cong HF((L_2,b_2),(L_1,b_1))$,
which is a part of cyclic symmetry and is induced by the 
Poincar\'e duality.
In the Calabi-Yau case, this is consistent with  homological 
mirror symmetry. 
However in non-Calabi-Yau case this means that we 
need to study non-compact Lagrangian submanifolds. 
Seidel \cite{Sei2} defined a directed $A_{\infty}$ category 
from the pair $(M^{\vee},W)$ (in the exact situation) 
using a certain set of non-compact Lagrangian submanifolds.\footnote{
Seidel mentioned that such a theory is proposed by Kontsevich. 
See \cite{Ko2}. The work \cite{HIV} by physicists is also an 
origin of this construction.}

\begin{defn}\label{Directed}
A directed $A_{\infty}$ category $\mathscr C$ over $R$ 
is an $A_{\infty}$ category which is strict, unital and linear over $R$
and such that:
\begin{enumerate}
\item The set of its objects consist of finitely many elements $c_1,\dots,c_n$, 
indexed by a totally ordered set $i=1,\dots,n$.
\item 
If $i<j$ the set of morphisms $\mathscr C(c_j,c_i)$ is $\{0\}$.
\item 
The set of endmorphisms $\mathscr C(c_i,c_i)$ is $R$  whose basis is the unit.
\item
If $i<j$ the set of morphisms $\mathscr C(c_i,c_j)$ is a 
finite dimensional free $R$ module.
\end{enumerate}
\end{defn}
In complex geometry, an origin of such a directed category 
is Beilinson's paper \cite{Be}, where the category of coherent sheaves 
on the projective space $\mathbb P^n$ is studied. In that case $c_i = \mathscr O(i)$ for $i=0,\dots,n$ satisfies this condition for $\mathscr C(c_i,c_j): = H_{\overline\partial}
(Hom(c_i,c_j))$. (Here $\mathscr O(i)$ is a line bundle of degree $i$.)
This set is also a generator of the derived category of coherent sheaves.
Such a generator is called a full (strong) exceptional collection
(of the category of coherent sheaves).
This notion is defined in \cite{GR}.
\par
Seidel's construction is the case when the symplectic manifold $M^{\vee}$ is 
exact and $W: M^{\vee} \to \C$ has a mild singularity, that is, 
for each critical value $q$, the fiber $M^{\vee}_q = W^{-1}(q)$ is an 
immersed symplectic submanifold with one transversal self-intersection point.\footnote{In the case of quartic surface Seidel studied, this situation appears 
after removing a certain divisor.}
The non-compact Lagrangian submanifold Seidel studied is 
a kind of unstable manifold of $W$ which is called 
a Lefschetz thimble.
For $p \in M^{\vee}$ such that $D_pW \ne 0$, 
decompose the tangent space $T_pM^{\vee}$ as 
$
T_pM^{\vee} = T^v_pM^{\vee} \oplus T^h_pM^{\vee}
$
where $T^v_pM^{\vee} = T_p W^{-1}(W(p))$ and 
$T^h_pM^{\vee}$ is its orthogonal complement with respect to the 
symplectic form. 
Using this decomposition we can define a horizontal lift of 
a path in $\C$. 
For each critical value $q_i \in \C$ take a path 
$\gamma_i:[0,\infty) \to \C$ such that $\gamma_i(0) = q_i$ and that,  
if $t$ is sufficiently large, then $\gamma_i(t) = t + \sqrt{-1}C_i$ for a 
certain $C_i \in \R$,
such that it does not intersect with critical values for $t\ne 0$.
For $x \in W^{-1}(\gamma_i(t_0))$ we take a horizontal lift $\widetilde{\gamma}_x:[0,t_0]
\to M^{\vee}$
of $\gamma_i$ such that $\widetilde{\gamma}_x(t_0) = x$.
We consider the set $L_i$ of all $x \in W^{-1}(\gamma_i([0,\infty))$ such that 
$\widetilde{\gamma}_x(0)$ is a critical point of $M_{q_i} = W^{-1}(q_i)$.
It is known that $L_i$ is a Lagrangian submanifold of $M^{\vee}$ and 
is called a Lefschetz thimble.
For each $\gamma_i(t_0)$, $t_0\ne 0$ the fiber 
$W^{-1}(\gamma_i(t_0)) \cap L_i$ is known to be a Lagrangian sphere of 
the symplectic manifold $W^{-1}(\gamma_i(t_0))$ and is called a vanishing 
cycle.
\par
We assume that the path $\gamma_i$ and $\gamma_j$ do not intersect each other 
for $i\ne j$ and that $C_i < C_j$ for $i<j$.
The Lefschetz thimbles $L_i$ $i=1,\dots,n$ do not intersect each other.
We will perturb it slightly near the infinity so that they intersect as follows.
We take a Hamiltonian $H: M^{\vee} \to \R$ such that $H(x) = \epsilon\Vert W(x)\Vert^2$.
We consider the time one map $\varphi_H$ 
of the Hamiltonian vector field $\frak X_{H}$. The map $\varphi_H$ rotates the image of  $\gamma_i$
slightly to the counter clockwise direction so that $\varphi_H(L_i) \cap L_j \ne \emptyset$ if and only if $i<j$.
We put
$$
\mathscr{S}(L_i,L_j) = CF(\varphi_H(L_i),L_j),
$$
where the right hand side is the usual Floer complex.
This is zero if $j<i$ and so (\ref{Directed}) is 
satisfied. For $i=j$ we define $\mathscr{S}(L_i,L_i) = R$, the ground ring.
Under a certain exactness assumption on the symplectic form 
and a certain condition on the behavior of $W$ near infinity,
Seidel defined a directed $A_{\infty}$ category\footnote{Sometime called Fukaya-Seidel 
category.} 
$\mathscr{S}(M^{\vee},W)$ the set of whose objects are $\{L_i\mid i= 1,\dots,n\}$
and the morphism complex is $\mathscr{S}(L_i,L_j)$.
\par
It is possible to increase objects by requiring the Lagrangian 
submanifold to have an asymptotic behavior similar to $L_i$'s.
Namely we can consider $L$ such that outside a compact set 
$L$ is obtained by a parallel transport along an arc $t \mapsto t+C\sqrt{-1}$.
(The case of a compact Lagrangian submanifold $L$ is included.)
\begin{rem}
In place of rotating Lagrangian submanifolds a bit
one can take a Hamiltonian such as $\Vert W(x)\Vert^4$ so that 
the rotation angle goes $\infty$ as one goes to the infinity of $\C$.
Abouzaid-Seidel \cite{ASe} defined such an $A_{\infty}$ category 
which is called a wrapped (Fukaya) category.  A wrapped category is related to 
the symplectic homology in the same way as $\frak{Fukst}(X)$ is related 
to the quantum cohomology of $X$.\footnote{Namely the latter is the 
Hochshild cohomology of the former.}
Auroux proposed to generalize the notion of a wrapping 
to a partial wrapping. Ganatra-Pardon-Shande \cite{GPS1,GPS2,GPS3}
realized this program suceessfully.
\end{rem}

The homological mirror symmetry conjectures that a fully exceptional 
collection of a, say, Fano manifold $M$
is isomorphic to a version of $\mathscr{S}(M^{\vee},W)$.
Seidel first proved it in the case when $M=\C P^2$.
It then is generalized to the case of $2$-dimensional weighted projective space 
(Auroux-Kazarkov-Orlov \cite{AKO2}) and blow up of $\C P^2$ at 
$k$ point ($k\le 3$  by Ueda \cite{Ue} and then $k\le 8$ by 
Auroux-Kazarkov-Orlov \cite{AKO1}).
It is also proved when  $M$ is a toric Fano manifold (Abouzaid \cite{Ab1} 
and Fang-Liu-Treumann-Zaslov \cite{FLTZ}\footnote{The latter research uses 
$D$-module rather than Floer homology.}). 
\par
Compared to the story when $M$ is the symplectic side, 
an issue is that it is not clear how to obtain the function $W$ in a 
conceptional way.
The approach by \cite{AAK} (See also \cite{CLL}) may give 
a way to do so.
 
\section{Lagrangian correspondence and Gauge theory.}
\label{sec:COR+GAUGE}

Before discussing further the topological field theory 
of 2-3-4 dimensional manifolds based on Yang-Mills theory, 
we describe its `finite dimensional analogue', 
that is, a Lagrangian correspondence and its relation to 
Lagrangian Floer theory.
\par
Let $(X,\omega_X),(Y,\omega_Y)$ be symplectic manifolds.
We consider the product $X \times Y$ together 
with the symplectic form
$-\pi_1^*\omega_X + \pi_2^*\omega_Y$.
We write $-X \times Y: = (X \times Y,-\pi_1^*\omega_X + \pi_2^*\omega_Y)$.
Weinstein \cite{Wei} proposed to regard a Lagrangian 
submanifold of the product $-X \times Y$
as a morphism $X \to Y$ between symplectic manifolds.
We call a Lagrangian 
submanifold of $-X \times Y$ a Lagrangian correspondence.
\par
There are several reasons behind this proposal.
\begin{enumerate}
\item[(Wei1)]
We consider a symplectic diffeomorphism 
$\varphi: X \to Y$ (that is, a diffeomorphism $\varphi$
such that $\varphi^*\omega_Y = \omega_X$.)
Then its graph $\{(x,\varphi(x)) \mid x \in X\}$
is a Lagrangian submanifold of $-X \times Y$.
\item[(Wei2)]
Let $X_i$ be a symplectic manifold for $i=1,2,3$ and 
$L_{i(i+1)} \subset -X_i \times X_{i+1}$ be a 
Lagrangian submanifold for $i=1,2$.
We consider the fiber product:
$$
L_{13} = \{(x,y) \in L_{12} \times L_{23} \mid \pi_2(x) = \pi_1(y)\},
$$  
where for $z = (a,b) \in -X_i \times X_{i+1}$
we denote $\pi_1(z) = a$, $\pi_2(z) = b$.
In the generic case, that is, the case when
the fiber product is transversal, it is easy to see that 
$L_{13}$ is a smooth manifold and 
$L_{13} \to -X_1 \times X_3$, $(x,y) \mapsto (\pi_1(x),\pi_2(y))$
is a Lagrangian immersion.
Thus we can `generically' compose Lagrangian correspondences.
\item[(Wei3)]
Another example of a Lagrangian correspondence is obtained 
from the symplectic quotient.
Suppose $X$ is a symplectic manifold on which a compact Lie 
group $G$ acts preserving the symplectic structure.
We assume that there exists a moment map $\mu: X \to \frak g^*$,
where $\frak g^*$ is the Lie algebra of $G$.\footnote{
This means that for $V \in \frak g$ the vector field 
$V_*$ on $X$ obtained by $G$ action satisfies
$
\omega(V_*(p),W) = \langle D_p\mu(W),V\rangle
$. Here $W \in T_p(X)$ and $D_p\mu: T_p(X) \to \frak g^*$ is the 
derivative of $\mu$.
}
The symplectic quotient $X/\!/G$ is by definition $\mu^{-1}(0)/G$
(\cite{MW}). If $X/\!/G$ is smooth, it is known that it carries a symplectic form 
$\overline{\omega}$ whose pullback to  $\mu^{-1}(0)$ 
coincides with the restriction of the symplectic form of $X$.
\par
Now we consider
$$
L = \{(x,y) \in X \times X/\!/G \mid \mu(x) = 0, \,\, [x] = y\}.
$$
This is a Lagrangian submanifold of $-X \times X/\!/G$.
\end{enumerate}

Since this proposal by Weinstein looks so natural, 
there has been attempts to associate a functor 
$\frak F_{\frak L} : \frak{Fukst}(X) \to \frak{Fukst}(Y)$
to a Lagrangian correspondence $\frak L$.
A possible naive idea to do so is the following.
Let $L$ be a Lagrangian submanifold of $X$.
Instead of associating an object of  $\frak{Fukst}(Y)$ to $L$,
we try to define a right $\frak{Fukst}(Y)$-module $\frak F_{\frak L}(L)$.
In the cohomology level, $\frak F_{\frak L}(L)$ can be defined 
by associating the Floer homology
$HF(\frak L;L\times L')$ in the product $-X \times Y$ 
to a Lagrangian submanifold $L'$ of $Y$.
Actually we can construct an $A_{\infty}$ functor:
\begin{equation}\label{corr71}
\frak F_{\frak L}: \frak{Fukst}(X) \to \mathcal{RMOD}(\frak{Fukst}(Y))
\end{equation}
in this way. (See \cite[Section 5]{cor} for its rigorous construction.)
As we explained in Section \ref{sec:Ainfinity} of this article, an object of 
$\mathcal{RMOD}(\frak{Fukst}(Y))$ can be 
regarded as an `extended object' of $\frak{Fukst}(Y)$
(via Yoneda embedding). Thus (\ref{corr71}) 
could be regarded as a version of 
$\frak F_{\frak L} : \frak{Fukst}(X) \to \frak{Fukst}(Y)$.
However the problem is in this formulation it is difficult 
to compose $\frak F_{\frak L_{12}}$ and $\frak F_{\frak L_{23}}$
where $\frak L_{i(i+1)}$ is a Lagrangian submanifold of 
$-X_i \times X_{i+1}$.\footnote{This point is mentioned also in the first page of 
\cite{MWW}.}
This situation is somewhat similar to the following:
If we are given a  current $S$ (that is, a Schwartz kernel) on $M \times N$ and a smooth differential 
form $u$ on $M$ then we obtain a current $S_*(u)$  on $N$
by the equality
\begin{equation}\label{Shkernel}
S_*(u)(v) = S(u \times \pi_2^*(v)).
\end{equation}
In other words, an object of $\frak{Fukst}$ is an analogy of a 
smooth differential form and an object of $\mathcal{RMOD}(\frak{Fukst}(Y))$ 
is an analogy of a distributional form.
We use a `test Lagrangian' in place of a `test function'.
It is difficult to compose two operators of the form (\ref{Shkernel}).
In this sense, Theorem \ref{thm71} could be regarded as a kind of `regularity 
theorem'.
\par
In \cite{MWW,WW,WW2,WW3,WW4}, Wehrheim-Woodward-Ma'u 
used the following idea to go around this problem.
For a given symplectic manifold $X$, 
they consider a series of Lagrangian correspondences
$
L_{i} \subset -X_i \times X_{i+1} 
$
such that $X_0$ is a point and $X_n = X$.
They regard such a system $(L_0,\dots,L_n)$ as 
an object of expanded category $\frak{Fuk}(X)^+$.
Then, if $\frak L \subset -X \times Y$ is a 
Lagrangian correspondence,  
one can define
$$
(\mathcal{W}_{{\frak L}})_{\rm ob}: 
\frak{OB}(\frak{Fuk}(X)^+) \to \frak{OB}(\frak{Fuk}(Y)^+),
$$
by
$(L_0,\dots,L_n) \mapsto (L_0,\dots,L_n,\frak L)$.
\par
To define the category $\frak{Fuk}(X)^+$, one needs to define the 
Floer homology 
 between extended objects 
$(L_0,\dots,L_n)$,  $(L'_0,\dots,L'_{n'})$,
where $L_{i} \subset -X_i \times X_{i+1}$
and $L'_{i} \subset -X'_i \times X'_{i+1}$,
$X_0, X'_0$ are points and $X_n = X'_{n'} = X$.
They denote this Floer homology by $HF(L_0,\dots,L_n,L'_{n'},\dots,L'_0)$.
Wehrheim-Woodward-Ma'u used the notion of 
a pseudo-holomorphic quilt to define it.
Actually their definition is equivalent to the following:
\begin{equation}\label{formula72}
HF(L_0,\dots,L_n,L'_{n'},\dots,L'_0)
:= 
HF(L_0\times\dots\times L_n\times L'_0\times \dots\times L'_{n'};
\Delta).
\end{equation}
Here 
\begin{equation}\label{formula73}
\Delta \subset \left(\prod_{i=1}^{n-1} (-X_i\times X_i)\right) \times 
\left(\prod_{i=1}^{n'-1} (-X'_{i'}\times X'_{i'})\right) \times (-X \times X)
\end{equation}
is the product of diagonals. The right hand side of (\ref{formula72})
is the Floer homology of two Lagrangian submanifolds 
in the symplectic manifold given in (\ref{formula73}).
\par
In the simplest case, a pseudo-holomorphic quilt used to 
define $HF(L_0,L_1,L'_0)$ is a map $u = (u_1,u_2)$ from the 
domain depicted in Figure \ref{Figure6} below.
\begin{figure}[h]
\centering
\includegraphics[scale=0.3]{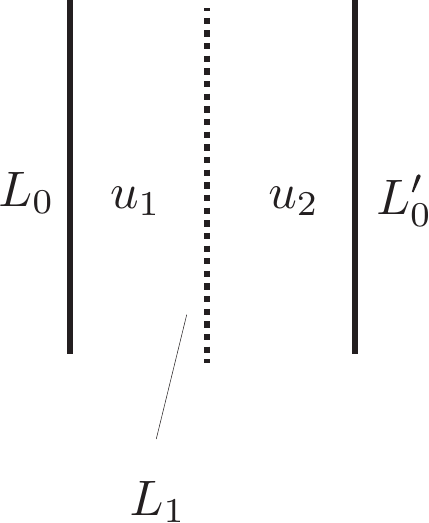}
\caption{A simplest pseudo-holomorphic quilt.}
\label{Figure6}
\end{figure}
Here the domain is divided into two parts. The map $u_1$ 
(resp. $u_2$) is defined on the 
left (resp. right) part  of the domain and is a pseudo-holomorphic map to $-X_1$
(resp. $X$).
The boundary conditions are required on the three lines, the left most, the right most 
and the middle (dotted) lines.
The boundary condition for the left most line is $u_1(z) \in L_0$,
one for the right most line is $u_2(z) \in L'_0$.
The boundary condition for the middle line is 
$(u_1(z),u_2(z)) \in L_1$. We remark that $L_1$ is a Lagrangian 
submanifold of $-X_1 \times X$. 
The middle line is called a seam.
\par
By reflexion principle, $(\overline u_1,u_2): [0,1] \to -X_1 \times X$ 
is a pseudo-holomorphic map which satisfies the boundary 
condition given by Lagrangain submanifolds $L_1$ and $L_0 \times L'_0$.
The moduli space of such pseudo-holomorphic maps is used to define
$HF(L_1,L_0 \times L'_0)$.
\par
The pseudo-holomorphic quilt used to define 
$HF(L_0,\dots,L_n,L'_{n'},\dots,L'_0)$
is as in Figure \ref{Figure7} below.
\begin{figure}[h]
\centering
\includegraphics[scale=0.3]{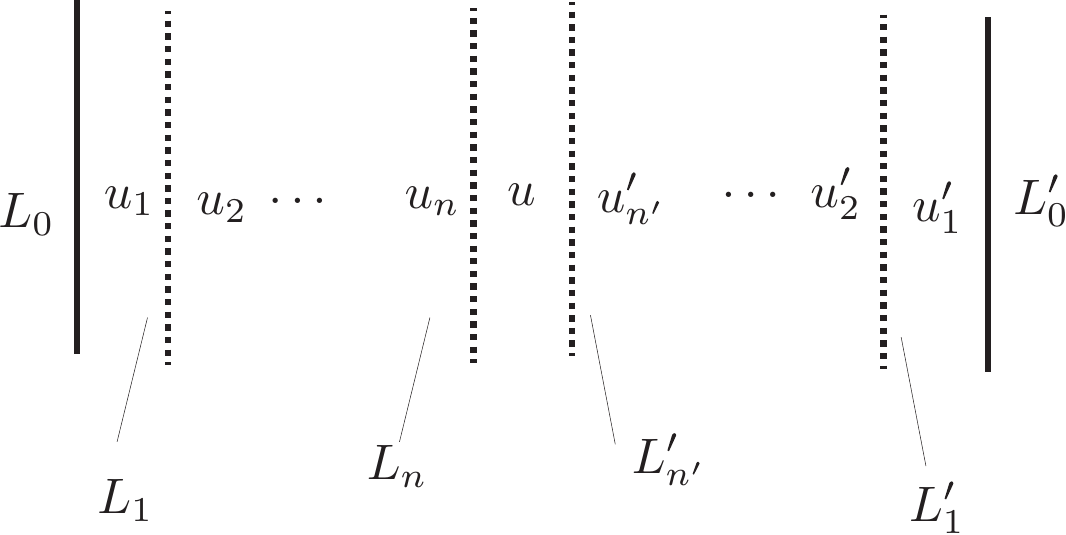}
\caption{A  pseudo-holomorphic quilt.}
\label{Figure7}
\end{figure}
Here $u_i$ (resp. $u'_i$) is a pseudo-holomorphic map to $-X_i$ (resp. $X'_i$)
and $u$ is a pseudo-holomorphic map to $X$.
\par
Wehrheim-Woodward-Ma'u studied the case when all the 
Lagrangian submanifolds $L_i$, $L'_i$ are monotone.
So we can use Oh's work \cite{Oh} to define it.
\par
In this way \cite{MWW} defined a version of an 
$A_{\infty}$ bi-functor\footnote{See \cite[Subsection 5.1]{cor} for 
the definition of $A_{\infty}$ bi-functor.}
\begin{equation}\label{MWW}
\mathcal{MWW}: \frak{Fukst}(-X\times Y) \times \frak{Fukst}(X)
\to \frak{Fukst}(Y)
\end{equation}
where all the Lagrangian submanifolds involved are assumed to 
be monotone and embedded and 
$\frak{Fukst}(\dots)$ is replaced by $\frak{Fukst}(\dots)^+$.
\par
\begin{figure}[h]
\centering
\includegraphics[scale=0.4]{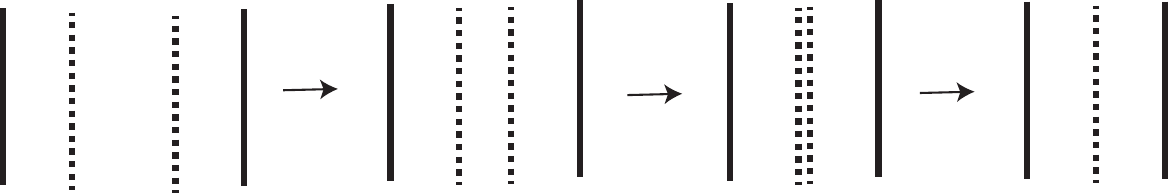}
\caption{Strip shrinking.}
\label{Figure8}
\end{figure}

The advantage to use peudo-holomorphic quilt rather than 
Floer homology in the direct product (as in (\ref{formula72}))
lies in the fact that, then, one can use `strip shrinking' to prove the 
next isomorphism
\begin{equation}\label{form74}
\aligned
&HF(L_0,\dots,L_n,L'_{n'},\dots,L'_0) \\
&\cong
HF(L_0,\dots,L_{n-1},L_n\times_XL'_{n'},L'_{n'-1},\dots,L'_0).
\endaligned
\end{equation}
The strip shrinking is a process to change the width between two seams
until it becomes $0$.
(See Figure \ref{Figure8}.)
Note that the method of using reflexion principle to replace 
Wehrheim-Woodward's definition by (\ref{formula72}) works only 
in the case when all the strips have the same width. 
Therefore it is not consistent with strip shrinking.
\par
\cite{WW,WW2,WW3,WW4} proved the isomorphism (\ref{form74})
under the assumption that all the Lagrangian submanifolds 
involved (including the fiber product $L_n\times_XL'_{n'}$) 
are embedded and monotone.\footnote{See \cite{BW} for an attempt 
to remove this condition.}
The isomorphism (\ref{form74}) is a version of 
composability of filtered $A_{\infty}$ functors 
associated to the composition of Lagrangian correspondences.
Later Lekili and Lipyanskiy \cite{LL} found an alternative 
method, using cobordism argument instead of strip shrinking.\footnote{
Using the virtual fundamental chain technique together 
with Lekili-Lipyanskiy's method we can prove (\ref{form74}) 
without assuming monotonicity or embedded-ness of Lagrangian 
submanifolds.}
For example Lekili-Lipyanskiy's proof of 
$HF(L_0,L_1,L'_1,L'_0) \cong HF(L_0,L_1 \times_X L'_1,L'_0)$ uses the moduli space 
of pseudo-holomorphic curves with the domain depicted by 
Figure \ref{Figure9} below.
\begin{figure}[h]
\centering
\includegraphics[scale=0.3]{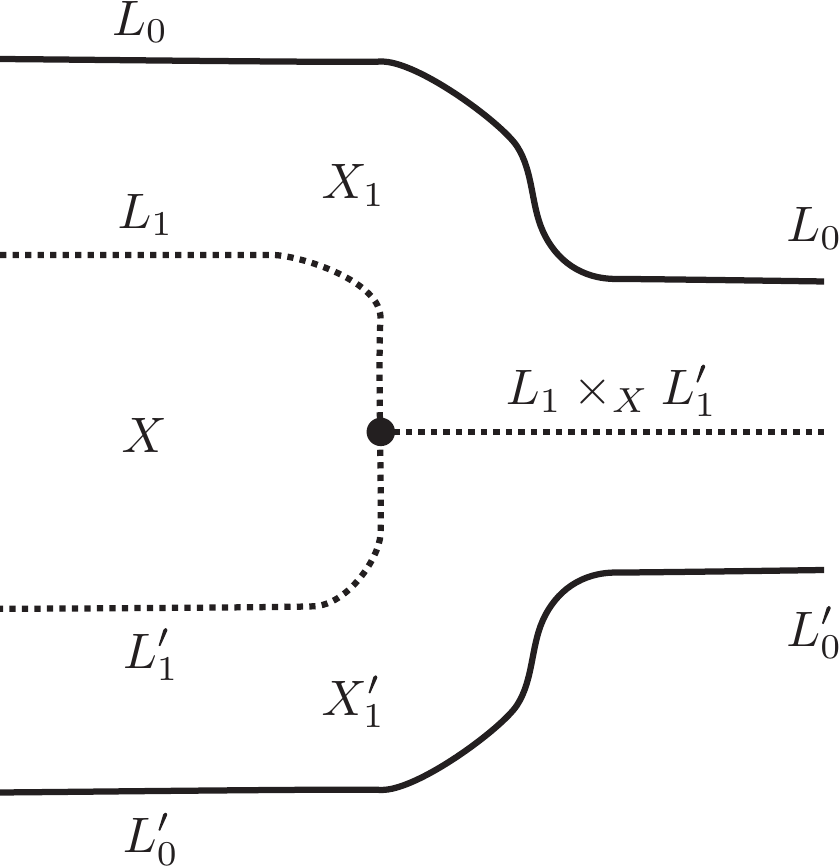}
\caption{Y Diagram.}
\label{Figure9}
\end{figure}
The three dotted lines (seams) are required to 
be mapped to $L_1$, $L'_1$, $L_1 \times_X L'_1$, respectively.
The other two curves are required to be mapped to 
$L_0$ and $L'_0$. 
The three domains are required to be mapped to 
$-X_1$, $X'_1$, $X$ as depicted. The maps on those domains
are required to be pseudo-holomorphic.
In the left end the pseudo-holomorphic quilt used to 
define $HF(L_0,L_1,L'_1,L'_0)$ appears. 
In the right end the pseudo-holomorphic quilt used to 
define $HF(L_0,L_1\times_XL'_1,L'_0)$ appears. 
Lekili and Lipyanskiy called Figure \ref{Figure9} 
the Y-diagram.
\par
After those works had been done, the author found that the 
naive idea in 1990's which we described at the beginning 
of this section, can be realized using the technology
developed in these 20 years and the cobordism argument 
due to Lekili and Lipyanskiy, as follows.
\par
We consider the filtered $A_{\infty}$ functor (\ref{corr71})
and also the Yoneda functor:
\begin{equation}
\frak{YON}: \frak{Fukst}(Y) \to
\mathcal{RMOD}(\frak{Fukst}(Y)).
\end{equation}
\begin{thm}\label{thm71}{\rm (\cite{cor})}
Let $(\frak L,\frak b) \in {\frak{OB}}(\frak{Fukst}(-X\times Y))$
and $(L,b) \in \frak{OB}(\frak{Fukst}(X))$
be objects such that the fiber product 
$L \times_X \frak L$ is transversal.
Then there exists a bounding cochain $b'$
of the immersed Lagrangian submanifold 
$L': = L  \times_X \frak L$ of $Y$ such that 
$(\frak{YON})_{\rm ob}(L',b')$ is homotopy equivalent to 
$(\frak F_{\frak L})_{\rm ob}(L,b)$.
\end{thm}
Theorem \ref{thm71} enables us to obtain a 
filtered $A_{\infty}$ functor 
$$
\mathcal{W}_{(\frak L,\frak b)} : \frak{Fukst}(X) \to \frak{Fukst}(Y)
$$
such that its composition with the Yoneda functor is 
homotopy equivalent to
the functor $\frak F_{\frak L}$ in (\ref{corr71}).
It is also proved in \cite{cor}
that $\mathcal{W}_{(\frak L,\frak b)}$ is induced from 
a filtered $A_{\infty}$ bi-functor
(\ref{MWW}), where objects of the 
categories involved are immersed and unobstructed.
(It is unnecessary to assume them to be monotone or embedded.)
\par
The association $(\frak L,\frak b) \mapsto \mathcal{W}_{(\frak L,\frak b)}$ is functorial.
Namely the following is proved in \cite{cor}.
Let $(\frak L_{i(i+1)},\frak b_{i(i+1)}) \in {\frak{OB}}(\frak{Fukst}(-X_i\times X_{i+1}))$
for $i=1,2$. We assume that the fiber product
$\frak L_{12} \times_{X_2} \frak L_{23}$ is transversal 
and put $\frak L_{13} = \frak L_{12} \times_{X_2} \frak L_{23}$, which is
an immersed Lagrangian submanifold of $-X_1\times X_{3}$.
Then there exists a bounding cochain $\frak b_{13}$ of $\frak L_{13}$
such that:
$$
\mathcal{W}_{(\frak L_{13},\frak b_{13})} \sim \mathcal{W}_{(\frak L_{23},\frak b_{23})} 
\circ \mathcal{W}_{(\frak L_{12},\frak b_{12})}
$$
where $\sim$ means homotopy equivalent.
\par
In  the situation of (Wei3) the Lagrangian 
correspondence from $X$ to $X /\!/G$
is expected to induce a functor from an
equivariant version of $\frak{Fukst}(X)$ 
to $\frak{Fukst}(X/\!/G)$.
Such a functor is studied by  Woodward and Xu in \cite{WX}.
It is an `open string analogue' of \cite{Wo},
based on gauged sigma model (See \cite{CGRS}).
\par\medskip
An infinite dimensional version of  the situation of (Wei3) 
appears in gauge theory as follows.
Let $\Sigma$ be a Riemann surface and $\mathcal E_{\Sigma}$
an $SO(3)$ or $SU(2)$ principal bundle on it.
We define $\mathcal A(\Sigma;\mathcal E_{\Sigma})$
to be the set of all $SO(3)$ or $SU(2)$ connections of 
$\mathcal E_{\Sigma}$.\footnote{The discussion here is 
formal or heuristic. So we do not specify how much 
regularity we require for an element of $\mathcal A(\Sigma;\mathcal E_{\Sigma})$.}
The space $\mathcal A(\Sigma;\mathcal E_{\Sigma})$ is an affine space 
and the tangent space of each point is the 
vector space $\Gamma(\Sigma,\Lambda^1\otimes ad \mathcal E_{\Sigma})$
of sections. Here $\Lambda^1$ is the bundle of 1-forms and 
$ad \mathcal E_{\Sigma}$ is the adjoint bundle associated to 
$\mathcal E_{\Sigma}$.
For $V,W \in \Gamma(\Sigma,\Lambda^1\otimes ad \mathcal E_{\Sigma})$
we put
$$
\Omega(V,W) =  -\frac{1}{8\pi^2}\int_{\Sigma} {\rm Tr}(V \wedge W),
$$
which defines a symplectic structure on $\mathcal A(\Sigma;\mathcal E_{\Sigma})$.
\par
Let $\mathcal G(\Sigma;\mathcal E_{\Sigma})$ be the gauge transformation group.
It is easy to see that the $\mathcal G(\Sigma;\mathcal E_{\Sigma})$ action 
on $\Gamma(\Sigma,\Lambda^1\otimes ad \mathcal E_{\Sigma})$ preserves the 
symplectic structure.
Moreover there is a moment map of this action, which is nothing but the curvature,
in the following way.
The Lie algebra of the gauge transformation group is $\Gamma(\Sigma,ad \mathcal E_{\Sigma})$.
We identify its dual with $\Gamma(\Sigma,\Lambda^2\otimes ad \mathcal E_{\Sigma})$
by
$$
V(W) = -\frac{1}{8\pi^2}\int_{\Sigma} {\rm Tr}(V \wedge W)
$$
where $W \in \Gamma(\Sigma,ad \mathcal E_{\Sigma})$ and 
$V \in \Gamma(\Sigma,\Lambda^2\otimes ad \mathcal E_{\Sigma})$ in the right hand side
and $V$ is a linear map $\Gamma(\Sigma,ad \mathcal E_{\Sigma}) \to \R$ in the left hand side.
Thus the moment map is a map
$
\mu: \mathcal A(\Sigma;\mathcal E_{\Sigma}) \to \Gamma(\Sigma,\Lambda^2\otimes 
ad \mathcal E_{\Sigma}).
$
One can check that 
$$
\mu(A) = -\frac{1}{8\pi^2} F_A.
$$
It implies that the symplectic quotient of the 
$\mathcal G(\Sigma;\mathcal E_{\Sigma})$ action 
on $\mathcal A(\Sigma;\mathcal E_{\Sigma})$
is the set of gauge equivalence classes of the flat connections.
\par
One can also observe that a pseudo-holomorphic map $u$
from a domain $D$ of $\C$ to $\mathcal A(\Sigma;\mathcal E_{\Sigma})$
is related to a solution of the ASD-equation on $D \times \Sigma$ as follows.
Let $z = s + \sqrt{-1} t$ be the standard coordinate of $D \subset \C$.
We write a connection of $\mathcal E_{D\times \Sigma}: = D \times \mathcal E_{\Sigma}$
as 
$$
\frak A = A(s,t) + \Phi(s,t) ds + \Psi(s,t) dt.
$$
Here $A(s,t)$ is a connection of $\mathcal E_{\Sigma}$ and 
$\Phi(s,t)$, $\Psi(s,t)$ are sections of 
$\Lambda^1_{\Sigma} \otimes ad \mathcal E_{\Sigma}$, for each $(s,t)$.
The ASD-equation (\ref{ASDeq}) for $\frak A$ can be written 
as
\begin{equation}\label{ASDprod}
\aligned
\frac{\partial A}{\partial t} + *_{\Sigma} \frac{\partial A}{\partial s}
&= d_{A(s,t)} \Psi(s,t) - *_{\Sigma} d_{A(s,t)} \Phi(s,t), \\
F_{A(s,t)} &= *_{\Sigma}\left(\frac{\partial\Phi}{\partial t} 
- \frac{\partial\Psi}{\partial s} - 
[\Phi,\Psi] \right).
\endaligned
\end{equation}
(See \cite{DS}.)
Here we use the product metric on $D \times \Sigma$.
Note that $*_{\Sigma}: \Gamma(\Sigma,\Lambda^1\otimes ad \mathcal E_{\Sigma})
\to \Gamma(\Sigma,\Lambda^1\otimes ad\mathcal E_{\Sigma})$
is a complex structure of $\mathcal A(\Sigma;\mathcal E_{\Sigma})$.
So if we define $\frak u: D^2 \to \mathcal A(\Sigma;\mathcal E_{\Sigma})$
by $\frak u(s,t) = A(s,t)$, then the first equation of (\ref{ASDprod})
can be regarded as 
\begin{equation}\label{holomorphicsss}
\frac{\partial \frak u}{\partial s} \equiv J \frac{\partial \frak u}{\partial t}
\mod \text{\rm Im} d_{A(s,t)} + *_{\Sigma}\text{\rm Im}d_{A(s,t)}.
\end{equation}
We consider the metric $g_D \oplus \epsilon g_{\Sigma}$.
Then the first equation does not change but the second equation becomes 
$$
F_{A(s,t)} = {\epsilon}*_{\Sigma}\left(\frac{\partial\Phi}{\partial t} 
- \frac{\partial\Psi}{\partial s} - [\Phi,\Psi] \right).
$$
Thus, in the limit $\epsilon \to 0$, this equation becomes 
$F_{A(s,t)} \equiv 0$, that is to say, $A(s,t)$ is flat.
Therefore $(s,t) \mapsto [\frak u(s,t)]$ defines a map 
$D \to R(\Sigma,\mathcal E_{\Sigma})$, which we write $u$.
Then, the equation (\ref{holomorphicsss}) says that 
$u$ is a (pseudo)holomorphic map.
Thus in the limit $\epsilon \to 0$ the equation (\ref{ASDprod}) 
becomes the (non-linear) Cauchy-Riemann equation of a map 
$u : D \to R(\Sigma,\mathcal E_{\Sigma})$.
This fact is the main motivation of Conjecture \ref{AFconj}.
\par
In this infinite dimensional version, 
the space corresponding to the Lagrangian correspondence 
is a subspace of 
$\mathcal A(\Sigma;\mathcal E_{\Sigma}) 
\times R(\Sigma,\mathcal E_{\Sigma})$
consisting of $(a,x)$ such that:
\begin{enumerate}
\item
[(mat1)] $a$ is a flat connection of $\mathcal E_{\Sigma}$.
\item
[(mat2)] 
$x$ is a point of $R(\Sigma,\mathcal E_{\Sigma})$ 
which is represented by the connection $a$. 
\end{enumerate}
This condition is called the matching condition and is introduced 
by Lipyanskiy \cite{Li}.\footnote{A related moduli space 
was introduced in \cite{gafa}. We remark that the line where 
the equation changes from the ASD-equation to the non-linear Cauchy-Riemann 
equation in \cite{gafa,Li} plays a similar role as the seams 
appearing in the pseudo-holomorphic quilt.}
The set of such $(a,x)$ is a Lagrangian submanifold of 
$\mathcal A(\Sigma;\mathcal E_{\Sigma}) 
\times R(\Sigma,\mathcal E_{\Sigma})$.
\par
Let us try to define a  moduli space similar to that of 
simple pseudo-holomorphic quilts depicted in Figure \ref{Figure6}, 
in our infinite dimensional situation.
The correspondence $L_1$ is the set of $(A,x)$ defined above.
$L'_0$ is a Lagrangian submanifold of $R(\Sigma,\mathcal E_{\Sigma})$.
$L_0$ is supposed to be a `Lagrangian submanifold' 
of the infinite dimensional symplectic manifold 
$\mathcal A(\Sigma;\mathcal E_{\Sigma})$.
Here however the idea is to regard a 
pair $(M,\mathcal E_M)$ of a 3-dimensional manifold $M$ 
with boundary $\Sigma$
and a bundle $\mathcal E_{M}$ such that 
its restriction to $\Sigma$ is $\mathcal E_{\Sigma}$
as a `Lagrangian submanifold' of 
$\mathcal A(\Sigma;\mathcal E_{\Sigma})$.
\par
More precisely, we consider the following 
moduli space.
We first put a Riemannian metric of $M$ such that 
the neighborhood of its boundary is 
isometric to (and is identified with)
$\Sigma \times [-1,0]$, where 
$\partial M = \Sigma \times \{0\}$.
We then consider a pair $(A,u)$
such that:
\begin{enumerate}
\item[(mix1)]
$A$ is a connection on $M \times \R$.
We require that $A$ satisfies the ASD-equation (\ref{ASDeq}).
\item[(mix2)]
$u : [0,1] \times \R \to R(\Sigma,\mathcal E_{\Sigma})$
is a map. We require that $u$ satisfies the non-linear 
Cauchy-Riemann equation (\ref{phceq}).
\item[(mix3)]
Let $\tau \in \R$. We consider the restriction 
$A(0,\tau)$ of $A$ to $\Sigma \times \{(0,\tau)\}$.
We then require that $(A(0,\tau),u(0,\tau))$  satisfies 
Conditions (mat1),(mat2).
\item[(mix4)]
We require that $u(1,\tau)$ is contained in the 
given Lagrangian submanifold $L$ of $R(\Sigma,\mathcal E_{\Sigma})$.
\end{enumerate}
\begin{figure}[h]
\centering
\includegraphics[scale=0.25]{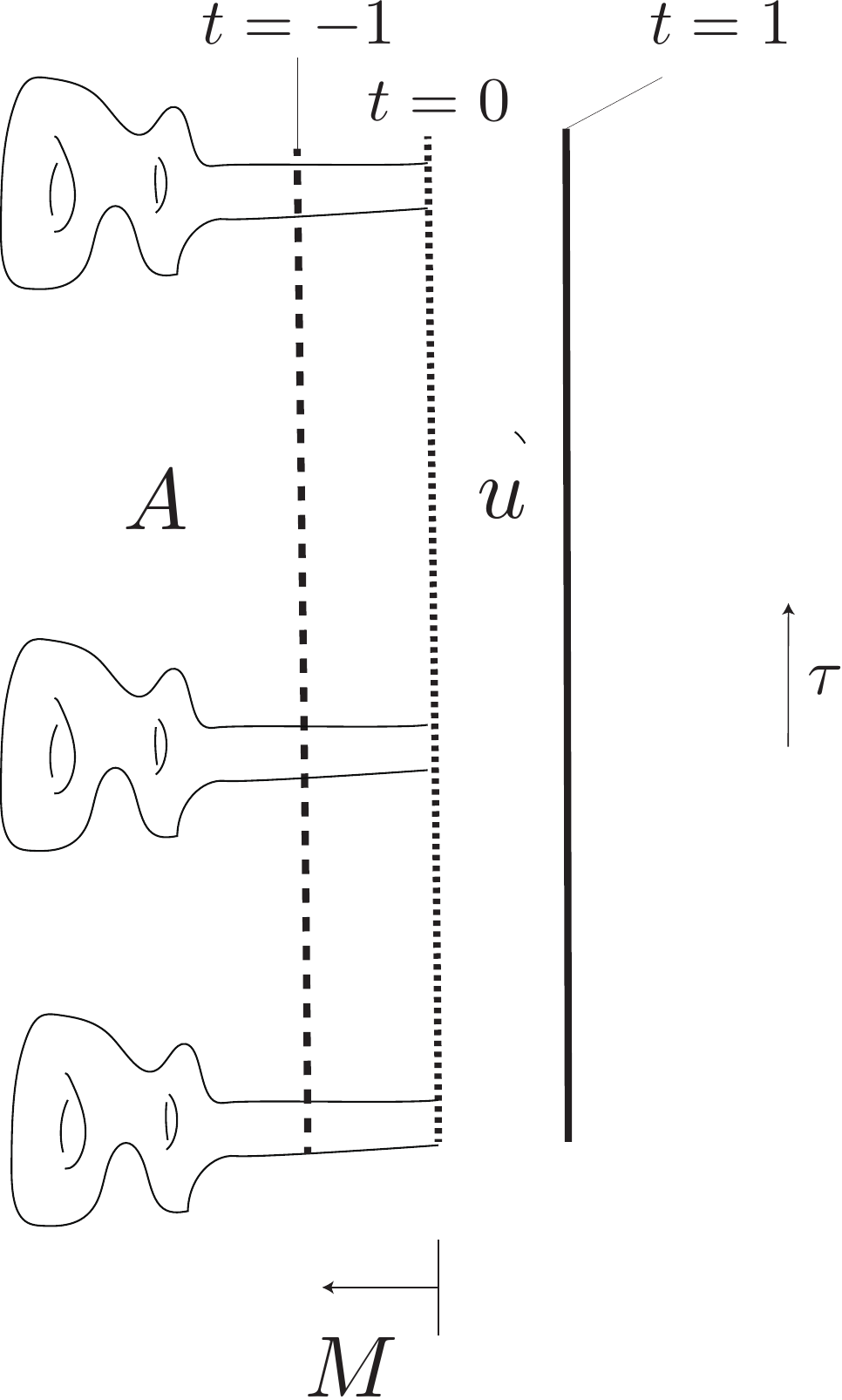}
\caption{Mixed equation.}
\label{Figure22}
\end{figure}

Using the moduli space of such pairs $(A,u)$ satisfying 
appropriate asymptotic boundary conditions, 
for $\tau \to \pm \infty$, we can define 
$HF(M;(L,b))$ the Floer homology with boundary condition $(L,b)$.
(Here $b$ is a bounding cochain of $L$.)
Note that we assume that $w_2(\mathcal E_{\Sigma}) = [\Sigma]$.
The basic analytic package to study such 
moduli spaces is established in \cite{Li,DFL1,DFL2}.
\par
Moreover using several Lagrangian submanifolds in place of 
a single $L$ in (mix4),
we can extend $(L, b) \mapsto HF(M;(L,b))$ to a
right filtered $A_{\infty}$ module over the 
filtered $A_{\infty}$ category $\frak{Fukst}(R(\Sigma,\mathcal E_{\Sigma}))$.
\par
Furthermore, the proof of Theorem \ref{thm71} can be generalized 
in this gauge theoretical situation, and we can show that 
there is a bounding cochain $b_M$ of the immersed 
Lagrangian submanifold $R(K;\mathcal E_{M})$ 
such that $(R(M;\mathcal E_{M}),b_M)$
via Yoneda embedding is sent to the right filtered $A_{\infty}$
module, $(L, b) \mapsto HF(M;(L,b))$.
\par
Finally by a cobordism argument as in 
\cite{fu1,LL,DFL2} we can 
show
\begin{equation}\label{SO3AF}
HF((R(M_1;\mathcal E_{M_1}),b_{M_1}),(R(M_2;\mathcal E_{M_2}),b_{M_2}))
\cong
I(M,\mathcal E_M)
\end{equation}
where $\partial M_1 = \Sigma = - \partial M_2$
and $M$ is obtained from $M_1$ and $M_2$ by gluing them along $\Sigma$.
\par
In the case when $R(M_i;\mathcal E_{M_i})$ is an embedded Lagrangian 
submanifold for $i=1,2$, (\ref{SO3AF}) is proved in \cite{DFL2}, where $b_{M_i} = 0$.
The general case is on the way being written.

\par\bigskip
{\small The author emphasizes that the reference below is far from being complete.
The references in the symplectic geometry 
are more emphasized than those in the complex or algebraic geometry.}

\bibliographystyle{amsalpha}

\end{document}